\documentclass[review,english,sn-mathphys]{sn-jnl}
\pdfoutput=1
\usepackage{amsmath}
\usepackage{amsfonts,amsthm}
\usepackage{amssymb}
\usepackage{graphicx}
\usepackage{placeins}
\usepackage[capitalize]{cleveref}
\usepackage{babel}
\usepackage{subfig}
\usepackage{upgreek,latexsym}
\usepackage{hhline}
\usepackage{epstopdf}
\usepackage[normalem]{ulem}

\newcommand{\innerp}[1]{\left\langle #1 \right\rangle}
\newcommand{\innerpp}[2]{\left\langle #1 \right\rangle _{#2}}
\newcommand{\bb}[1]{\mathbf{#1}}
\newcommand{\calbb}[1]{\boldsymbol{\mathcal{#1}}}

\newcommand{\abs}[1]{\left\vert{#1}\right\vert}

\newcommand{\ovl}[1]{\overline{#1}}
\newcommand{\norm}[2]{{\left\|#1\right\|}_{#2}}
\newcommand{\normsq}[2]{{\left\|#1\right\|^{2}}_{#2}}

\newtheorem{corollary}[]{Corollary}

\def\RR{{\mathbb{R}}}

\theoremstyle{plain}
\newtheorem{thm}{\protect\theoremname}[section]
\theoremstyle{plain}
\newtheorem{prop}[thm]{\protect\propositionname}
\theoremstyle{plain}
\newtheorem{lem}[thm]{\protect\lemmaname}
\theoremstyle{plain}
\theoremstyle{remark}
\newtheorem{rem}[thm]{\protect\remarkname}
\theoremstyle{definition}
\newtheorem{defn}[thm]{\protect\definitionname}

\makeatother

\providecommand{\definitionname}{Definition}
\providecommand{\theoremname}{Theorem}
\providecommand{\lemmaname}{Lemma}
\providecommand{\propositionname}{Proposition}
\providecommand{\remarkname}{Remark}

\begin{document}

\title[Identification of a boundary obstacle in a Stokes fluid]{Identification of a boundary obstacle in a Stokes fluid with Dirichlet--Navier boundary conditions: external measurements}

\author*[1]{\fnm{Louis} \sur{Breton}}\email{louis.breton@ciencias.unam.mx}

\author[2]{\fnm{Cristhian} \sur{Montoya}}\email{cristhian.montoya@unidu.hr}
\equalcont{These authors contributed equally to this work.}

\author[3]{\fnm{Pedro} \sur{Gonz\'alez-Casanova}}\email{casanova@matem.unam.mx}
\equalcont{These authors contributed equally to this work.}

\author[4]{\fnm{Jes\'us} \sur{L\'opez Estrada}}\email{jelpze@ciencias.unam.mx}
\equalcont{These authors contributed equally to this work.}

\affil*[1]{\orgdiv{Facultad de Ciencias}, \orgname{UNAM}, \orgaddress{\street{Av. Universidad 3000}, \city{Mexico}, \postcode{04510}, \country{Mexico}}}

\affil[2]{\orgdiv{Electrical Engineering and Computing}, \orgname{University of Dubrovnik}, \orgaddress{\street{\'Cira Cari\'ca 4}, \city{Dubrovnik}, \postcode{20000}, \country{Croatia}}}

\affil[3]{\orgdiv{Instituto de Matem\'aticas}, \orgname{UNAM}, \orgaddress{\street{Av. Universidad 3000}, \city{City}, \postcode{04510}, \state{State}, \country{Croatia}}}

\abstract{The problem of identifying an obstruction into a fluid
duct has several applications, one of them, for example in medicine the presence of
Stenosis in coronary vessels is a life threatening disease.
In this paper we formulate a continuous setting and study from a numerical perspective
the inverse problem of identifying an obstruction contained in a 2D
duct where a Stokes flow hits on the boundary
(Dirichlet and Navier--slip boundary conditions), generating
an acoustic wave. To be precise,
by using acoustic wave measurements
at certain points
on the exterior of the duct, we can identify
the location, extension and height of the obstruction. Thus, our framework
constitutes an external approach to solving this obstacle-inverse problem.  Synthetic examples are used to verify
the effectiveness of the proposed numerical formulation.}


\maketitle


\section{Introduction}
	Coronary artery disease (CAD), or coronary heart disease is the
	most common form of heart disease. In fact,  myocardial infarction (heart attack) 
	is one of the main causes of death in several countries of the first world, with a growing 
	incidence in developing countries. In the United States, at least 360,000 deaths 
	occur each year (see Lewandowoski and Cinquegrani, 
	{\it Coronary Heart Disease} in \cite[Sect.2 Chap.8]{wing2021cecil} and  
	\cite{Quarteroni2019} and references therein). 
	This disease is due to a very complex process of formation of atherosclerotic 
	plaques within the walls of the coronary arteries as a consequence of the gradual 
	accumulation of cholesterol, fatty acids, calcium, and fibrous connective tissue, 
	among others (see \cite[Sect.2 Chap.8]{wing2021cecil} and 
	\cite{Formaggia2009}).
	 
The presence of these atherosclerotic plaques leads to a local obstruction
	(called stenosis) of the blood flow with fatal consequences, such as
	myocardial infarction. Although the symptoms of CAD
	are pronounced in the later stages, it is very difficult to diagnose the disease before the first onset of symptoms, usually a sudden myocardial infarction with
	often causing death \cite{Rosamond2007}.	
		
   Coronary stenosis can be detected by fluoroscopy femoral catheterization or
   multi detector computed tomography scanning specially to quantify coronary calcium
   deposits that correlate with a significant obstructive lesion 
   \cite{wing2021cecil}. But both techniques are invasive.
   Therefore, the development of an alternative non-invasive technique of "auscultation" for
   the early unveiling of stenosis in coronary vessels is of great relevance.

	In addition, it is well established that the presence of an obstruction on 
	the wall of a coronary artery produces 
	an acoustic wave (murmur), which propagates from the wall of the coronary vessel 
    through thorax cavity reaching the chest surface, where it can be recorded with
    the help of special sensors (\cite{Lees1970}, \cite{Banks2009}). 
    What provides the fundamentals for a non-invasive technique of 
    ``auscultation'' to the early unveiling of stenosis into coronary vessels, 
    relatively simple to operate, inexpensive, and applicable to 
    hospital or office environments where background noise is inevitably present
    \cite{Lees1970}.
    
	From a physical perspective, the early detection problem of coronary 
	stenosis by non-invasive procedures leads to an inverse fluid-structure 
	problem with external data.

	Although the modelling and simulation of blood flow, direct problem, 
	through different techniques, have been carried out by many groups in order 
	to provide knowledge of flow behavior in order to clinical applications 
	\cite{2017-quarteroni}, \cite{Quarteroni2019}, there is 
	a growing interest in connecting the computational efforts and its theoretical 
	analysis in a unique way. However, due to the complexity arising from such 
	a combination, it is in an early stage even to the direct problem, i.e., 
	the well--posedness problem in the sense of Hadamard. 
	
	In this paper we present a 2D geometrical inverse problem about the 
	identificability and detection of stenosis in a coronary duct by using 
	``records'' from acoustic waves. 
	In our setting, the detection problem means position, extension, and coronary 
	lumen reduction. 
	In other words, we treat an inverse boundary obstacle problem through 
	``external measurements'' of local acoustic waves measurements, 
	which are located outside (at a suitable distance) from the coronary artery.  
	In this way, the words ``record'' and ``external measurements''
	turn out to be equivalent. 
	Now, this technique is not simple because of its non-invasiveness and the fact that it is used daily in 
	clinical practice through ultrasound techniques
	\cite{baron2009simulation,okita2011development} and acoustic wave radiation in biological materials 
	\cite{narasimhan2004high}. Additionally, considering complex and coupled models with patient--specific data makes a both 
	theoretical and computational efforts exorbitant. 	

	For these reasons, in this paper we consider a simple rectangular geometry $\Omega$
	in order to model the blood flow  as a viscous Stokes fluid in 2D with both 
	Dirichlet and Navier--slip conditions on the arterial walls. Henceforth, these 
	kind of boundary conditions are called mixed boundary conditions.
   	The domain has an obstacle $\mathcal{O}$ representing stenosis and such as 
    $\ovl{\mathcal{O}}$ intersects the fluid domain at only one part of 
	the boundary, i.e., $\ovl{\mathcal{O}}\cap \partial\Omega\neq\emptyset$. 
    The Dirichlet boundary conditions are applied to the velocity flow
	at the inlet and exit; meanwhile, Navier--slip boundary conditions are used for 
	the top and bottom walls.
	Simultaneously, the transmission and propagation of the sound through the 
	surrounding biological tissues is modeled with a linear wave equation with 
	state function $\omega$. The acoustic domain is covered by Cartesian grid whose 
	coupling zone with the fluid occurs on a wall (top or bottom) with Navier--slip 
	condition. 
	On the surface of the epidermis, a Higdon-absorbing boundary condition 
	\cite{1987Higdon} is applied (see Figure
	\ref{fig:Representacion onda w}). 
	Since Navier--slip conditions describe 
	the fluid on the arterial walls (top and bottom walls), the coupling term among 
	fluid and wave 
	corresponds to the normal component of the Cauchy tensor. Taking into account this information, our inverse problem consists in detecting 
	the obstacle $\mathcal{O}$ from the knowledge of 
	external data given by the wave--state $\omega$ on a subpart of its boundary 
	contactless with the fluid during a certain interval of time $(0,T)$. 
	In its current form, this problem and 
	its associated geometry might be considered as a toy problem in 2D, 
	however, there exist some relevant aspects to highlight: 
	\begin{enumerate}
	\item [a)] \textit{the condition 
	    $\ovl{\mathcal{O}}\cap \partial\Omega\neq\emptyset$.} 
		Roughly speaking, the above condition reflects an inverse boundary obstacle 
		problem for a linear two--dimensional fluid. 
		Although there is a wide literature on inverse obstacle problems involving 
		the obstruction inside the domain 
		( i.e., $\ovl{\mathcal{O}}\subset\subset \partial\Omega$) with different techniques
		 \cite{2019Caubet,1991-Allaire,alvarez2005identification,badra2011detecting,2010Caubet,caubet2016detection,2002Alessandrini}, as far as we know,  
		a first analysis for the case in which $\ovl{\mathcal{O}}\cap \partial\Omega\neq\emptyset$
		has not been reported.  
 
	\item [b)] \textit{boundary conditions and transmission term.} Following 
	   a recently theoretical setting for Navier--Stokes fluids with Navier--slip boundary conditions
	   \cite{amrouche2014lp,acevedo2019stokes}, we prove the existence and uniqueness of a weak solution 
	   for Stokes fluids by considering mixed boundary conditions; see \eqref{eq:Stokes-navier-slip} below. The transmission term follows
	   the ideas formulated in \cite{Lions1969}, and thus the information given by the Cauchy tensor is located on a boundary 
	   condition for a wave equation.
	   	
	\item [c)] \textit{external measurements to the inverse problem}. 
	   Our approach addresses the boundary obstacle fluid--structure (rigid structure) inverse problem from exterior measurements given
	   by local data from sound wave. We note that a computational study of a direct problem for a different flow--acoustic model was
	   done in  \cite{seo2012coupled}. Nevertheless, our model has a different structure, which relies on Navier-slip boundary conditions, the acoustic wave coupling with the fluid through the Cauchy tensor.  
	   Regarding external observations, the recent article \cite{2020Karageorghis} uses the method of fundamental 
	   solutions with Tikhonov's regularization to solve the identification of obstacles immersed in a steady Oseen fluid (with Dirichlet boundary 
	   conditions), and whose exterior measurements are given by the velocity fluid, the stress force or the pressure gradient.

\end{enumerate}
	
This article is organized as follows. In Section 2, we introduce the mathematical formulation of the direct problem, a fluid-sound wave system. In Section 3, a brief discussion on the problem of identification of the detection of obstruction in a viscous fluid flow is given. In Section 4, we present the RBF free divergence hybrid method and numerical simulations of the direct problem. In Section 5, we numerically solve the inverse geometric problem of identifying fluid obstruction from sound wave measurement data. Section 6 is devoted to the final remarks and conclusions.

\section{The problem setting}
\begin{figure}[htbp]
	\begin{centering}
		\subfloat[]{\includegraphics*[bb = 50 200 2400 800,scale=0.15]{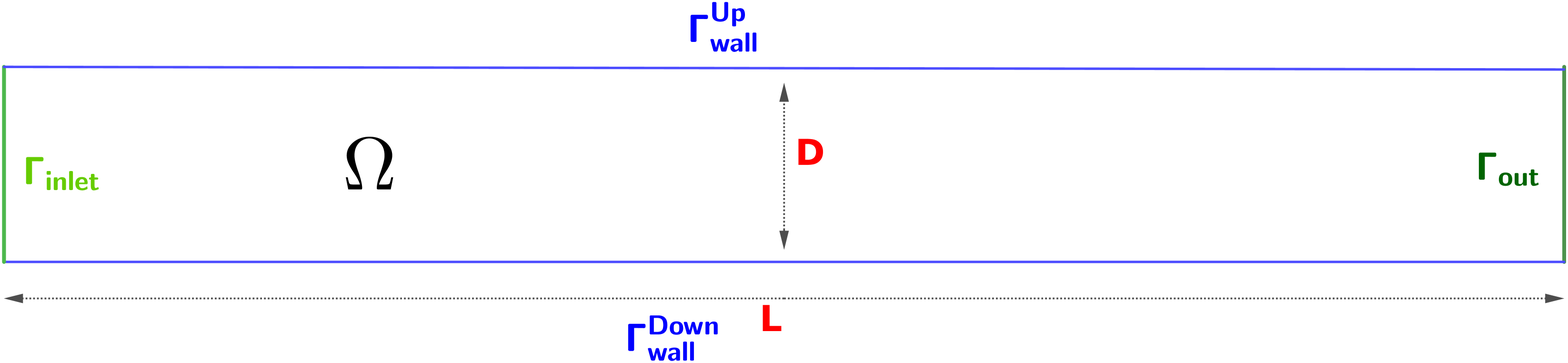}}\\
		\subfloat[]{\includegraphics*[bb = 50 500 2400 1000,scale=0.15]{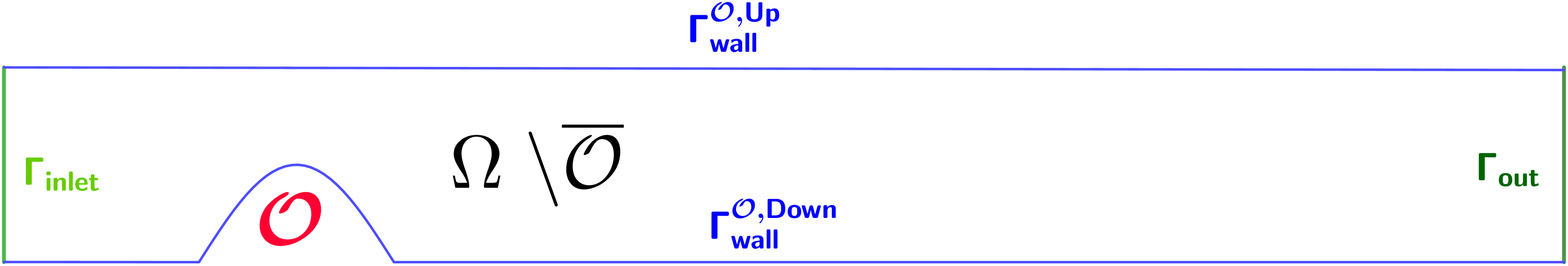}}  
	\par\end{centering}
	\centering{}\caption{Geometric representation of the domain: a) $\Omega$ , b) $\Omega_{{\mathcal {O}}}:=\Omega\backslash \ovl{{\mathcal {O}}}$.}
	\label{fig:Representacion-geomtrica-Dominio}
	\end{figure}

	We define the spatial domain of the
	fluid flow by $\Omega:=(0,L)\times(0,D)$,  where $L>0$ can be understood as the length of the blood vessel
	and $D>0$ its diameter. Assume that $\partial\Omega$ is divided into four parts $\Gamma_{inlet}, \Gamma_{wall}^{up}, 
	\Gamma_{wall}^{down}$ and $\Gamma_{out}$ of a non-vanishing measure such that 
	$\Gamma_{inlet}\bigcap\Gamma_{wall}^{up}\bigcap\Gamma_{wall}^{down}\bigcap\Gamma_{out}=\emptyset$. Throughout
	the paper, we adopt the convention that a boldface character denotes a vector or a tensor.  
	Furthermore, we consider an obstacle 
	${\mathcal{O}}\subset\ovl{\Omega}$ with non-empty interior such that  
	$\Omega_{{\mathcal {O}}}:=\Omega\backslash \ovl{{\mathcal{O}}}$ only has one 
	connected component, which is split into four parts 
	$\Gamma_{inlet}, \Gamma_{wall}^{\mathcal{O},up}, 	\Gamma_{wall}^{\mathcal{O},down}$, 
	$\Gamma_{out}$ and the boundary $\partial\mathcal {O}$ 
	of $\mathcal {O}$ satisfies (see Figure \ref{fig:Representacion-geomtrica-Dominio}):
	\begin{equation}\label{cond.borde}
	\begin{cases}
		\partial{\mathcal{O}}\cap\Gamma_{inlet} &=\emptyset,\\
		\partial{\mathcal{O}}\cap\Gamma_{out}   &=\emptyset,\\
		\partial{\mathcal{O}}\cap\Gamma_{wall}^{up} &\neq\emptyset,
	\end{cases}\quad \mbox{or}\quad 
	\begin{cases}
		\partial{\mathcal{O}}\cap\Gamma_{inlet} &=\emptyset,\\
		\partial{\mathcal{O}}\cap\Gamma_{out} &=\emptyset,\\
	\partial{\mathcal{O}}\cap\Gamma_{wall}^{down} &\neq\emptyset.
	\end{cases}
	\end{equation}

	The velocity vector $\bb{u}$ and the scalar pressure $p$ of the fluid in the presence of the obstacle 
	${\mathcal {O}}$ are modeled by the Stokes system with mixed boundary conditions:

	\begin{equation}\label{eq:Stokes-navier-slip}
	\begin{cases}
	 	\bb{u}_{t}-div(\sigma(\bb{u},p))=0 &\mbox{in }  \Omega_{{\mathcal {O}}}\times(0,T),\\
	 	
		div(\bb{u})=0 &\mbox{in } \Omega_{{\mathcal {O}}}\times(0,T),\\
		
	\bb{u}=\bb{g_{in}} &\mbox{on } \Gamma_{inlet}\times(0,T),\\
	
		\bb{u}=\bb{g_{out}} &\mbox{on } \Gamma_{out}\times(0,T),\\
		
		\bb{u}\cdot\bb{n}=0,\, \left[\sigma(\bb{u},p)\bb{n}\right]_{tg}=0 &\mbox{on }  \Gamma_{wall}^{\mathcal{O}}\times(0,T),\\
		
		\bb{u}(\cdot,0)=\bb{u}_{0}(\cdot) &\mbox{in } \Omega_{{\mathcal {O}}},
	\end{cases}
	\end{equation}
	where $\bb{n}$ is the unit outward normal vector of $\partial\Omega_{{\mathcal {O}}}$, $\bb{g_{in}} $ and $\bb{g_{out}}$ 
	are given nontrivial Dirichlet data. Besides,
	$\sigma(\bb{u},p)$ is the Cauchy stress tensor, i.e., 
	$\sigma(\bb{u},p)=2\bb{D}(\bb{u})-Ip=2\mu\frac{\nabla+\nabla^{t}}{2}\bb{u}-Ip$, 
	$I$ is the identity matrix of order $2\times 2$ and
	$\mu>0$ is the kinetic viscosity coefficient. The subscript $\textit{tg}$ denotes the tangential component of the corresponding vector field,
	which can be defined as (see \cite{1823naviermemoire}) $\bb{v}_{tg}:=\bb{v}-(\bb{v}\cdot\bb{n})\bb{n}$.

\begin{figure}[htbp]
	\begin{centering}
		\includegraphics[scale=0.25]{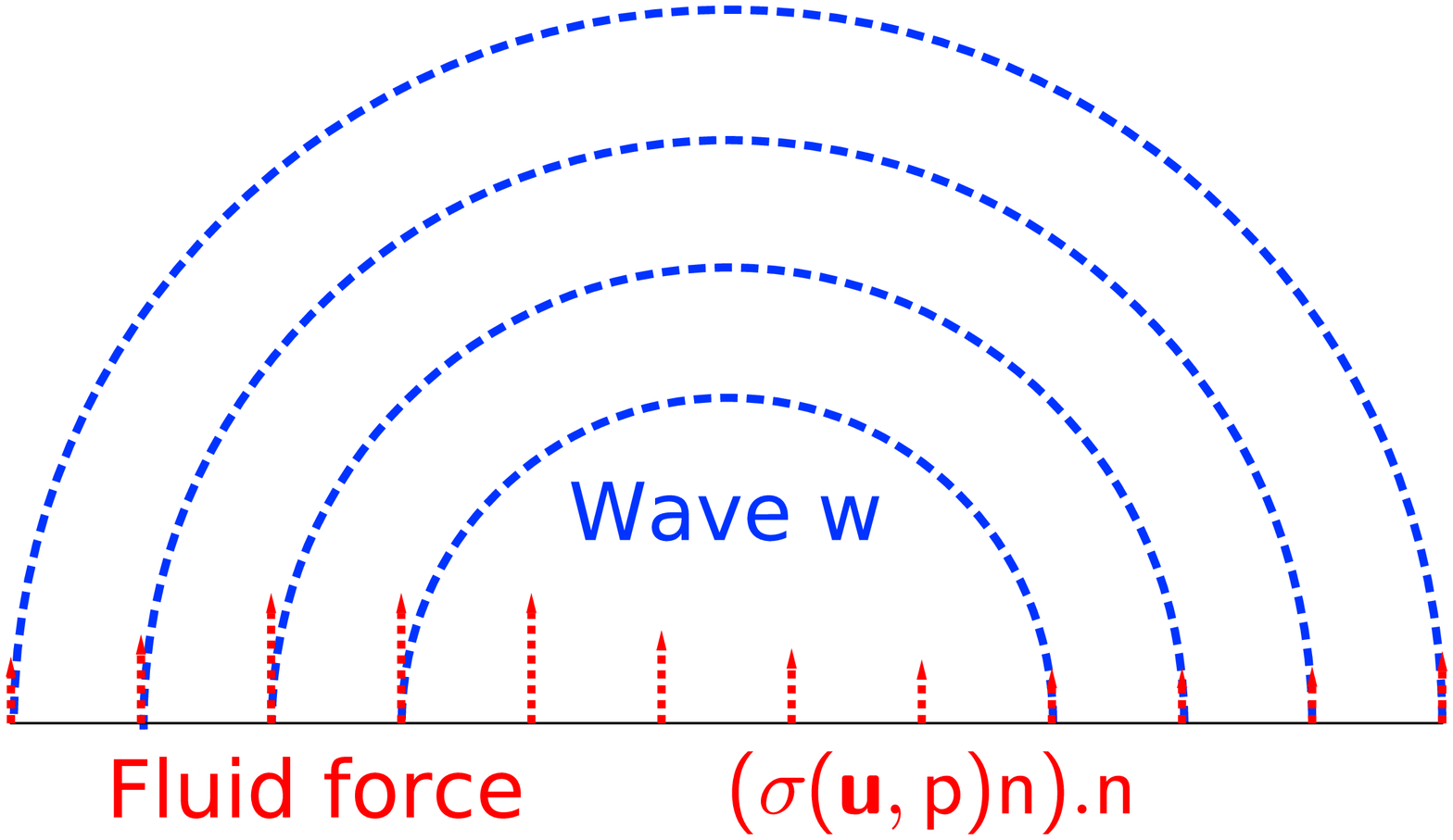}
		\includegraphics*[bb = 0 325 2000 1500,scale=0.1]{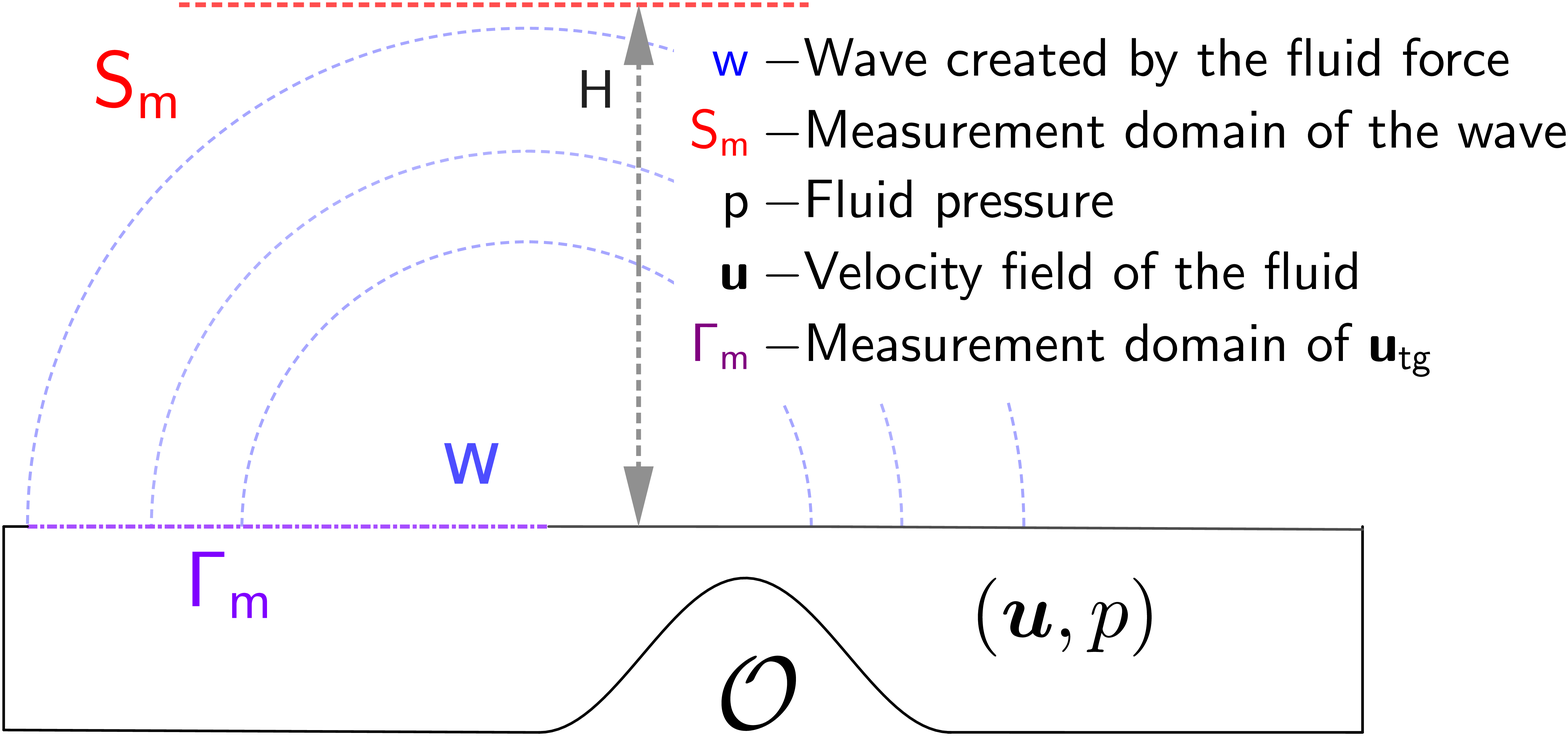} 
		\par\end{centering}
		\caption{\label{fig:Representacion onda w}Geometric representation of the
		wave $w$.}
	\end{figure}

	On the other hand, sound propagation through the arterial wall and surrounding tissue is modeled
	via a linear wave equation in a finite domain $S:=(0,L)\times (D,H)$  as follows:
	
	\begin{equation}\label{eq:Onda}
\begin{cases}
w_{tt}-c^{2}\Delta{w}=0 & \mbox{in }S\times(0,T),\\
w=(\sigma(\bb u,p)\bb{\nu})\cdot\bb{\nu} & \mbox{on }\left([0,L]\mathbb{\times}\{D\}\right)\times(0,T),\\
\frac{\partial w}{\partial t}+c\frac{\partial w}{\partial n}=0 & \mbox{on }\left(\partial S\backslash\left([0,L]\mathbb{\times}\{D\}\right)\right)\times(0,T),\\
w_{t}(\cdot,0)=w(\cdot,0)=0 & \mbox{in }S,
\end{cases}
	\end{equation}
	where $\bb{\nu}:=(0,1)^{t}$, $\frac{\partial}{\partial n}$ represents the normal partial derivative operator and $c$
	is the given wave speed. From a mechanical point of view, 
	$H>0$ can be viewed as the distance between the blood vessel and the epidermis, where the acoustic wave 
	will be observed, as shown in Figure \ref{fig:Representacion onda w}, while the term $(\sigma(\bb{u},p)\bb{\nu})\cdot\bb{\nu}$
	physically corresponds to the normal stress exerted by the fluid on the upper wall of the domain $\Omega_{{\mathcal {O}}}$. 
	We also emphasize that the value $(\sigma(\bb{u},p)\bb{\nu})\cdot\bb{\nu}$ 
	constitutes the coupling term in an explicit form, that is, once we have determined
	the solution to \eqref{eq:Stokes-navier-slip}, such a boundary condition of the wave equation \eqref{eq:Onda} is 
	explicitly obtained. 
	Note that the second boundary condition corresponds to the Higton absorbing boundary condition 
	of order 1,  \cite{givoli2006finite,1987Higdon}.  It is well known that the imposition of absorbing boundary conditions is a technique 
	used to reduce the necessary spatial domain when numerically solving partial differential equations admit traveling waves. The well--posedness
	of the initial boundary value problem of the absorbing boundary conditions, coupled to the wave equation have been discussed in 
	several papers, for instance, \cite{1986Trefethen, 1980Bayliss, 1986Kosloff}.

	As mentioned, there are no papers dealing with the effective reconstruction of an obstacle in contact with a subset
	of $\partial\Omega$, and using some type of external measurement for unsteady Stokes fluids with mixed 
	boundary conditions. The aim of this work consists in determining an obstacle  $\mathcal{O}\subset\ovl{\Omega}$ 
	(time--independent and satisfying \eqref{cond.borde})  using external measurements from  
	the acoustic wave $w$ (see \eqref{eq:Onda}) in the set $S_{m}\times (0,T_{0})$, where 
	$S_{m}:=[k_1,k_2]\times \{H\} \subset \partial S$ with $0<k_1<k_2<L$, and $T_{0}<T$. 
	Our idea consists in uncoupling the system \eqref{eq:Stokes-navier-slip}--\eqref{eq:Onda} and analysing the
	inverse obstacle problem throughout two inverse sub-problems, namely: 
	one inverse problem for the wave equation related to recover a boundary datum from $w\rvert_{S_{m}\times(0,T_{0})}$ and,
	another one connecting the recovered boundary datum for the wave equation with partial information of the  
	Cauchy tensor given by the Stokes fluid; see Theorem \ref{identification_theo}. 
	To consider a wave domain $S$ independent of the obstruction $\mathcal {O}$, we impose the geometrical condition
	$\partial{\mathcal {O}}\bigcap \ovl{\Gamma_{wall}^{up}}=\emptyset$, as well as a priori knowledge of one of the 
	configurations given in \eqref{cond.borde}. 
\section{Theoretical results for the direct problem}
	In this section, we give a detailed uniqueness and existence proof of weak solutions for the 
	two--dimensional Stokes system with mixed boundary conditions, i.e. the Stokes system with 
	Dirichlet and Navier--slip boundary conditions.  Although these mixed boundary conditions were mentioned 
	in the paper \cite{2016AbdaKhayat}, there is no rigorous proof in that article. For that reason and for 
	the sake of self--containment, we use as a starting point the recent work \cite{acevedo2019stokes,2016Bada,amrouche2014lp}  
	associated with both the unsteady and steady Stokes system with Navier boundary conditions. In the literature, these types of problems are commonly referred to as Zaremba's problem. 
	For simplicity, we set the viscosity at $\mu=1$.
	
	Until further notice, we assume that $\Omega\subset\mathbb{R}^{2}$
is a bounded curvilinear polygone of class $C^{1,1}$ in the sense
of \cite[Definition 1.4.5.1]{grisvard2011elliptic}. We assume that
there is $M\in\mathbb{N}$ such that:
\[
\Gamma=\partial\Omega=\bigcup_{j=1}^{2M}\overline{\Gamma_{j}}\quad,
\]
were each $\Gamma_{j}$ is a curve of class $C^{k,1}$, $\Gamma_{j+1}$
follows $\Gamma_{j}$ according to the positive orientation and we
have that $\Gamma_{j}\subset\Gamma$ , $\Gamma_{j}\bigcap\Gamma_{j+1}=\emptyset$. We define $\Gamma_{D},\Gamma_{N}\subset\Gamma$ as:
\[
\Gamma_{D}=\bigcup_{j=2}^{M}\Gamma_{(2j)}\text{ and }\Gamma_{N}=\bigcup_{j=1}^{M}\Gamma_{(2j-1)}\quad.
\]

\noindent
 As mentioned before, we denote by $\bb{v}_{tg}$ as 
	the tangential component of $\bb{v}$, i.e., $\bb{v}_{tg}=\bb{v}-(\bb{v}\cdot\bb{n})\bb{n}$.	Moreover, we consider the following spaces equipper with theirs usual respectives norms (for instance see \cite{amrouche2014lp}  ) :
	\begin{align*}
		\bb{V}(\Omega) & =\{\bb{u}\in\bb{H}^{1}(\Omega):\ div(\bb{u})=0,\bb{u}\lvert_{\Gamma_{D}}=0,\bb{u}\cdot \bb{n}=0\ on\ \partial \Omega \},\\
		\bb{H}_{0,\Gamma_{D}}(\Omega) & =\{\bb{u}\in\bb{H^{1}}(\Omega):\ \bb{u}\lvert_{\Gamma_{D}}=0\},\\
		\bb{H}(div,\Omega) & =\{\bb{u}\in\bb{L}^{2}(\Omega):\ div(\bb{u})\in\bb{L}^{2}(\Omega)\},\\
		\bb{H}_{0}(div,\Omega) & =\{\bb{u}\in\bb{L}^{2}(\Omega):\ div(\bb{u})\in\bb{L}^{2}(\Omega),\bb{u}\cdot \bb{n}=0\ on\ \partial \Omega \},\\
		\bb{E}(\Omega) & =\{\bb{u}\in\bb{H}^{1}(\Omega):\ \Delta\bb{u}\in\bb{H}_{0}(div,\Omega)'\}.
	\end{align*}

The main novelty of this this section with respected to \cite{amrouche2014lp}, is our mixed boundary condition. Therefore to treat this rigorously we have to make use of the following Lions-Magenes space (see \cite[Chp 11]{lions2012non}, \cite[Section 1.5.2]{grisvard2011elliptic} for futher details):

\begin{equation*}
 \bb{H}_{00}^{1/2}(\Gamma_i)=\left\{ \bb{u}\in \bb{H}^{1/2}(\Gamma_i): \ d(x,\partial\Gamma_i)^{-1/2}\bb{u}\in \bb{L}^{2}(\Gamma_i)\right\} \ , 
\end{equation*}
\noindent
equipped with the norm: $\normsq{\bb u}{\bb H_{00}^{1/2}(\Gamma_{i})}=\normsq{\bb u}{\bb{H^{1/2}}(\Gamma_{i})}+\normsq{d^{-1/2} \bb u}{\bb{L}^{2}(\Gamma_{i})}$, where $d(x,\partial\Gamma_i)$ is the geodesic distance from $x$ to $\partial \Gamma_{i}$.
Following this definition we set:
\[
\bb{H}_{00}^{1/2}(\Gamma_{D})=\prod_{j=2}^{M}\bb{H}_{00}^{1/2}(\Gamma_{(2j)})\text{ and }\bb{H}_{00}^{1/2}(\Gamma_{N})=\prod_{j=1}^{M}\bb{H}_{00}^{1/2}(\Gamma_{(2j-1)})
\]
In other to avoid heavy notation, for this section, we will denote:
\begin{align*}
& \innerpp{\cdot,\cdot}{\Gamma_i} \text{ as the dual pairing between $ \bb{H}_{00}^{1/2}(\Gamma_i)'$ and $\bb{H}_{00}^{1/2}(\Gamma_i)$} \ , \\ 
& \innerpp{\cdot,\cdot}{\partial \Omega}  \text{ as the dual pairing between $ \bb{H}^{-1/2}(\partial \Omega)$ and $\bb{H}^{1/2}(\partial \Omega)$}\ , \\ 
& \innerpp{\cdot,\cdot}{\Omega}  \text{ as the dual pairing between $\bb{H}_{0}(div,\Omega)' $ and $\bb{H}_{0}(div,\Omega)$.} 
\end{align*}

Too proceed in the followings subsections we shall stated the following Green identity which we will be proven in the appendix \ref{Appendix.Green_identity}.

\begin{thm}\label{Th:Green_identity}
Let $\Omega$ be a bounded open subset of $\mathbb{R}^{2}$ whose
boundary is a curvilinear polygon of class at least $C^{1}$. Let
$\gamma_{j}(\bb{u})$ denote the following trace mapping:
\[
\gamma_{j}(\bb{u}):\bb{u}\rightarrow \bb{u}\lvert_{\Gamma_{j}}\quad .
\]
Furthermore assume that $\gamma_{j}(\bb u)\in \bb{H}_{00}^{1/2}(\Gamma_{j})$
for all $j\in\{1,..,2M\}$. Then we have the following Green identity:

\begin{equation}
-\left\langle \Delta\bb u,\bb{\varphi}\right\rangle _{\Omega}=\int_{\Omega}2\bb D(\bb u):\bb D(\bb{\varphi})\,dx-\sum_{j=1}^{2M}{\left\langle 2[\bb D(\bb u) \bb{n}]_{tg},\gamma_{j}(\bb u)\right\rangle _{\Gamma_{j}}}\label{eq:Green_identity}
\end{equation}

for all $\bb u\in\bb E(\Omega)$ and $\varphi\in\{\bb{u} \in \bb H^{1}(\Omega) \lvert div(\bb u)=0\}$
\end{thm}

\begin{rem}
In Theorem \ref{Th:Green_identity} the spaces $\bb{H}_{00}^{1/2}$ are necessary to decompose the dual dual pairing $\innerpp{\cdot,\cdot}{\partial \Omega}$ into the sum of the dual pairing. This is because we lack regularity of the tangential component of the Cauchy stress tensor $[\bb D(\bb u) \bb{n}]_{tg}$ when $\bb u\in \bb{E}(\Omega)$.
\end{rem}

\subsection{Steady Stokes system}
Although we are only interested in proving the existence and uniqueness of the evolutionary Stokes flow with mixed boundary conditions, it is important that we first prove the stationary case. Thus, the aim of this subsection is to prove the existence and uniqueness of a weak solution for the steady Stokes system with 
	mixed boundary conditions:
	\begin{equation}\label{eq:Stokes_stacionary}
	\begin{cases}
		-div(\sigma(\bb{u},p))=\bb{f} &\mbox{in }  \Omega,\\
		div(\bb{u})=0 &\mbox{in }  \Omega,\\
		\bb{u}\cdot\bb{n}=0 &\mbox{in }  \partial\Omega,\\
		\bb{u}=\bb{g} &\mbox{on }  \Gamma_{D},\\{}
		[2\bb{D}(\bb{u})\bb{n}]_{tg}=\bb{h} &\mbox{on }  \Gamma_{N}.
		\end{cases}
	\end{equation}
	The main result of this section is given in the following theorem. 
	\begin{thm}\label{th.steadycase}	
		Let $\bb{f}\in\bb{H}_{0}(div,\Omega)'$, $\bb{h}\in (\bb{H}^{1/2}_{00}(\Gamma_{N}))'$
		and $\bb{g}\in \bb{H}^{1/2}_{00}(\Gamma_{D}) $ such that:

\begin{equation}\label{cond:compatiblity_g}
\bb{g}\cdot \bb{n}\lvert_{\Gamma_{D}}=0
\text{ and }
\bb{h}\cdot \bb{n}\lvert_{\Gamma_{N}}=0 \ .
\end{equation}

		Then the problem of finding $(\bb{u},p)\in\bb{H}^{1}(\Omega)\times L_{0}^{2}(\Omega)$
		satisfying \eqref{eq:Stokes_stacionary} in the distribution
		sense has a unique solution. Furthermore, we have the following estimate:
		\begin{equation}
			\norm{\bb{u}}{\bb{H}^{1}(\Omega)}+\norm{p}{L^{2}(\Omega)}\leq C\left(\norm{\bb{f}}{\bb{H}_{0}(div,\Omega)'}
			+\norm{\bb{g}}{\bb{H}^{1/2}_{00}(\Gamma_{D})}+\norm{\bb{h}}{(\bb{H}^{1/2}_{00}(\Gamma_{N}))'}\right). \label{eq:Estimation_stationary-1}
		\end{equation}
	\end{thm}
	Before presenting the proof of Theorem \ref{th.steadycase}, some preliminaries results are needed. 
	\begin{prop}\label{Prop:Steady_equivalence}
		Assume $\bb{g}=0$ in \eqref{eq:Stokes_stacionary}. Let $\bb{f}\in\bb{H}_{0}(div,\Omega)'$,
		$\bb{h}\in{(\bb{H}^{1/2}_{00}(\Gamma_{N}))'}$ such that:

\begin{equation}\label{cond:compatiblity_h}
\bb{h}\cdot \bb{n}\lvert_{\Gamma_{N}}=0 \ .
\end{equation}		
		
		Then, the problem of finding a pair $(\bb{u},p)\in\bb{H}^{1}(\Omega)\times L^{2}(\Omega)$
		satisfying \eqref{eq:Stokes_stacionary} in the distribution sense is equivalent to:
		\begin{equation}
		\begin{cases}
			\text{Find \ensuremath{\bb{u}\in\bb{V}(\Omega)} such that}\\
			\forall\bb{v}\in\bb{V}(\Omega), \  2\int_{\Omega}\bb{D}(\bb{u}):\bb{D}(\bb{v})\,dx
			=\innerpp{\bb{f},\bb{v}}{\Omega}+\innerpp{ \bb{h},\bb{v}}{\Gamma_{N}}.
		\end{cases}\label{eq:Var_stationary}
		\end{equation}
	\end{prop}

	\begin{proof}[Proof of Proposition \ref{Prop:Steady_equivalence}]
		Let $(\bb{u},p)\in\bb{H}^{1}(\Omega)\times L^{2}(\Omega)$ be a solution of \eqref{eq:Stokes_stacionary}, and let 
		$\bb{v}\in\bb{V}(\Omega)$. Integrating by part
		with the help of \cite[Lemma 2.4]{amrouche2014lp}, we obtain:
		\begin{equation}
		\int_{\Omega}2\bb D(\bb u):\bb D(\bb v)\,dx=\innerpp{f,v}{\Omega}+\innerpp{2\bb D(\bb u)\bb n,v}{\partial \Omega}\ .
\label{eq:weak_varform}		
		\end{equation}
		
	 Now, since $\bb{v}\lvert_{\Gamma_{D}}=0$, it follows from Theorem \ref{th.trace_h_00} that $\bb{v}\lvert_{\Gamma_{N}} \in \bb{H}^{1/2}_{00}(\Gamma_{N})$, thus using Theorem \ref{Th:Green_identity} we are able to perform the following decomposition:

		\begin{align}
		\begin{split}
			\int_{\partial\Omega}2\bb{D}(\bb{u})\bb{n}\cdot \bb{v}\,dS & =\int_{\Gamma_{D}}2\bb{D}(\bb{u})\bb{n}\cdot \bb{v}\,dS
			+\int_{\Gamma_{N}}2\bb{D}(\bb{u})\bb{n}\cdot \bb{v}\,dS\\
 			& =\int_{\Gamma_{N}}\left[\left(2\bb{D}(\bb{u})\bb{n}\right)\cdot \bb{n}\right]\bb{n}\cdot \bb{v}\,dS
			+\int_{\Gamma_{N}}[2\bb{D}(\bb{u})\bb{n}]_{tg}\cdot \bb{v}\,dS\\
 			& =\int_{\Gamma_{N}}[2\bb{D}(\bb{u})\bb{n}]_{tg}\cdot \bb{v}\,dS=\int_{\Gamma_{N}}\bb{h}\cdot \bb{v}\,dS \ ,  \label{eq:dual_decomposition} \\
 		\end{split}
		\end{align}
		where the integrals are to be understood as dual pairing. Therefore, from equations \eqref{eq:weak_varform} and \eqref{eq:dual_decomposition} we have the following:

	\begin{equation}
		\int_{\Omega}2\bb D(\bb u):\bb D(\bb \varphi)\,dx=\innerpp{f,v}{\Omega}+\innerpp{2\bb D(\bb u)\bb n,\bb \varphi}{\Gamma_{N}}\ . \label{eq:vareq}
	\end{equation}
	Conversely, let $\bb{u}\in\bb{V}(\Omega)$ solution of \eqref{eq:Var_stationary}, and let 
	$\bb{\varphi}\in{\mathcal{\bb{D}}}_{\sigma}(\Omega)=\{\bb{\varphi}\in\bb{{\mathcal D}}(\Omega)\lvert div(\bb{\varphi})=0\}$,
	and notice that we have that (see \cite{amrouche2014lp}):
\[
2\int_{\Omega}\bb D(\bb u):\bb D(\bb \varphi)\,dx=\int_{\Omega}\nabla\bb u:\nabla\bb \varphi\,dx\ .
\]
As a consequence, we have:
\[
\forall\bb{\varphi}\in\bb{{\mathcal D}}_{\sigma}(\Omega),\quad\innerp{ -\Delta\bb{u}-\bb{f},\bb{\varphi}} _{\bb{{\mathcal D}}(\Omega)'\times{\mathcal{\bb{D}}}(\Omega)}=0 \ .
\]
Therefore by De Rham theorem \cite[Theorem 2.1]{amrouche1994decomposition}, there exists a distribution $p\in{\mathcal D}(\Omega)'$
defined uniquely up to an additive constant such that:
\begin{equation}
-\Delta\bb{u}+\nabla p=\bb{f} \ . \label{eq:distr}
\end{equation}
Note that since $\bb{H}_{0}(div,\Omega)'$ is embedded in
$\bb{H}^{-1}(\Omega)$ we can say that $p\in L^{2}(\Omega)$  (see \cite[Proposition 2.10]{amrouche1994decomposition}). Let
$\bb{v}\in\bb{V}(\Omega)$, using \eqref{eq:vareq}
and \eqref{eq:distr} we obtain that:

\begin{equation}\label{eq:dual_H_00}
\innerpp{[2\bb D(\bb u)\bb n]_{tg}-\bb h,\bb v}{\Gamma_{N}}=0
\end{equation}

\noindent
Now let $\bb{\alpha} \in \bb{H}^{1/2}_{00}(\Gamma_{N})$ and let $\widetilde{\bb{\alpha}}$ be its extension by zero to $\partial \Omega$, it follows from Theorem \ref{th.extension} that $\widetilde{\bb{\alpha}} \in \bb{H}^{1/2}(\partial \Omega)$. So, there exists $\bb{\varphi} \in \bb{H}^{1}(\Omega)$ such that:

\begin{equation*}
\begin{cases}
div(\bb{\varphi})=0 & \text{ in }\Omega \ ,\\
\varphi=\widetilde{\bb{\alpha}}_{tg} & \text{ on }\Gamma \ .
\end{cases}
\end{equation*}

\noindent
It is clear that $\bb{\varphi} \in \bb{V}(\Omega)$, therefore, using the compatibility condition \eqref{cond:compatiblity_h} and the equation \eqref{eq:dual_H_00}, we have that:

\begin{align*}
\innerpp{[2\bb D(\bb u)\bb n]_{tg}-\bb h,\bb{\alpha}}{\Gamma_{N}} & =\innerpp{[2\bb D(\bb u)\bb n]_{tg}-\bb h,\bb{\alpha}_{tg}}{\Gamma_{N}}\\
 & =\innerpp{[2\bb D(\bb u)\bb n]_{tg}-\bb h,\bb{\varphi}}{\Gamma_{N}} \\
 & =0 \ .
\end{align*} 
Therefore:
\[
[2\bb{D}(\bb{u})\bb{n}]_{tg}=\bb{h}\quad on\,\ \Gamma_{N} \ .
\]
\end{proof}
	To prove existence we shall introduce a Korn--type inequality, whose proof can be done following the scheme given, 
	for instance, in \cite[Theorem 5.3.4]{allaire2007numerical}. 
	\begin{prop}\label{prop:Korn} 
		Let $\Omega$ be a Lipschitz bounded domain. Then there exists a constant $C>0$ depending  
		only on $\Omega$ such that
		\begin{equation}\label{eq:Pointcarr=0000E9Korn}
			\norm{\bb{u}}{\bb{H}^{1}(\Omega)}^{2}\leq C\norm{\bb{D}(\bb{u})}{L^{2}(\Omega)}^{2} 
			   \quad\forall\bb{u}\in\bb{H}_{0,\Gamma_{D}}(\Omega).
		\end{equation}
	\end{prop}

	\begin{prop}\label{thm:Steady_existence} Suppose that $\bb{g}=\bb{0}$,
and let $\bb{f}\in\bb{H}_{0}(div,\Omega)'$, $\bb{h}\in (\bb{H}^{1/2}_{00}(\Gamma_{N}))'$.
Then the problem of finding $(\bb{u},p)\in\bb{H}^{1}(\Omega)\times L_{0}^{2}(\Omega)$
satisfying \eqref{eq:Stokes_stacionary} in the distribution
sense has a unique solution. Furthermore, we have the following estimate:
\begin{equation}
\norm{\bb{u}}{\bb{H}^{1}(\Omega)}+\norm{p}{L_{0}^{2}(\Omega)}\leq C\left(\norm{\bb{f}}{\bb{H}_{0}(div,\Omega)'}+\norm{\bb{h}}{(\bb{H}^{1/2}_{00}(\Gamma_{N}))'}\right) \ . \label{eq:Estimation_stationary}
\end{equation}
\end{prop}

\begin{proof}[Proof of Proposition \ref{thm:Steady_existence}]
Using the proposition \ref{Prop:Steady_equivalence}, we just have to
prove the existence and uniqueness of the following variational
problem:

\[
\begin{cases}
\text{Find \ensuremath{\bb u\in\bb V(\Omega)} such that}\\
\forall\bb v\in\bb V(\Omega),\,\int_{\Omega}2\bb D(\bb u):\bb D(\bb v)\,dx=\innerpp{\bb f,\bb v}{\Omega}+\innerpp{\bb h,\bb v}{\Gamma_{N}}\ .
\end{cases}
\]
We set
\[
a(\bb{u},\bb{v})=2\int_{\Omega}\bb{D}(\bb{u}):\bb{D}(\bb{v})\,dx \ ,
\]
and 
\[
L(\bb v)=\innerpp{\bb f,\bb v}{\Omega}+\innerpp{\bb h,\bb v}{\Gamma_{N}}\ .
\]
Notice that the bi-linear form $a(.\,,.)$ is continuous in $\bb{V}(\Omega)$ since:
\[
\abs{\int_{\Omega}\bb{2}\bb{D}(\bb{u}):\bb{D}(\bb{v})\,dx}=\abs{\int_{\Omega}\nabla\bb{u}:\nabla\bb{v}\,dx}\leq\norm{\nabla\bb{u}}{L^{2}(\Omega)}\,\norm{\nabla\bb{u}}{L^{2}(\Omega)} \ .
\]
Using proposition \ref{prop:Korn} and the fact
that $\bb{V}(\Omega)\subset\bb{H}_{0,\Gamma}(\Omega)$
we also obtain the coercivity of the bilinear form $a(.\,,.)$ on $\bb{V}(\Omega)$.
Now since:
\begin{align}
\begin{split}
\innerpp{\bb f,\bb v}{\Omega} & \leq\norm{\bb f}{\bb H_{0}(div,\Omega)'}\norm{\bb v}{\bb H_{0}(div,\Omega)}\\
&\leq\norm{\bb f}{\bb H_{0}(div,\Omega)'}\norm{\bb v}{\bb H^{1}(\Omega)}
\end{split}
\label{inq:h_div}
\end{align}
\noindent
Using the fact that $v\lvert_{\Gamma_D} =0 $ from theorem \ref{th.trace_h_00} we know that the mapping  $\bb v\lvert_{\Gamma_{N}}:\bb V(\Omega)\to\bb H_{00}^{1/2}(\Gamma_{N})$ is continuous, thus there exist a constant $C_0 >0$ such that:

\begin{equation}
\innerpp{\bb h,\bb v}{\Gamma_{N}}\leq C_{0}\norm{\bb h}{\bb H_{00}^{1/2}(\Gamma_{N})'}\norm v{\bb H^{1}(\Omega)} \ ,
\label{inq:h_00}
\end{equation}

\noindent
Therefore $L:V\rightarrow\mathbb{R}$ is continuous, and thus by
the Lax-Milgram theorem and the proposition \ref{Prop:Steady_equivalence},
we can guarantee the existence and uniqueness of $(\bb{u},p)\in\bb{H}^{1}(\Omega)\times L_{0}^{2}(\Omega)$
satisfying \eqref{eq:Stokes_stacionary} in the distribution
sense. 

To obtain continuity estimates $(\ref{eq:Estimation_stationary})$ notice that by proposition \ref{prop:Korn} and inequalities \eqref{inq:h_div},\eqref{inq:h_00} there exist a constant $C_{1}>0$ such that:
\begin{align}
\norm{\bb{u}}{\bb{H}^{1}(\Omega)} & \leq C_{1}\left(\norm{\bb{f}}{\bb{H}_{0}(div,\Omega)'}+\norm{\bb{h}}{ (\bb{H}^{1/2}_{00}(\Gamma_{N}))'}\right) \ . \label{eq:vel_estimation_stationary}
\end{align}
Now notice that from \eqref{eq:distr} we have:

\begin{align*}
\norm{\nabla p}{\bb{H}^{-1}(\Omega)} & \leq\norm{\bb{f}}{\bb{H}^{-1}(\Omega)}+\norm{\Delta\bb{u}}{\bb{H}^{-1}(\Omega)}\\
 & \leq C_{1}\norm{\bb{f}}{\bb{H}_{0}(div,\Omega)'}+C_{2}\norm{\bb{u}}{\bb{H}^{1}(\Omega)} \ .
\end{align*}
Therefore, from \cite[Proposition 2.10]{amrouche1994decomposition} we obtain that:
\begin{equation}
\norm{p}{L_{0}^{2}(\Omega)}\leq C_{3}\left(\norm{\bb{f}}{\bb{H}_{0}(div,\Omega)'}+\norm{\bb{u}}{\bb{H}^{1}(\Omega)}\right)\ .\label{eq:pressure estimation}
\end{equation}
Using $(\ref{eq:vel_estimation_stationary})$ and $(\ref{eq:pressure estimation})$
we obtain the desired estimate.
\end{proof}

	Now we are in a position to prove the main result of this subsection, Theorem \ref{th.steadycase}.

\begin{proof}[Proof of theorem \ref{th.steadycase}]

Let $\widetilde{\bb g}$ be the extension by zero to $\partial\Omega$
of $\bb g$, from Theorem \ref{th.extension} it follows that $\widetilde{\bb{\bb g}}\in\bb H^{1/2}(\partial\Omega)$
and that there is a constant $C>0$ such that:

\begin{equation}\label{inq:g_extension}
\norm{\widetilde{\bb g}}{\bb H^{1/2}(\partial\Omega)} \leq C \norm{\bb g}{\bb H_{00}^{1/2}(\Gamma_{D})}
\end{equation}

\noindent
Using the compatibility condition \eqref{cond:compatiblity_g}, let $(\bb{w,}\pi)\in\bb{H}^{1}(\Omega)\times L_{0}^{2}(\Omega)$
be the solution of
\begin{equation}
\begin{cases}
-div(\sigma(\bb{w},\pi))=0 &\mbox{in }  \Omega,\\
div(\bb{w})=0 &\mbox{in }  \Omega,\\
\bb{w}=\widetilde{\bb g} &\mbox{on }  \partial \Omega.\\
\end{cases}
\end{equation}

Therefore, using the classical estimate and inequality \eqref{inq:g_extension} there is a constant $C_{0}>0$ such that:

\[
\norm{\bb{w}}{\bb{H}^{1}(\Omega)}+\norm{\pi}{L^{2}(\Omega)}\leq C_{0}\norm{\bb{g}}{ \bb{H}^{1/2}_{00}(\Gamma_{D}) } \ .
\]
Now notice that $\bb{w}\in \bb{E}(\Omega)$ since $\pi\in L^{2}(\Omega)$
therefore $[2\bb{D}(\bb{w})\bb{n}]_{tg}\lvert_{\Gamma_{N}}\in\bb{H}^{1/2}_{00}(\Gamma_{N})'$
due to Theorem \ref{th.dual_tensor}. Thus we obtain the following inequality: 
\begin{align}\begin{split}\norm{[2\bb D(\bb w)\bb n]_{tg}}{(\bb H_{00}^{1/2}(\Gamma_{N}))'} & \leq C_{1}\norm{\bb w}{E(\Omega)}\\
 & \leq C_{2}\left(\norm{\bb{w}}{\bb{H}^{1}(\Omega)}+\norm{\pi}{L^{2}(\Omega)}\right)\\
 & \leq C_{3}\norm{\bb g}{\bb H_{00}^{1/2}(\Gamma_{D})}\ .
\end{split}
\label{eq:tensor_bound}\end{align}
Now since $\bb{h}-[2\bb{D}(\bb{w})\bb{n}]_{tg}\in (\bb{H}^{1/2}_{00}(\Gamma_{N}))'$,
using theorem \ref{thm:Steady_existence}, we denote $(\bb{\widetilde{u}},\widetilde{\pi})\in\bb{H}^{1}(\Omega)\times L_{0}^{2}(\Omega)$
to be a solution of
	\begin{equation}
	\left\{
	\begin{array}{lll}
		-div(\sigma(\bb{\widetilde{u}},\widetilde{\pi}))=\bb{f} &\mbox{in } & \Omega,\\
		div(\bb{\widetilde{u}})=0 &\mbox{in }  &\Omega,\\
		\bb{\widetilde{u}}\cdot \bb{n}=0 &\mbox{in } & \partial\Omega,\\
		\bb{\widetilde{u}}=0 &\mbox{in }&  \Gamma_{D},\\{}
		[2\bb{D}(\bb{\widetilde{u}})\bb{n}]_{tg}=\bb{h}-[2\bb{D}(\bb{w})\bb{n}]_{tg} &\mbox{in } & \Gamma_{N}.
	\end{array}\right.
	\end{equation}
Defining $\bb{u}=\bb{\widetilde{u}}+\bb{w}$, $p=\widetilde{\pi}+\pi$,
we obtain the desired existence and uniqueness result. Using estimation
\eqref{eq:Estimation_stationary} and equation \eqref{eq:tensor_bound}
and the triangle inequality we obtain the desired estimates.
\end{proof}
\subsection{Unsteady Stokes system}
The aim of this subsection is to prove the existence and uniqueness of weak solution for the unsteady 
Stokes system with mixed boundary conditions:
	\begin{equation}\label{eq:Stokes_unsteady}
	\left\{
	\begin{array}{lll}
		\bb{u}_{t}-div(\sigma(\bb{u},p))=\bb{f} &\mbox{in }& \Omega\times(0,T),\\
		div(\bb{u})=0 &\mbox{in } & \Omega\times(0,T),\\
		\bb{u}\cdot \bb{n}=0 &\mbox{on } & \partial\Omega\times(0,T),\\
		\bb{u}=\bb{g} &\mbox{on } & \Gamma_{D}\times(0,T),\\{}
		[2\bb{D}(\bb{u})\bb{n}]_{tg}=\bb{h} &\mbox{on } & \Gamma_{N}\times(0,T),\\
		\bb{u}(0,x)=\bb{u}_{0}(x) &\mbox{in } & \Omega \ .
	\end{array}\right.
	\end{equation}
Since we are dealing with a non-stationary case, we will need an additional space:
\[
\bb{H}=\{\bb{u}\in\bb{L}^{2}(\Omega)\,:\, div(\bb{u})=0\}.
\]	

\noindent	
The main result of this section is given in the following theorem. 
%
%
%
\begin{thm}
\label{thm:Unsteady_existence}Let $\bb{f}\in L^{2}(0,T;\bb{H}_{0}(div,\Omega)')$, $\bb{h}\in L^{2}(0,T;(\bb{H}^{1/2}_{00}(\Gamma_{N}))')$,\\
$\bb{g}\in L^{2}(0,T;\bb{H}^{1/2}_{00}(\Gamma_{D}))$ and $\bb{u_{0}}\in\boldsymbol{H}$
such that: 

\begin{equation}\label{eq:normal_boundary_condition} 
\bb{g}(t)\cdot \bb{n}\lvert_{\Gamma_{D}}=0
\text{ and }
\bb{h}(t)\cdot \bb{n}\lvert_{\Gamma_{N}}=0 \quad \forall t \in (0,T)\ .
\end{equation}
Then the problem of finding $(\bb{u},p)\in L^{2}(0,T;\bb{V}(\Omega))\cap C([0,T],\bb{H}(div,\Omega))\times L^{2}(0,T;L_{0}^{2}(\Omega))$
satisfying \eqref{eq:Stokes_unsteady} in the distribution
sense has a unique solution. 
\end{thm}

\begin{rem}
Notice that if we define $\bb{V}_{0}(\Omega)=\{\bb{u}\in H_{0}^{1}(\Omega)\,:\,div(\bb{u})=0\}$,
it is clear that $\bb{V}_{0}(\Omega)\subset\bb{V}(\Omega)\subset\bb{H}$
and we know by \cite{temam2001navier} that $\bb{V}_{0}(\Omega)$
is dense in $\bb{H}$ therefore $\bb{V}(\Omega)$
is dense in $\bb{H}$.
\end{rem}

\noindent
Before presenting the proof of theorem \ref{th.steadycase}, some preliminaries results are needed. The following proposition characterizes the distributional solution of \eqref{eq:Stokes_unsteady}
in terms of weak solutions.

\begin{prop}
\label{Prop:Unsteady_equivalence}Suppose that $\bb{g}=0$,
and let $\bb{f}\in L^{2}\left(0,T;\bb{H}_{0}(div,\Omega)'\right)$, \\ $\bb{h}\in L^{2}(0,T;(\bb{H}^{1/2}_{00}(\Gamma_{N}))')$ and $\bb{u_{0}}\in\bb{H}$  such that:

\[
\bb{h}(t)\cdot \bb{n}\lvert_{\Gamma_{N}}=0 \quad \forall t \in (0,T)\ .
\]

Then the problem of finding $(\bb{u},p)\in L^{2}(0,T;\bb{V}(\Omega))\cap C([0,T],\bb{H}(div,\Omega))\times L^{2}(0,T;L^{2}_{0}(\Omega))$
satisfying \eqref{eq:Stokes_unsteady} in the distribution
sense is equivalent to:
\begin{equation}
\begin{cases}
\text{Find }u\in L^{2}(0,T;\bb V(\Omega))\cap C([0,T],\bb H(div,\Omega))\,\mbox{such that }\forall\bb v\in\bb V(\Omega),\\
\frac{d}{dt}\int_{\Omega}\bb u\cdot\bb v\,dx+\int_{\Omega}2\bb D(\bb u):\bb D(\bb v)\,dx=\innerpp{\bb f,\bb v}{\Omega}+\innerpp{\bb h,\bb v}{\Gamma_{N}}\\
\text{in the scalar distributional sense on }(0,T),\mbox{ and }\\
\bb u(0)=\bb{u_{0}}.
\end{cases}\label{eq:Var_unsteady}
\end{equation}
\end{prop}
\begin{proof}
Let $\bb{u}\in L^{2}(0,T;\bb{V}(\Omega))\cap C([0,T],\bb{H})$
and $p\in L^{2}(0,T;L^{2}_{0}(\Omega))$ be a solution of (\ref{eq:Stokes_unsteady}),
and let $\bb{v}\in\bb{V}(\Omega)$, integrating by
part as in Proposition \ref{Prop:Steady_equivalence} we obtain:
\begin{equation}
\frac{d}{dt}\int_{\Omega}\bb u.\bb v\,dx+\int_{\Omega}2\bb D(\bb u):\bb D(\bb v)\,dx=\innerpp{\bb f,\bb v}{\Omega}+\innerpp{\bb h,\bb v}{\Gamma_{N}}\ .
  \label{eq:Var_unsteady_1}
\end{equation}
 Conversely, let $\bb{u}\in L^{2}(0,T;\bb{V}(\Omega))\cap C([0,T],\bb{H}(div,\Omega))$
be a solution of \eqref{eq:Var_unsteady}. Let $\bb{\varphi}\in {\calbb{D}}_{\sigma}(\Omega)=\{\bb{\varphi}\in\bb{{\mathcal D}}(\Omega)\,: \,div(\bb{\varphi})=0\}$,
and notice that we have the following:

\[ 
\int_{\Omega}\frac{d\bb u}{dt}(t)\cdot\bb{\varphi}\,dx+\int_{\Omega}\nabla(\bb u(t)):\nabla(\bb{\varphi})\,dx=\innerpp{\bb f(t),\bb \varphi}{\Omega}\ .
\]
As a consequence, we have:
\[
\forall\bb{\varphi}\in\bb{{\mathcal{D}}}_{\sigma}(\Omega),\quad\innerp{ \frac{d\bb{u}}{dt}(t)-\Delta\bb{u}(t)-\bb{f}(t),\bb{\varphi}} _{{\calbb{D}}(\Omega)'\times{\calbb{D}}(\Omega)}=0 \ .
\]
Therefore, thanks to  De Rham Theorem \cite[Theorem 2.1]{amrouche1994decomposition}, there exists a distribution $p(t)\in{\mathcal D}(\Omega)'$
defined uniquely up to an additive constant such that:
\begin{equation}
\frac{d\bb{u}}{dt}(t)-\Delta\bb{u}(t)+\nabla p(t)=\bb{f}(t)\ . \label{eq:distr_unsteady}
\end{equation}

\noindent
Note that since $L^{2}(0,T;\bb{H}_{0}(div,\Omega)')$ is
embedded in $L^{2}(0,T;\bb{H}^{-1}(\Omega))$(i.e. $\bb{H}_{0}(div,\Omega)'\hookrightarrow \bb{H}^{-1}(\Omega)$),
we have $p \in L^{2}(0,T;L^{2}_{0}(\Omega))$.

Now let $\bb{v}\in\bb{V}(\Omega)$, using equations (\ref{eq:Var_unsteady_1}),\eqref{eq:distr_unsteady}, and the same argument as in proposition \ref{Prop:Steady_equivalence}
we obtain
\[
\innerpp{[2\bb D(\bb u(t))\bb n]_{tg}-\bb h(t),\bb v}{\Gamma_{N}}=0 \ .
\]

\noindent
By assumption $\bb{u}(0)=\bb{u}_{0}$ and 
$\bb{u}\in L^{2}(0,T;\bb{V}(\Omega))\cap C([0,T],\bb{H}(div,\Omega))$.
This completes the proof of Proposition \ref{Prop:Unsteady_equivalence}. 
\end{proof}
To proceed we will prove the existence of the solution (\ref{eq:Var_unsteady}). In order to do that, we use 
\cite[Chap 3, Theorem 4.1]{lions2012non}. 
\begin{lem}
\label{lem:Unsteady_existence}Suppose that $\bb{g}=0$, and
let $\bb{f}\in L^{2}(0,T;\bb{H}_{0}(div,\Omega)')$, $\bb{h}\in L^{2}(0,T; (\bb{H}^{1/2}_{00}(\Gamma_{N}))')$, $\bb{u_{0}}\in\bb{H}$ such that:

\[
\bb{h}(t)\cdot \bb{n}\lvert_{\Gamma_{N}}=0 \quad \forall t \in (0,T)\ .
\]

Then the problem of finding $(\bb{u},p)\in L^{2}(0,T;\bb{V}(\Omega))\cap C([0,T],\bb{H}(div,\Omega))\times L^{2}(0,T;L_{0}^{2}(\Omega))$
satisfying \eqref{eq:Stokes_unsteady} in the distribution
sense has a unique solution. 
\end{lem}

\begin{proof}
Thanks to  proposition \ref{Prop:Unsteady_equivalence} , we only need to prove the existence and uniqueness of the following initial variational problem:
\[
\begin{cases}
\text{Find }\bb u\in L^{2}(0,T;\bb V(\Omega))\cap C([0,T],\bb H)\mbox{such that}\forall\bb v\in V(\Omega)\\
\frac{d}{dt}\int_{\Omega}\bb u\cdot\bb v\,dx+\int_{\Omega}2\bb D(\bb u):\bb D(\bb v)\,dx=\innerpp{\bb f,\bb v}{\Omega}+\innerpp{\bb h,\bb v}{\Gamma_{N}}\\
\text{in the scalar distributional sense on }(0,T),\mbox{ and }\\
\bb u(0)=\bb{u_{0}}.
\end{cases}
\]
We will use \cite[Chap 3, Theorem 4.1]{lions2012non}. In that spirit, note that
we have $\bb{V}(\Omega)\subset\bb{H}$ and also we
know that $\bb{V}(\Omega)$ is dense $\bb{H}$. Since
we have already proven in proposition (\ref{thm:Steady_existence}) that the bi-linear
form $a(u,v)=2\int_{\Omega}\bb{D}(\bb{u}):\bb{D}(\bb{v})\,dx$
is continuous and coercive over $\bb{V}(\Omega)$. Notice
also that $L(t;v)=\innerpp{\bb f(t),\bb v}{\Omega}+\innerpp{\bb h(t),\bb v}{\Gamma_{N}}$
is a continuous linear function belonging to $L^{2}(0,T;\bb{V}(\Omega)')$.
Therefore by \cite[Chap 3, Theorem 4.1]{lions2012non} and proposition \ref{Prop:Unsteady_equivalence}, there exists a unique
$\ensuremath{(\bb{u},p)\in L^{2}(0,T;\bb{V}(\Omega))\cap C([0,T],\bb{H})\times L^{2}(0,T;L_{0}^{2}(\Omega))}$
solution of (\ref{eq:Var_unsteady}), furthermore,  $\ensuremath{\frac{du}{dt}\in L^{2}(0,T;\bb{V}(\Omega)')}$.
\end{proof}

\noindent
Now we are in position to prove the main result of this subsection, theorem \ref{thm:Unsteady_existence}.

\begin{proof}[Proof of theorem \ref{thm:Unsteady_existence}]
Let $(\bb w,\pi) \in L^{2}(0,T;\bb{V}(\Omega))\bigcap C([0,T],\bb{H}(div,\Omega))\times L^{2}(0,T;L_{0}^{2}(\Omega))$ be the 
unique solution of (thanks to the compatibliby condition of $\bb g$ in equation (\ref{eq:normal_boundary_condition}))
\[
\begin{cases}
u_{t}-div(\sigma(\bb{w},\pi))=0 & \mbox{in }\Omega\times(0,T),\\
div(\bb{w})=0 & \mbox{in }\Omega\times(0,T),\\
\bb{w}=\bb{g} & \mbox{on }\Gamma_{D}\times(0,T),\\
\bb{w}=0 & \mbox{on }\Gamma_{N}\times(0,T),\\
\bb{u}(0,x)=0 & \mbox{in }\Omega \ .
\end{cases}
\]

Using an analogous argumentation as in theorem \ref{th.steadycase} we have that $\bb{h}-[2\bb{D}(\bb{w})\bb{n}]_{tg}\lvert_{\Gamma_N}\in L^2(0,T;(\bb{H}^{1/2}_{00}(\Gamma_{N})))'$.
From lemma \ref{lem:Unsteady_existence}, let $(\widetilde{\bb{u}},\widetilde{\pi})\in L^{2}(0,T;\bb{V}(\Omega))\bigcap C([0,T],\bb{H}(div,\Omega))\times L^{2}(0,T;L_{0}^{2}(\Omega))$ be the unique solution to
\[
\begin{cases}
\widetilde{\bb{u}}_{t}-div(\sigma(\widetilde{\bb{u}},\widetilde{\pi}))=0 & \mbox{in }\Omega\times(0,T),\\
div(\widetilde{\bb{u}})=0 & \mbox{in }\Omega\times(0,T),\\
\widetilde{\bb{u}}\cdot n=\bb{0} & \mbox{on }\partial\Omega\times(0,T),\\
\widetilde{\bb{u}}=0 & \mbox{on }\Gamma_{D}\times(0,T),\\{}
[2\bb{D}(\widetilde{\bb{u}})\bb{n}]_{tg}=\bb{h}-[2\bb{D}(\bb{w})\bb{n}]_{tg} & \mbox{on }\Gamma_{N}\times(0,T),\\
\widetilde{\bb{u}}(0,x)=\bb{u}_{0} & \mbox{in }\Omega \ .
\end{cases}
\]
Defining $\bb{u}=\widetilde{\bb{u}}+\bb{w}$,
$p=\widetilde{\pi}+\pi$, we obtain the desired existence and uniqueness
result. 

\end{proof}

\begin{rem}
As far as we know the regular solution with mixed boundary condition, as in \eqref{eq:Stokes_unsteady}, for
general Lipschitz domains is an open problem. However some results for the polygonal case
and for similar boundary conditions are available in \cite{brown2010mixed} using potential theory. The generalization of
these result to our particular case is out of the scope of this paper and will be subject to further work.	
\end{rem}

\section{Theoretical results to the identification problem}
	In this section we prove the identifiability of the obstacle $\mathcal{O}$ associated to the Stokes system with mixed boundary 
	\eqref{eq:Stokes-navier-slip}. Since in our case the obstacle (rigid body) is in contact with part
	of the boundary $\partial\Omega$ (see \eqref{cond.borde}), we first define the admissible 
	deformation of the domain.
	Secondly, the proof of the main result is presented, theorem \ref{identification_theo}. Our proof 
	is based on the arguments given 
	in \cite{alvarez2005identification}, where the authors proved the identification of immersed 
	obstacles from boundary 
	data for the Navier--Stokes system with  Dirichlet boundary conditions. It is worth mentioning that our main novelties rely on 
	the geometrical configuration, since the obstacle will intersect the flow boundary, also the mixed boundary conditions are considered instead of Dirichlet boundary conditions.

	\begin{defn}\label{def:admisible_defo}
		Let $\Omega\subset\mathbb{R}^{2}$ be a simple connected, bounded and $C^{1,1}$ domain, 
		$\partial\Omega=\Gamma_{inlet}\cup\Gamma_{wall}\cup\Gamma_{out}$.
		A domain $\Omega_{{\mathcal O}}$ is called an admissible deformation
		of $\Omega$ if and only if: 
		\begin{enumerate}
			\item [i)] $\Omega_{{\mathcal O}}\subset\Omega$ is simple a connected and $C^{1,1}$ domain.
			\item [ii)]  $\partial\Omega_{{\mathcal O}}:=\Gamma_{inlet}\cup\Gamma^{{\mathcal O}}_{wall}\cup\Gamma_{out}$
				and $\Gamma_{wall}\cap\Gamma^{{\mathcal O}}_{wall}$ are non--empty
				relatively open sets of $\Gamma_{wall}$. 
			
			\item [iii)] There exist a relative open set 
				$W \subset\Omega$ satisfying  $\Omega\backslash \ovl{\Omega_{{\mathcal O}}}\subset W$,
				there exists a diffeomorphism $\psi:\ovl{\Omega}\to\ovl{\Omega}$ 
				such that $\psi(\Omega)=\Omega_{{\mathcal O}}$ and 
				$\psi=I$ in $\Omega\backslash\ovl{W}$ (see Figure \ref{fig:Diffeo_example}).
		\end{enumerate}

	\begin{figure}[htbp]
	\begin{centering}
		\includegraphics[scale=0.15]{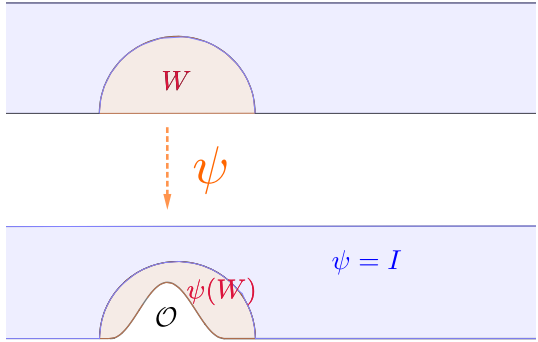}
		\par\end{centering}
		\caption{\label{fig:Diffeo_example} Diffeomorphism example }
	\end{figure}

	\end{defn}
	\begin{rem}
		From definition \ref{def:admisible_defo}, given an admissible deformation $\Omega_{{\mathcal O}}$, an admissible obstruction 
		can be define as ${\mathcal{O}}=\Omega \backslash \ovl{\Omega_{{\mathcal O}}}$.
	\end{rem}

The following theorem is due to Fabre and Lebeau in \cite{fabre1996prolongement}. It is worth pointing out that this result is independent of the boundary conditions and thus it can be used in particular for Navier-slip boundary conditions.

\begin{thm}\label{theo.fabre} 
Let $\Omega\subset\mathbb{R}^{N}$ be a connected open set, $N \geq 2$
and $T>0$. Let $a\in L_{loc}^{\infty}(\Omega\times(0,T))^{N}$ and
$c\in C([0,T];L_{loc}^{r}(\Omega,\mathbb{R}^{N\times N}))$ be a matrix-value
function with $r>N$. If $(\bb{v},p)\in L^{2}(0,T_{;}H_{loc}^{1}(\Omega))\times L_{loc}^{2}(\Omega\times(0,T))$
is a solution of:
	\begin{equation}
	\left\{
	\begin{array}{lll}
  		\bb{v}_{t}-\Delta \bb{v}+(a.\nabla)\bb{v}+c\bb{v}+\nabla p=0 &\mbox{in }& \Omega\times(0,T),\\
 		 div(\bb{v})=0 &\mbox{in }& \Omega\times(0,T),
	\end{array}\right.
	\end{equation}
with $\bb{v}=0$ in $\omega_{0}\times(0,T)$, where $\omega_{0}$ is and
open set of $\Omega$.
Then $\bb{v}=0$ in $\Omega\times(0,T)$ and $p$ is a constant.
\end{thm} 
\medskip

\noindent Using the above theorem, we have the following corollary, which is essential for the following identification theorem. We note that this result can be also in  \cite[Corollary 2.4]{alvarez2005identification} and a similar one in \cite{boulakia2012stability}. Since no prove is given \cite{alvarez2005identification} we present a short prove.

\begin{corollary} \label{cor.fabeau}
Let $\Omega\subset\mathbb{R}^{N}$ be a connected open lipschitz domain, $N \geq 2$ and $T>0$. 
If $(\bb{u},p)\in L^{2}(0,T;H^{1}(\Omega)^{N})\times L^{2}(\Omega\times(0,T))$
is a solution of:
	\begin{equation}\label{cor.fabeau.stokes}
	\left\{
	\begin{array}{lll}
  		 \bb{u}_{t}-div(\sigma(\bb{u},p))=0 & \mbox{in }& \Omega\times(0,T),\\
   		div(\bb{u})=0 &\mbox{in }&  \Omega\times(0,T),
	\end{array}\right.
	\end{equation}
satisfying $\bb{u}=\sigma(\bb{u},p)\bb{n}=0 \ \,on\ \,\Gamma\times(0,T)$ 
where $\Gamma\subset\partial\Omega$ are relatively open non-empty subsets, 
then
\[
   \bb{u}=0\quad in\quad\Omega\times(0,T) \ .
\]
\end{corollary}

\begin{proof}[Proof of corollary \ref{cor.fabeau}]
Let $\bb{x_0}\in \Gamma$ and $r>0$ such that $\mathcal{\bb{B}}(\bb{x_0},r) \cap \partial \Omega \subset \Gamma $. Let $B=\mathcal{\bb{B}}(\bb{x_0},r)\cap\Omega^{c}$ and $\widetilde{\Omega}=\Omega \cup B$, we define the extension by 0 of $\bb{u}$ and $p$ in $\widetilde{\Omega}$ by:

\[
\widetilde{\bb{u}}\text{ (resp \ensuremath{\widetilde{p}})}=\begin{cases}
\bb{u}\text{ (resp \ensuremath{p})} & \text{in }\Omega\times(0,T),\\
0 & \text{in }B\times(0,T).
\end{cases}
\]

\noindent
Our goal now is to prove that $(\widetilde{\bb{u}},\widetilde{p})$ is solution of \eqref{cor.fabeau.stokes} in $\widetilde{\Omega}$, for this let $\bb{\varphi}\in\bb{D}(\widetilde{\Omega})$ and notice that
\begin{align*}
\left\langle \left(\widetilde{\bb{u}}_{t}-div(\sigma(\widetilde{\bb{u}},\widetilde{p}))\right)\chi_{B},\bb{\varphi}\right\rangle _{\bb{D'}(\widetilde{\Omega})\times\bb{D}(\widetilde{\Omega})} & = \int_{B}\widetilde{\bb{u}}_{t}\,dx+2\int_{B}\bb{D}(\widetilde{\bb{u}}):\nabla \bb{\varphi} \,dx -\int_{B}\widetilde{p}\,div(\bb{\varphi}) \\
&- 2\left\langle \sigma(\widetilde{\bb{u}},\widetilde{p})\bb{n},\bb{\varphi}\right\rangle _{\bb{H}^{-1/2}(\partial B)\times\bb{H}^{1/2}(\partial B)}\\
 & =\int_{B}\widetilde{\bb{u}}_{t}\,dx+2\int_{B}\bb{D}(\widetilde{\bb{u}}):\nabla \bb{\varphi}\,dx-\int_{B}\widetilde{p}\,div(\bb{\varphi})\\
 & \quad -2\left\langle \sigma(\widetilde{\bb{u}},\widetilde{p})\bb{n},\bb{\varphi}\right\rangle _{\bb{H}^{-1/2}(\Gamma)\times\bb{H}^{1/2}(\Gamma)}\\
 & =-2\left\langle \sigma(\widetilde{\bb{u}},\widetilde{p})\bb{n},\bb{\varphi}\right\rangle _{\bb{H}^{-1/2}(\Gamma)\times\bb{H}^{1/2}(\Gamma)}\\
 & =0.
\end{align*}

\noindent
In a similar way we can also deduce that
\[
\left\langle \left(\widetilde{\bb{u}}_{t}-div(\sigma(\widetilde{\bb{u}},\widetilde{p}))\right)\chi_{\Omega},\bb{\varphi}\right\rangle _{\bb{D'}(\widetilde{\Omega})\times\bb{D}(\widetilde{\Omega})}=0.
\]

\noindent
Therefore, we obtain
\[
\left\langle \widetilde{\bb{u}}_{t}-div(\sigma(\widetilde{\bb{u}},\widetilde{p})),\bb{\varphi}\right\rangle _{\bb{D'}(\widetilde{\Omega})\times\bb{D}(\widetilde{\Omega})}=0.
\]

\noindent 
Using Theorem \ref{theo.fabre} for $(\widetilde{\bb{u}},\widetilde{p})$ we obtain the desired result.
\end{proof}

\begin{thm}\label{identification_theo}
Let $\Omega_{0}$,\,$\Omega_{1} \subset \mathbb{R}^{N}$ be two admissible deformation 
of $\Omega$ according to the definition \ref{def:admisible_defo} such that $\Omega_{0}\cap\Omega_{1}$ is a $C^{1,1}$ domain with a finite number of disjoint connected components. 
Let $(\bb{u}_{j},p_{j})\in L^{2}(0,T;H^{2}(\Omega_{j})^{N})\times L^{2}(0,T;L^{2}(\Omega_{j}))$
be solution of:
	\begin{equation}\left\{
	\begin{array}{lll}
		\frac{\partial \bb{u}_{j}}{\partial t}-div(\sigma(\bb{u}_{j},p_{j}))=0 &\mbox{in }&  \Omega_{j}\times(0,T),\\
		div(\bb{u}_{j})=0 &\mbox{in }&  \Omega_{j}\times(0,T),\\
		\bb{u}_{j}=\bb{g_{in}} &\mbox{on }&  \Gamma_{inlet}\times(0,T),\\
		\bb{u}_{j}=\bb{g_{out}} &\mbox{on }&  \Gamma_{out}\times(0,T),\\
		\bb{u}_{j}\cdot\bb{n}=0 ,\ \left(\sigma(\bb{u}_{j},p)\bb{n}\right)_{tg}=0 &\mbox{on }&  \Gamma_{wall}^{j}\times(0,T),\\
		\bb{u}(\cdot,0)=0 &\mbox{in }&  \Omega_{j}.
	\end{array}\right.
		\quad j=0,1
\label{eq:indentificacion.sol}
	\end{equation}
Assume that $\bb{g_{in}},\bb{g_{out}} \neq 0$ and that $\Gamma_{m}\subset \Gamma_{wall}$ is a non-empty relative open 
satisfying $\Gamma_{m}\subset\Gamma_{wall}^{0}\cap\Gamma_{wall}^{1}\subset\partial(\Omega_{0}\cap\Omega_{1})$,
such that 
\[
(\sigma(\bb{u}_{0},p_{0})\bb{n})\cdot\bb{n}=(\sigma(\bb{u}_{1},p_{1})\bb{n})\cdot\bb{n} \,\,\text{ and } \,\,\bb{u}_0\cdot\bb{\tau}=\bb{u}_1\cdot\bb{\tau} \quad\text{on }\Gamma_{m}\times(0,T),
\]
where $\bb{n}$ is external unit normal and $\bb{\tau}$ is tangent unit vector of $\Gamma_{m}$.
Then $\Omega_{0}\equiv\Omega_{1}$.

\end{thm}

\begin{figure}
\centering
\includegraphics[trim={0cm 3.5cm 0cm 3.2cm},clip,width=0.8\textwidth]{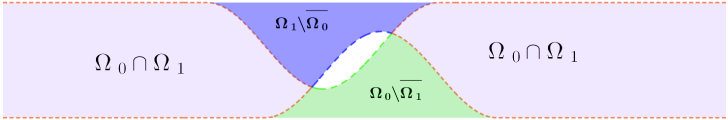}
\par
\caption{Example of domain intersection \label{fig.domain_intersection}}
\end{figure}

\begin{proof}[Proof of theorem \ref{identification_theo}]
First within the context of Theorem \ref{thm:Unsteady_existence} by setting $\Gamma_{D}= \Gamma_{inlet} \cup \Gamma_{out}$ and $\Gamma_{N}= \Gamma^{j}_{wall}$ we can assure the existences and uniqueness of the solutions $u_0,u_1$ of equations (\ref{eq:indentificacion.sol}).
\noindent
Let $\bb{u}=\bb{u}_{0}-\bb{u}_{1}$, $p=p_{0}-p_{1}$, thus we have:
	\begin{equation}
	\left\{
	\begin{array}{lll}
		\frac{\partial \bb{u}}{\partial t}-div(\sigma(\bb{u},p))=0 & \mbox{in }& \Omega_{0}\cap\Omega_{1}\times(0,T),\\
		div(\bb{u})=0 &\mbox{in }&  \Omega_{0}\cap\Omega_{1}\times(0,T),\\
		\bb{u}=0 & \mbox{on }& \Gamma_{inlet}\times(0,T),\\
		\bb{u}=0 & \mbox{on }& \Gamma_{out}\times(0,T),\\
		\bb{u}\cdot \bb{n}=0 ,\ \left(\sigma(\bb{u},p)\bb{n}\right)_{tg}=0 &\mbox{on }&  \Gamma_{m}\times(0,T),\\
		\bb{u}(.,0)=0 & \mbox{on }& \Omega_{0}\cap\Omega_{1}.
	\end{array}\right.
	\end{equation}
\noindent
Since $\sigma(\bb{u},p)\bb{n}=0$, $\bb{u}=0$ on $\Gamma_{m}\times(0,T)$ and $\Gamma_{m} \subset \partial (\Omega_{0}\cap \Omega_{1})$, it
follows from corollary \ref{cor.fabeau} applied to each connected component of $\Omega_{0}\cap \Omega_{1}$ that :

\begin{equation}\label{eq.igualdad_fabre}
\bb{u}_{0}=\bb{u}_{1}\text{ in }(\Omega_{0}\cap\Omega_{1})\times(0,T).
\end{equation}

\noindent
Now lets assume that $\Omega_{0}\backslash\ovl{\Omega_{1}}$
is a non-empty open subset of $\Omega_{0}$. Thus, we have:

\begin{equation}\label{eq:intersection_fluids}
\frac{\partial \bb{u}_{0}}{\partial t}-div(\sigma(\bb{u}_{0},p_{0})\bb{n})=0\quad in\quad(\Omega_{0}\backslash\ovl{\Omega_{1}} )\times (0,T) \ .
\end{equation}

\noindent
Noticing that for any two bounded regular open sets, we have (for an example, see Figure \ref{fig.domain_intersection}):
\[
\partial(\Omega_{0}\setminus\ovl{\Omega_{1}})=\left[\partial(\Omega_{0}\cap\Omega_{1})\cap\partial(\Omega_{0}\setminus\ovl{\Omega_{1}})\right]\cup\left[\partial(\Omega_{0})\cap\partial(\Omega_{0}\setminus\ovl{\Omega_{1}})\right]
\]

\noindent
Thus multiplying equation \eqref{eq:intersection_fluids} by $\bb{u}_{0}$ and integrating by parts in
$\Omega_{0}\backslash\ovl{\Omega_{1}}$ , we obtain:

\begin{align*}
\frac{d}{dt}\int_{\Omega_{0}\backslash\ovl{\Omega_{1}}}\lvert\bb u_{0}(x,t)\rvert^{2}dx &=-\int_{\Omega_{0}\backslash\ovl{\Omega_{1}}}\lvert\bb D(\bb u_{0})\rvert^{2}dx+\int_{\partial\left(\Omega_{0}\backslash\ovl{\Omega_{1}}\right)}[\bb D(\bb u_{0})\bb n]_{tg}\cdot\bb u_{0}\\
& +\int_{\partial\left(\Omega_{0}\backslash\ovl{\Omega_{1}}\right)}p_{0}\bb u_{0}\cdot n
\end{align*}

\noindent
Thus by equation \eqref{eq.igualdad_fabre} we have that:
\[
\begin{cases}
[\bb D(\bb u_{0})\bb n]_{tg}=[\bb D(\bb u_{1})\bb n]_{tg}=0 & \text{on }\partial(\Omega_{0}\cap\Omega_{1})\cap\partial(\Omega_{0}\setminus\ovl{\Omega_{1}})\\
\bb u_{0}\cdot\bb n=\bb u_{1}\cdot\bb n=0 
\end{cases}
\]

\noindent
Also by assumption, we know that $\bb u_{0}\cdot\bb n = [\bb D(\bb u_0)\bb n]_{tg}=0$ on $\partial \Omega_{0}$. Then we obtain:
\[
\frac{d}{dt}\int_{\Omega_{0}\backslash\ovl{\Omega_{1}}} \lvert \bb{u}_{0}(x,t) \rvert ^{2}dx=-\int_{\Omega_{0}\backslash\ovl{\Omega_{1}}} \lvert \bb{D}(\bb{u}_{0}) \rvert ^{2}dx \ .
\]
Therefore,
\[
E(t)=\int_{\Omega_{0}\backslash\ovl{\Omega_{1}}} \lvert \bb{u}_{0}(x,t) \rvert ^{2}dx \ 
\]
is a decreasing non-negative function. However since $\bb{u}_{0}(x,0)=0$, we have that $\bb{u}_{0}\lvert_{\Omega_{0}\backslash\ovl{\Omega_{1}}}=0$
for all $t\in(0,T)$. Therefore, from theorem \ref{theo.fabre} we conclude that:
\[
\bb{u}_{0}=0\,\,\,in\quad\Omega_{0} \ .
\]
Nevertheless, 
this is impossible since $\bb{u}_{0}\neq0$, because $\bb{g_{in}},\,\bb{g_{out}}\neq0$.
Therefore $\Omega_{0}\backslash\ovl{\Omega_{1}}=\emptyset$,
analogously we can deduce that $\Omega_{1}\backslash\ovl{\Omega_{0}}=\emptyset$
which implies that $\Omega_{0}=\Omega_{1}$.
\end{proof}

	As mentioned above, our main propose is to study and solve 
	the full inverse problem from a numerical point of view. In relation to the direct problem,  the smooth dependence of the 
	Cauchy forces for the Stokes system with mixed boundary conditions might be obtained by proving an extension of \cite[theorem 5.1]{alvarez2005identification}. 	Nevertheless,		its proof involves a deep analysis on the holomorphic semigroup associated to the main operator  of (\ref{eq:Stokes-navier-slip}), which as we know, involves mixed boundary conditions. Regarding the inverse wave equation problem it is require to obtain the regularity of the stokes fluid with mixed boundary condition which in it self is open problem. Furthermore to deal with hyperbolic inverse problem a different set of tools has to be used such as the Hilbert uniqueness method, we refer the reader to \cite{lions1988controlabilite,yamamoto1995stability}. Since these topics deserve a special treatment we did not attempt to board them in this article, further work is in progress.


\section{Numerical perspective for the direct problem}
\subsection{Direct problem and experimental setup}\label{subsec:Direct problem and experimental setup}
	In order to solve the evolutionary Stokes equations with Navier-slip boundary conditions 
	related to these problems, we decide to use radial divergence free kernels. The reason for this choice is that these methods 
	are mesh free, letting us to work with more complex geometries in the future, they have spectral convergence properties,
	 it is simple to implement linear boundary conditions such as Navier-slip and we can avoid the usual inf--sup conditions 	
	 \cite{wendland2009divergence}, \cite{keim2016high}. Never the less the down side of these methods are that the condition 
	 numbers of the resulting matrix are high and that spurious eigenvalues can appear for evolutionary problems, as it happens 
	 with some other spectral methods. To overcome this problem, we use divergency free hybrid kernels which has been 
	 introduced recently for scalar problems \cite{mishra2019stabilized} and for the Stokes system \cite{2021Bretoneatal}.

	To solve the numerical direct problem we decide to use hybrid divergence free kernels. which we defined in this section, and that 
	were recently introduced in \cite{mishra2018hybrid}. These kernels, in its scalar version, are a linear combination of Gaussian and 
	Polyharmonic splines. The Gaussian component contributes to attain exponential convergence, while the polyharmonic part control the 
	stability, namely the growth of the condition number, of the scheme. We refer the reader to article \cite{2021Bretoneatal} for 
	a detailed discussion of the material included in this subsection.
	\vskip 0.3cm
	\noindent Let us define $\bb{L}(\bb{u},p)=-\mu\Delta\bb{u}+\nabla p$ and consider the system:	
	\[
\begin{cases}
\bb{u}_{t}+\bb{L}(\bb{u},p)=\bb{f} & \text{ in }Q,\\
div(\bb{u})=0 & \text{ in }Q,\\
\bb{B}\bb{u}=\bb{g} & \text{ on }\Sigma,\\
\bb{u}(\cdot,0)=\bb{u}_{0}(\cdot) & \text{ in }Q,
\end{cases}
\]
where $\bb{B}$ is a given boundary operator. Using the method of lines combined with the general interpolation theorem, see for example \cite{wendland2009divergence}, we proposed the following ansatz.
\begin{equation}\label{anzats1}
\begin{split}
(\hat{\bb{u}},\hat{p})(\bb{x},t) &=\sum_{i=1}^{N}\sum_{j=1}^{N_{b}}\bb{B}_{i}^{\bb{\xi}}\bb{\Phi}(\bb{x-\bb{\xi}_{j}})\alpha_{(i-1)N_{b}+j}(t)\\
&+\sum_{i=1}^{N}\sum_{j=1}^{N_{in}}\bb{L}_{i}^{\bb{\xi}}\bb{\Phi}(\bb{x-\bb{\xi}_{N_{in}+j}})\alpha_{NN_{b}+(i-1)N_{in}+j}(t),
\end{split}
\end{equation}
where $N_b,\, N_{in}$ are the total number of boundary and interior nodes, respectively, and 
$\bb{\xi}_{1..N_{b}}^{b},\, \bb{\xi}_{N_{b}+1....N_{b}+N_{in}}^{in} \in \RR^N$ are the boundaries and interior centers.
	On the other hand, $\bb{B}_{i}^{\xi}\boldsymbol{\Phi}$, $\bb{L}_{i}^{\xi}\boldsymbol{\Phi}$ are vector--valued functions from $\RR^N \times \RR^N$ to $\RR^{N+1}$ defined as the application of the operator $\bb{B}_{i}^{\xi}$, $\bb{L}_{i}^{\xi}$ to each row of the kernel $\boldsymbol{\Phi}$.
\vskip 0.3cm
\noindent
In 2D, for equation (\ref{anzats1}), $\bb{\Phi}$ is the combined Div--free velocity------pressure kernel given by
\begin{equation}\label{hybrido}
\bb{\Phi}(\mathbf{x})=\left(\begin{array}{cc}
\bb{\Phi}_{Div}(\mathbf{x}) & 0\\
0 & {\mathrm{e}}^{-c_{2}\,r}+\gamma_{2}\,r^{2m+1}
\end{array}\right).
\end{equation}
where $\gamma_{2}$ is the corresponding weight relative to the hybrid kernel related to the pressure and $\bb{\Phi}_{Div}$ is the divergence-free hybrid kernel given by:
\[\bb{\Phi}_{Div}(\mathbf{x})=\Delta\times\nabla\times\psi=\{-\Delta I+\nabla\nabla^{T}\}(\psi_1(\mathbf{x})+\gamma_1 \psi_2(\mathbf{x})).\]
where $\mathbf{r} = \|\mathbf{x}\|$, $\psi_1(\mathbf{x}) = \exp(-c_{1}\,\mathbf{r}^2)$, 
$\psi_2(\mathbf{x})=r^{2n+1}$ and $\gamma_1$ is a positive real number. We note that by direct 
computation it is possible to obtain:

{\small
\begin{align*}
&\bb{\Phi}_{Div}(\mathbf{x}) =\left(\begin{array}{cc}
-2\,c_{1}\,{\mathrm{e}}^{-c_{1}\,r^{2}}\,\left(2\,c_{1}\,{\left(y_{1}-y_{2}\right)}^{2}-1\right) & 4\,c_{1}^{2}\,{\mathrm{e}}^{-c_{1}\,r^{2}}\,\left(x_{1}-x_{2}\right)\,\left(y_{1}-y_{2}\right)\\
4\,c_{1}^{2}\,{\mathrm{e}}^{-c_{1}\,r^{2}}\,\left(x_{1}-x_{2}\right)\,\left(y_{1}-y_{2}\right) & \quad\quad-2\,c_{1}\,{\mathrm{e}}^{-c_{1}\,r^{2}}\,\left(2\,c_{1}\,{\left(x_{1}-x_{2}\right)}^{2}-1\right)
\end{array}\right)\\
 & +\gamma_{1}\left(2n+1\right)\,{r}^{2n-3}\left(\begin{array}{cc}
-\left({\left(y_{1}-y_{2}\right)}^{2}\,\left(2n-1\right)+r^{2}\right)\,\quad\quad & \left(x_{1}-x_{2}\right)\,\left(y_{1}-y_{2}\right)\,\left(2n-1\right)\\
\left(x_{1}-x_{2}\right)\,\left(y_{1}-y_{2}\right)\,\left(2n-1\right) & -\left({\left(x_{1}-x_{2}\right)}^{2}\,\left(2n-1\right)+r^{2}\right)
\end{array}\right).
\end{align*}
}
\vskip 0.3cm
\noindent
Collocating the ansatz (\ref{anzats1}) in the Stokes system, we obtain 
the following ode system:

\begin{align}
M_{\phi}\frac{\ovl{\alpha}(t)}{dt}+M_{L\phi}\ovl{\alpha}(t) & =\mathbf{f} \ , \label{eq:ode_fbr}\\
M_{B_{\phi}}\ovl{\alpha}(t) & =\mathbf{g} \ , \nonumber 
\end{align}
where, for $\bb{B}_{i}^{\xi}\bb{\Phi}=(\phi_{1}^{{B}_{i}},\phi_{2}^{{B}_{i}},\phi_{3}^{{B}_{i}})^t$, $\bb{L}_{i}^{\xi}\bb{\Phi}=(\phi_{1}^{{L}_{i}},\phi_{2}^{{L}_{i}},\phi_{3}^{{L}_{i}})^{t}$ we have $M_{\phi},M_{L_{\phi}}\in \RR^{2N_{in}\times 
	(2N_{b}+2N_{in})}$ and are defined as follows:
	\begin{equation*}
		M_{\phi}=\left(
		\begin{array}{ccccc}
		\phi_{1}^{B_{1}} & \phi_{1}^{B_{2}} &  \phi_{1}^{L_{1}} & \phi_{1}^{L_{2}}\\
		\phi_{2}^{B_{1}} & \phi_{2}^{B_{2}} &  \phi_{2}^{L_{1}} & \phi_{2}^{L_{2}}
		\end{array}\right),
		\,\,\,
		M_{L_{\phi}}=\left(\begin{array}{ccccc}
		L_{1}\phi^{B_{1}} & L_{1}\phi^{B_{2}} & L_{1}\phi^{L_{1}} & L_{1}\phi^{L_{2}}\\
		L_{2}\phi^{B_{1}} & L_{2}\phi^{B_{2}} & L_{2}\phi^{L_{1}} & L_{2}\phi^{L_{2}} 
		\end{array}\right) \ ,
	\end{equation*}
	and $M_{B_{\phi}}\in  \RR^{2n_{b}\times(2n_{b}+2n_{in})}$
	\begin{equation*}
		M_{B_{\phi}}=\left(\begin{array}{ccccc}
		B_{1}\phi^{B_{1}} & B_{1}\phi^{B_{2}} & B_{1}\phi^{L_{1}} & B_{1}\phi^{L_{2}}\\
		B_{2}\phi^{B_{1}} & B_{2}\phi^{B_{2}} & B_{2}\phi^{L_{1}} & B_{2}\phi^{L_{2}}
		\end{array}\right).
	\end{equation*}
The numerical solution of system (\ref{eq:ode_fbr}) is solved by using backward differentiation formula methods (BDF's methods). It is important to note that, as shown numerically in \cite{2021Bretoneatal}, the eigenvalues of the Gram matrix of this system are all on the negative side of the complex plane. Thus, the scheme is stable for backward differentiation formulas.

\subsection{Numerical results for the direct problem}\label{subsec:Numerical_direct_problem}
To present numerical results for problem (\ref{eq:Stokes-navier-slip})--(\ref{eq:Onda}) using the hybrid divergence-free kernel method presented above,  we define the computational domain of the fluid as $\Omega\backslash{\mathcal{O}}=([0,8]\times[0,1])\backslash{\mathcal{O}}$, where the obstruction $\mathcal{O}$  which is part of its boundary is parameterized as follows:

\begin{equation}\label{eq:Dominio Parametrico}
\partial{\mathcal O}=\begin{cases}
x(s)=\theta_{1}+s & s\in[0,\theta_2],\\
y(s)=\frac{\theta_{3}}{2}\left(1.0-cos(\frac{2\pi s}{\theta_{2}})\right) & s\in[0,\theta_2].
\end{cases}
\end{equation}

\begin{figure}[hbt!]
\begin{centering}
\includegraphics*[bb = 50 200 2400 800,scale=0.15]{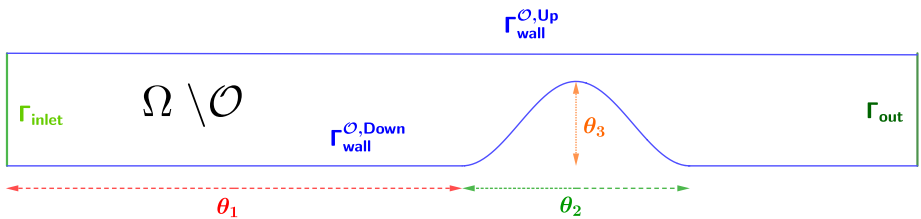}
\par\end{centering}
\centering{}\caption{\label{fig:Dominio Parametrico} Generic Parametric representation of the domain
$\Omega\backslash{\mathcal O}$}
\end{figure}
\noindent

In particular, for this numerical simulation we use the following parameters:

$$\theta_1=4, \quad \theta_2=1, \quad \theta_3=0.5 \ . $$
We recall that the obstruction parameters $\theta_1,\theta_2,\theta_3$ should be understood as the position, size, and percentage depth of the obstruction.
\noindent
For the computational domain of the wave, we set $S=[0,8]\times[1,5]$.

Since the couple problem (\ref{eq:Stokes-navier-slip})--(\ref{eq:Onda}) is indeed only in one direction, for simplicity, we decide to solve (\ref{eq:Stokes-navier-slip}) with hybrid RBF and store the normal component of the stress tensor to later simulate the wave equation using the finite element method due to its simple geometry domain. We note that in all of our examples, the finite element method attained a similar order of precision as the hybrid RBF technique. Thus, the coupled numerical results are consistent. 
For the numerical simulation (\ref{eq:Stokes-navier-slip}) we use a mesh of 1119 halton points, a BDF2 scheme to solve the ode system (\ref{eq:ode_fbr}) with a time step of $\Delta t=\frac{1}{200}$ in the interval $[0,5]$, a fixed viscosity constant $\mu=0.1$. It is important to note that for this simulation, we define $g_{in}(x,t)=g_{out}(x,t)=(cos(2 \pi t + \pi)+1)$. 
 
For the numerical simulation of the wave equation \eqref{eq:Onda} we use the element P1 with a mesh of 5963 points and a fixed propagation speed $c=1$.
 
It should be noted that we took the RBF hybrid parameter as $c_1=c_2=0.5$, $\gamma_{1}=10^{-4}$ and $\gamma_{2}=10^{-4}$ since experimentally this combination gives us a reasonable condition number of $10^{14}$ for the associated ode matrix system (\ref{eq:ode_fbr}) and an accuracy of $10^{-4}$ for the vector field and $10^{-3}$ for the pressure scalar field.

\begin{figure}[hbt!]
\centering{}
\subfloat[]{\includegraphics*[bb = 275 0 1240 850,scale=0.17]{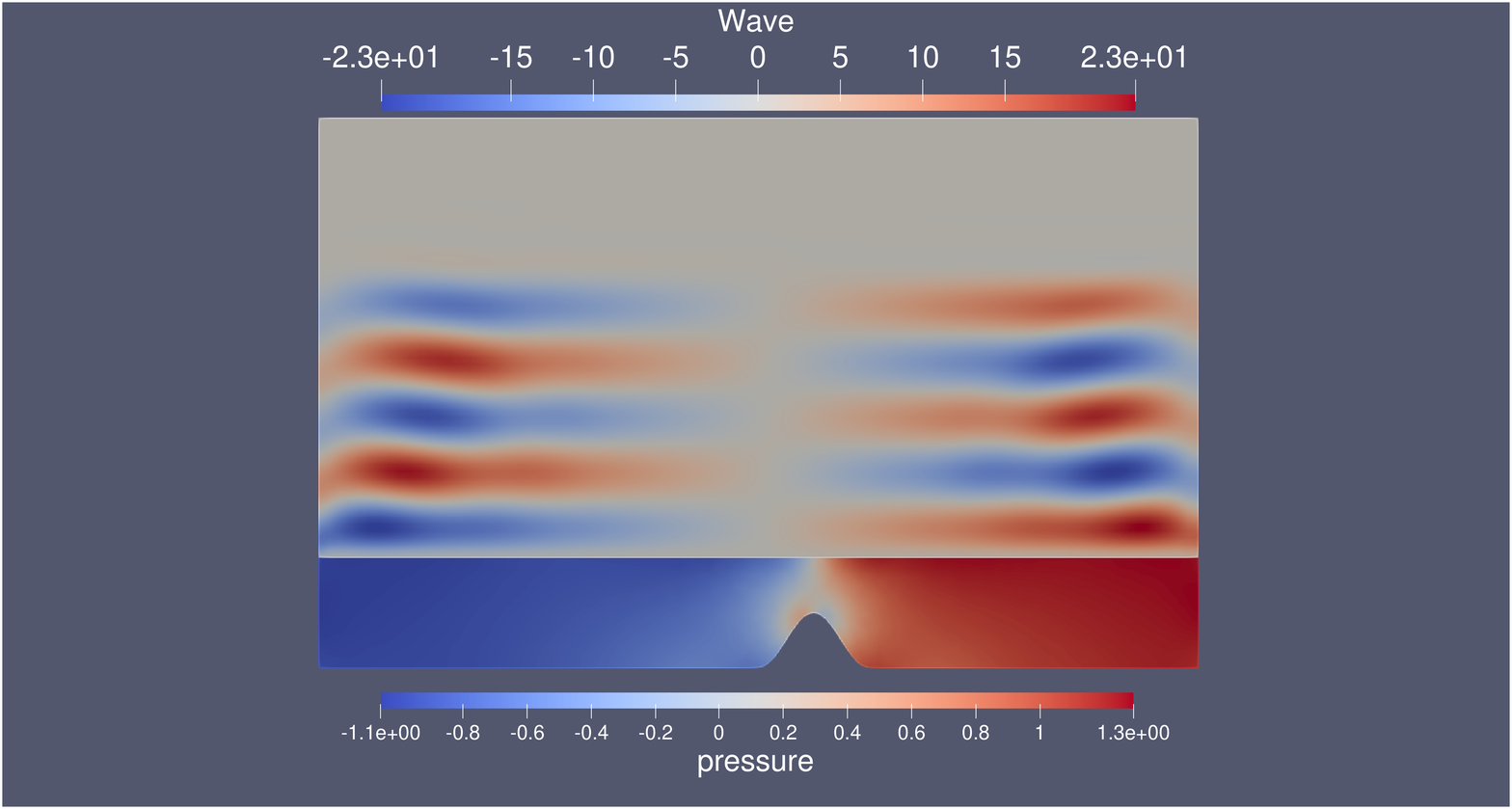}}
\subfloat[]{\includegraphics*[bb = 275 0 1240 850,scale=0.17]{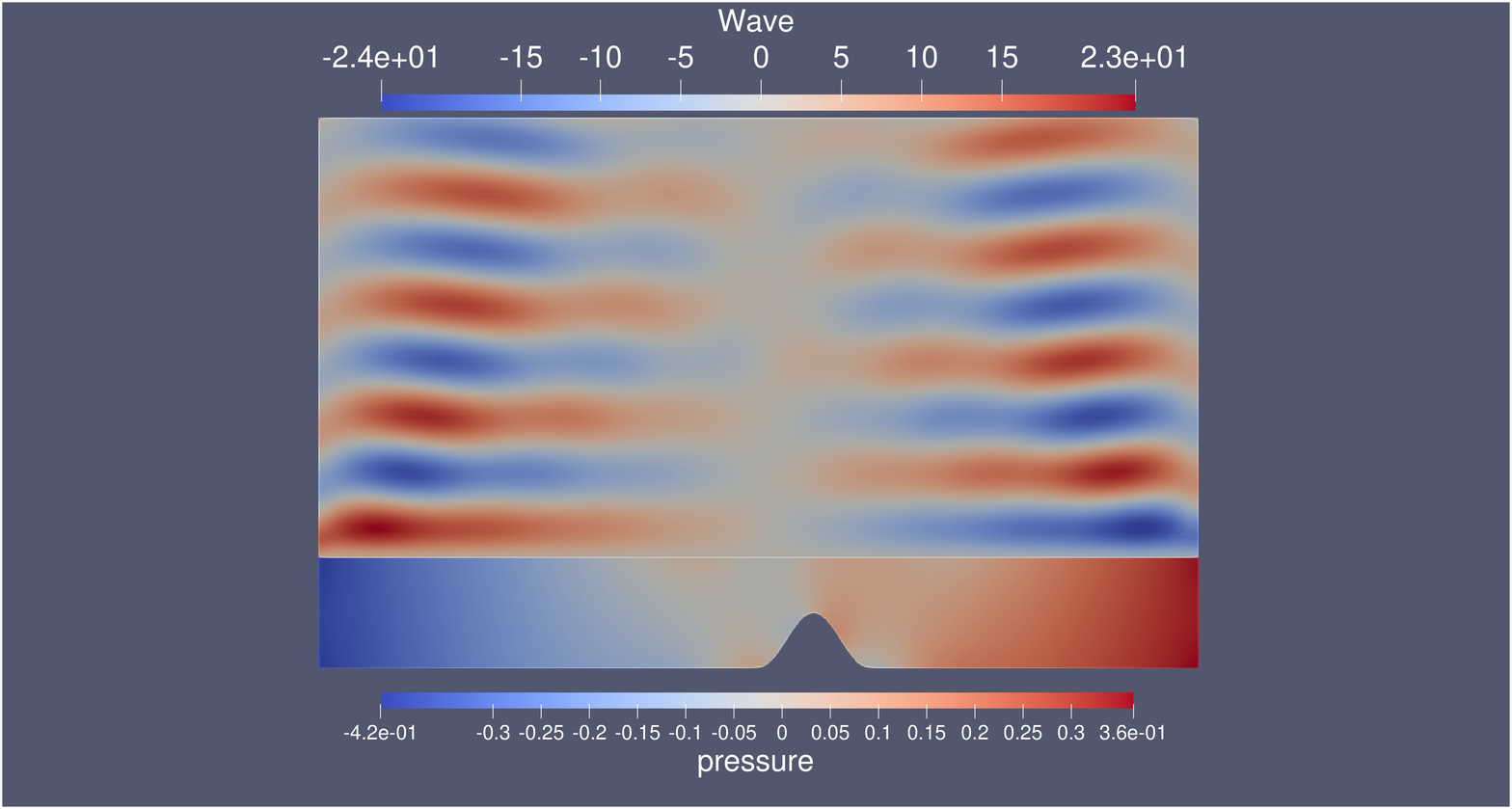}}
\caption{Norm of the pressure and the wave at $t=2.5$ (a), $t=5$ (b)}
\label{fig:Sim pres-onda}
\end{figure}

Figures \ref{fig:Sim vel-onda} and \ref{fig:Sim pres-onda} are snapshots of the simulations obtained from the velocity
field, the pressure and the acoustic wave:

\begin{figure}[hbt!]
\centering{}
\subfloat[]{\includegraphics*[bb = 300 0 1270 850,scale=0.17]{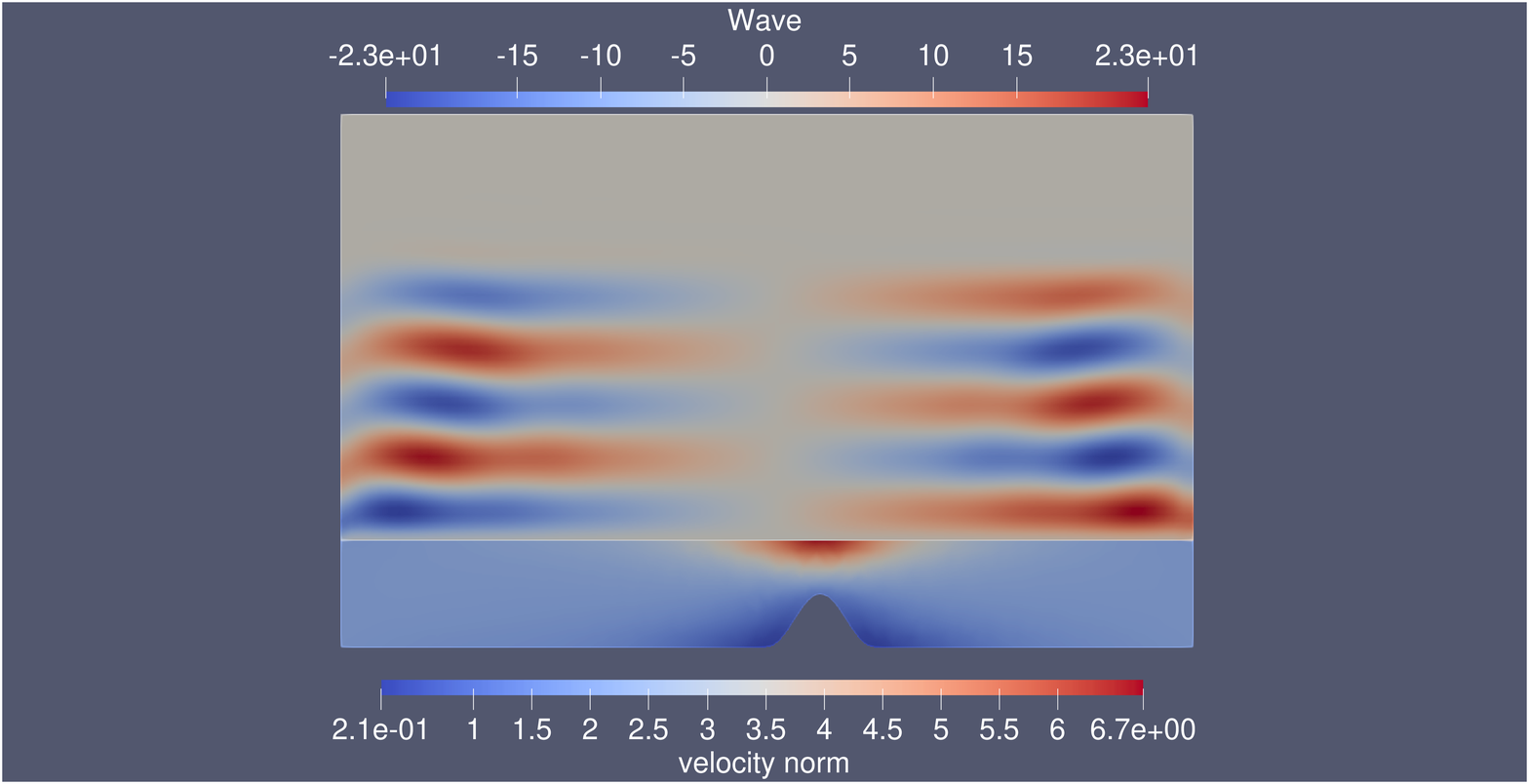}}
\subfloat[]{\includegraphics*[bb = 275 0 1240 850,scale=0.17]{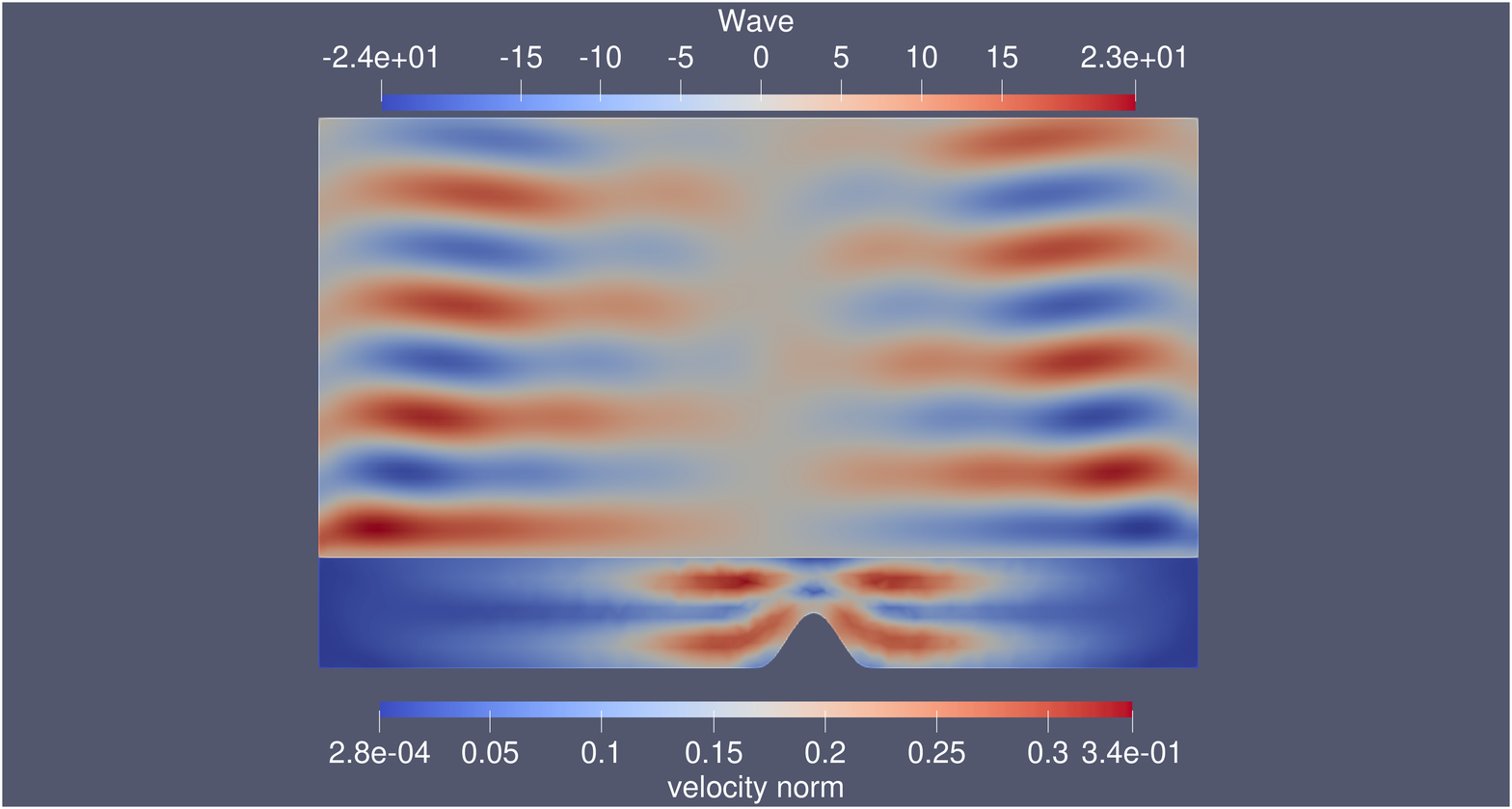}}
\caption{Norm of the velocity and the wave at $t=2.5$ (a), $t=5$ (b)}
\label{fig:Sim vel-onda}
\end{figure}

\section{An optimization process for the inverse problem}

As a first approximation and since our main goal is to determine the location, the deep and the size of the obstruction, we use the following class of obstacles to see if the inverse problem is feasible.

\begin{figure}[htbp]
	\begin{centering}
		\includegraphics*[bb = 0 325 2000 1500,scale=0.13]{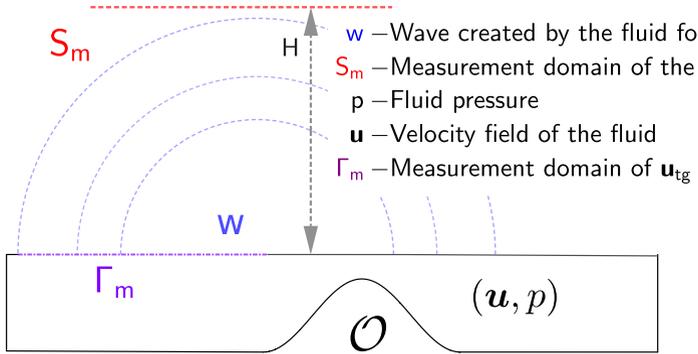} 
		\par\end{centering}
		\caption{\label{fig:mediciones} Example of measurement domains for the numerical inverse problems.}
	\end{figure}

Let ${\mathcal O}^{*}={\mathcal O}(\theta_{1}^{*},\theta_{2}^{*},\theta_{3}^{*})$ the obstruction that
we want to identify, and 
let $w_{m}(.)=w\left(.;{\mathcal O}(\theta_{1}^{*},\theta_{2}^{*},\theta_{3}^{*})\right)$
{\small{}be the measurements of the sounds waves generated
by the Cauchy tensor resulting from the Stokes flow}, (see systems (\ref{eq:Stokes-navier-slip}) 
--(\ref{eq:Onda})).
As stated previously, we are interested in reconstructing the obstruction
${\mathcal O}$ using measurements of the acoustic wave $w$ in the observable
set $S_{m}=[k_1,k_2]\times\{H\} \subset \partial S$ and measurements of the tangential velocity 
in $\Gamma_m$ a relative open set of the fluid domain border, see Figure \ref{fig:mediciones}.

We parameterize the obstruction as in the
previous subsection, that is:

\[
\partial{\mathcal O}=\begin{cases}
x(s)=\theta_{1}+s & s\in[0,\theta_2],\\
y(s)=\frac{\theta_{3}}{2}\left(1.0-cos(\frac{2\pi s}{\theta_{2}})\right) & s\in[0,\theta_2].
\end{cases}
\]

Using this simplification, the identification problem can be understood
as recovering the parameters $\theta_1,\theta_2,\theta_3$ that characterize the obstruction
${\mathcal O}$. In addition since the obstruction does not modify the
domain $S$ of the wave equation, it is possible to decompose the inverse
problem into two inverse problems as follows below.

\subsection{Inverse problem for the wave equation} \label{wave_subsection}
Using the assumption that $\partial{\mathcal {O}}\bigcap \ovl{\Gamma_{wall}^{up}}=\emptyset$ and 
equation (\ref{eq:Onda}) it can be observe that the measurement obtain at the set $S_{m}$ depends on 
the unknowns values of the boundary term defined on $\Gamma_{wall}^{up}\times(0,T)$.
Therefore, to estimate these unknown boundary values, we decide to minimize the following functional:
\begin{equation}\label{funcional_onda}
J_{1}(f)=\int_{S_{m}\times(0,T)}\norm{w(\cdot;f)-w_{m}}{2}^{2}\ dS\ dt,
\end{equation}
where $\|.\|_2$ is the Euclidean norm and $w_m$ are the wave measurements (see Figure \ref{fig:Representacion onda w}) and $w(.;f)\in C([0,T],H^{1}(S))$ is the solution to:

\begin{equation*}
\begin{cases}
w_{tt}-c^{2}\Delta{w}=0 & \mbox{ in } S\times(0,T),\\
w=f & \mbox{ on }\left([0,L] \mbox{ in } \mathbb{\times}\{D\}\right)\times(0,T)=\Gamma_{wall}^{up}\times(0,T),\\
\frac{\partial w}{\partial t}+c\frac{\partial w}{\partial n}=0 & \mbox{ on } \left(\partial S\backslash\left([0,L]\mathbb{\times}\{D\}\right)\right)\times(0,T),\\
w_{t}(\bb{x},0)=0 & \mbox{ in }S,\\
w(\bb{x},0)=0 & \mbox{ in }S.
\end{cases}
\end{equation*}

\noindent
Since it is well known that the wave equation has a finite propagation speed, we can only recover the unknown boundary datum at most within the time interval $(0,T-t_c)$, where $t_c$ represents the wave travel time between the boundary datum site and the measurement site. In our case, due to our simple geometry, the constant $t_{c}$ can be characterized by the following formula:
\begin{equation}
t_{c}=\frac{\sup_{x\in\Gamma_{wall}^{up}}d(x,S_{m})}{c} \ ,
\end{equation}
where $d(\bb{x},S_{m}):=inf\{d\in\mathbb{R}^{+}: \, d=\norm{\bb{x}-\bb{y}}{2},y\in\Gamma_{wall}^{up}\}$ 
and $c$ is the propagation speed.

\noindent
In summary, the first inverse problem corresponds to:
\begin{equation}\label{eq:Primer funcional inverso}
\mbox{Find }\hat{f}\in F_{ad}\mbox{ such that: }J_{1}(\hat{f})=\min_{f\in F_{ad}}J_{1}(f) \ , 
\end{equation}

\noindent
where $F_{ad}=\left\{ L^{2}\left(\Gamma_{wall}^{up}\times(0,T)\right): \,f=0\,\mbox{ on }\,\left(\Gamma_{wall}^{up}\times(t_{c},T)\right)\right\} $ is the admissible space of the boundary data.

\subsection{Obstacle inverse problem}
Consider $\hat{f}$ solution of (\ref{eq:Primer funcional inverso}), recalling that the functional (\ref{funcional_onda}) is a quadratic form and therefore convex, it is not unreasonable to assume that $\hat{f}$ is an estimation of the normal component of the Cauchy tensor $\sigma(\bb{u},p)\bb{n}\cdot\bb{n}$ on $\Gamma_{wall}^{up}\times(0,T)$. We also consider measurements of the tangential component of fluid $\bb{v_{m,\tau}}$ in a set $\Gamma_m \times(0,T) \subset \Gamma_{wall}^{up}\times(0,T)$. Thus,  we decide to minimize the following functional:

\begin{align}
\begin{split}
J_{2}(\theta) & =\int_{\Gamma_{wall}^{up}\times(0,T-t_{c})}\norm{\left(\sigma\left(\bb{u}(.;\Omega_{{\mathcal O}(\theta)}),p(.;\Omega_{{\mathcal O}(\theta)})\right)\bb{n}\right)\cdot\bb{n}-\hat{f}\,}{2}^{2}\ dS\ dt \\
&  + \int_{\Gamma_m \times (0,T)} \norm{\bb{v_{m,\tau}}-\bb{u}(.;\Omega_{{\mathcal O}(\theta)})\cdot\bb{\tau}}{2}^{2} \ dS\ dt \ ,
\end{split}
\label{funcional_fluido}
\end{align}

where $\theta \in \Theta\subset\mathbb{R}^{3}$ is the set of possible parameters
that represent an obstruction, $\Omega_{{\mathcal O}(\theta)} = \Omega \backslash \mathcal{O(\theta)}$  and $\bb{u}(.;\Omega_{{\mathcal O}(\theta)}),p(.;\Omega_{{\mathcal O}(\theta)})\in H^{1}(\Omega_{{\mathcal O}(\theta)})\times L^{2}(\Omega_{{\mathcal O}(\theta)})$
are solutions to:
\[
\begin{cases}
\bb{u}_{t}-div(\sigma(u,p))=0 & \mbox{ in } \Omega_{{\mathcal O}(\theta)}\times(0,T),\\
div(\bb{u})=0 & \mbox{ in } \Omega_{{\mathcal O}(\theta)}\times(0,T),\\
\bb{u}=g_{in} & \mbox{ on } \Gamma_{inlet}\times(0,T),\\
\bb{u}=g_{out} & \mbox{ on } \Gamma_{out}\times(0,T),\\
\bb{u}\cdot \bb{n}=0, \ (\sigma(\bb{u},p)\bb{n})_{tg}=0 & \mbox{ on } \Gamma_{wall}^{\mathcal{O}(\theta)}\times(0,T),\\
\bb{u}(.,0)=\bb{u}_{0} & \mbox{ in } \Omega_{{\mathcal O}(\theta)} \ .
\end{cases}
\]

It is important to remark that the use of the tangential velocity at the upper wall, which is included in the functional (\ref{funcional_fluido}) is necessary due to the identification result, Theorem \ref{identification_theo}. At a first glance this might seen as a contradiction with respect to our initial goal since no information of the tangential velocity is carry out by the wave equation. Nevertheless, in the following section we also consider the case without the tangential velocity ($\Gamma_m = \emptyset$ in (\ref{funcional_fluido})) producing a reasonable result (see experiment 4, Table \ref{Mcmc-table-1}). 

\subsection{Numerical results}
The optimization problem described in the previous section involves the solution of equations \eqref{eq:Stokes-navier-slip}-\eqref{eq:Onda} which are solve as it described in subsections (\ref{subsec:Direct problem and experimental setup})-(\ref{subsec:Numerical_direct_problem}). Nevertheless, to lower computational cost, we set $\Omega =[0,8] \times [0,1]$, $S=[0,8]\times[1,3]$, the final time $T=1$, the numerical time step as $\Delta t=1/500$ and the wave propagation speed as $c=\sqrt{30}$. Also, to increase the similarity of the simulated flow with respect to the blood, we set the fluid viscosity as $\mu=1/500$.

In order investigate the robustness of the optimization procedure we decide to perform 8 numerical experiments, with different obstructions geometries and different observable domains of the wave $S_m$ and of the tangential component of the velocity at the boundary $\Gamma_m$, see tables \ref{Mcmc-table-1}-\ref{Mcmc-table-3} for a summary. 

\begin{figure}[hbt!]
\centering
\includegraphics[trim={2.5cm 5.5cm 2.5cm 7.5cm},clip,width=0.8\textwidth]{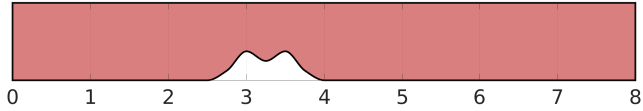}
\par
\caption{Original domain for experiments 7-8} \label{fig:Spline_shape}
\end{figure}

For the first 6 numerical experiments we generate synthetic data using obstructions $\mathcal{O}$ defined by the parameters $(\theta_1,\theta_2,\theta_3)$ (see equation (\ref{eq:Dominio Parametrico}) and Figure \ref{fig:Dominio Parametrico}) and for experiments 7-8 we generated the data using a different class of obstruction, defined by a cubic spline (see figure \ref{fig:Spline_shape}). In order to avoid the inverse crime we add a random error relative to the method order, given by a normal distribution of variance $\sigma=10^{-5}$ and mean $\mu=0$.

As state, to obtain an estimation of the normal component of the Cauchy tensor on the upper wall of the fluid domain, we minimize the functional (\ref{funcional_onda}). Since it is easy to see that the functional $J_1$ is convex, we decide to use the gradient method. 

\begin{figure}[hbt!]
\centering{}
\subfloat[]{\centering{}\includegraphics[width=0.5\textwidth]{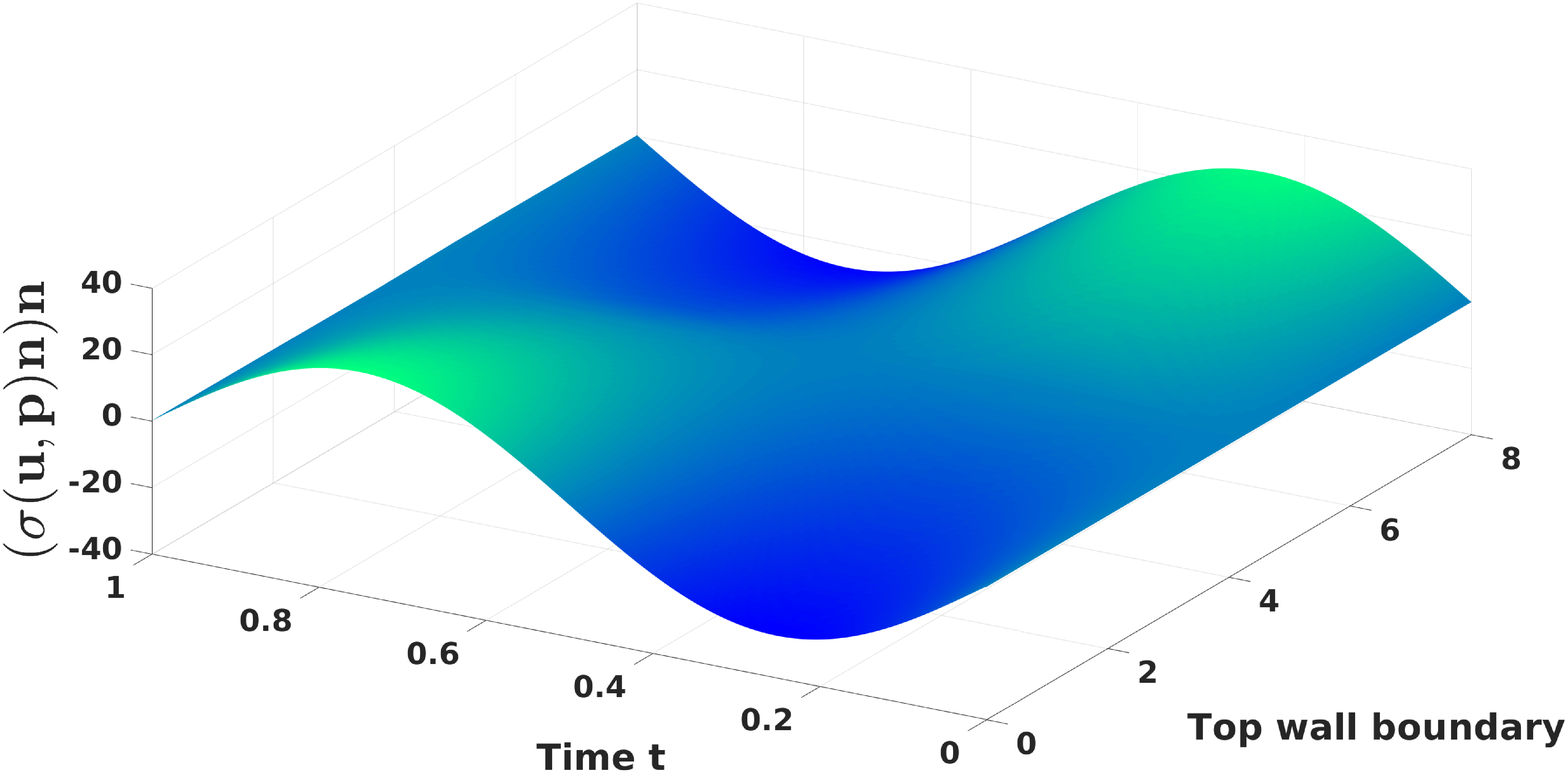}}
\subfloat[]{\centering{}\includegraphics[width=0.5\textwidth]{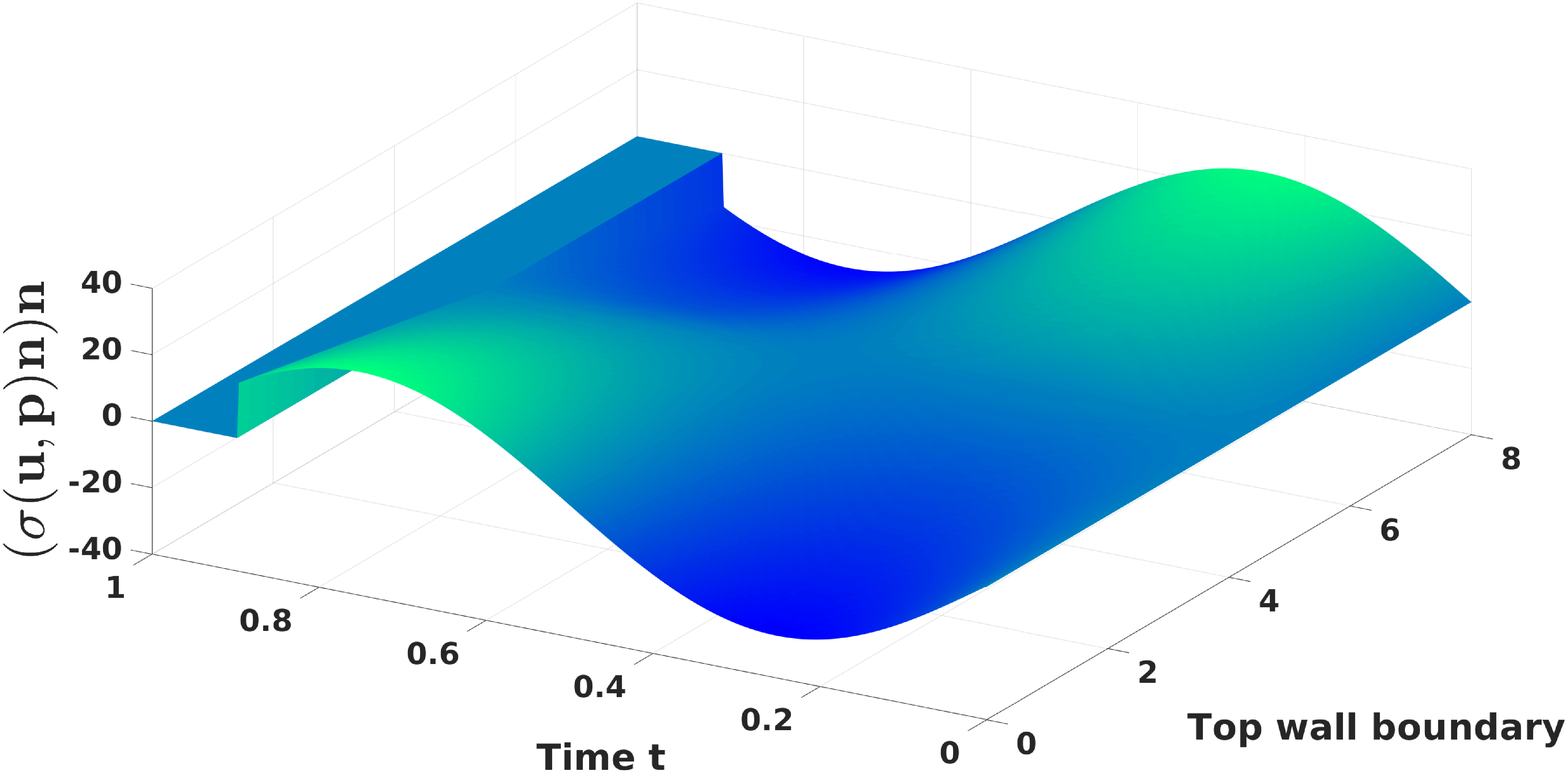}}\caption{a) Original wave boundary datum, b)
Reconstructed boundary datum.
\label{fig:Wave_inverse}
}
\end{figure}

Figure \ref{fig:Wave_inverse}, shows our achieved
reconstruction of the boundary data by the gradient descent method. Notice that we do not obtain a complete reconstruction of the boundary datum. This is a consequence of the fact that the wave has a time delay to arrive to the measure site as explained above.

\begin{figure}[hbt!]
\centering
\includegraphics[trim={2.5cm 5.5cm 2.5cm 7.5cm},clip,width=0.8\textwidth]{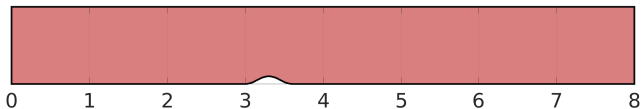}
\par
\caption{Initial guess shape for the MCMC method} \label{fig:Inicial_shape}
\end{figure}

On the other hand, for the geometry estimation of the obstruction, we decided to use an MCMC (Markov Chain Monte Carlo) method to minimize the functional (\ref{funcional_fluido}),  since as far as we know, no obstruction reconstruction has been done using this method. The reasons behind this choice are several; first we ignore the nature of this functional, it may be non-convex, second, taking into account that we deal with few parameters MCMC methods are effective in exploring the parameter space and furthermore they can be done in parallel thus making it possible to easily exploit the capability of multiple cpu cores, and also it has the advantage of providing an uncertainty range of the parameters, which is valuable in real application. It is worth pointing out that MCMC methods are not strongly dependent on the starting point, therefore, for all the following numerical experiments we use the following parameters as an initial guess: $\theta_{1}=3,\theta_{2}=0.6,\theta_3=0.1$ (see Figure \ref{fig:Inicial_shape}). 

We point out that the functional (\ref{funcional_fluido}) is a minimum square, thus the corresponding discretization is a sum of squares and therefore can be understood as a normal distribution of the error between the model and the data. This fact allows us to use common libraries of MCMC methods; in this case, we used the Matlab version presented in \cite{haario2006dram}. In our case, the a priori information given to the algorithm is that the parameters $(\theta_1,\theta_2,\theta_3)$ are uniformly distributed in $\Theta=[1, 6] \times [0.5, 1.5] \times [0.1, 0.9]$ for the first six numerical experiments and $\Theta=[1, 6] \times [0.5, 2.5] \times [0.1, 0.9]$ for the last two. This is because in the last two experiments (7-8) the obstacle size is larger (around 1.5).  For more information on these methods, we refer the reader to \cite{gelman2013bayesian}.

Tables \ref{Mcmc-table-1}-\ref{Mcmc-table-3} are a summary of the output of the MCMC algorithm. As can be seen, most experiments yield quite reasonable results and, in all cases, the original parameters stay within the standard deviation of the mean. Note that the tangential velocity information reduces the standard deviation of the estimate parameter, but acceptable results without this are obtained (Experiment 4).

\begin{figure}[hbt!]
\begin{centering}
\subfloat[]{\centering{}\includegraphics[width=0.5\textwidth]{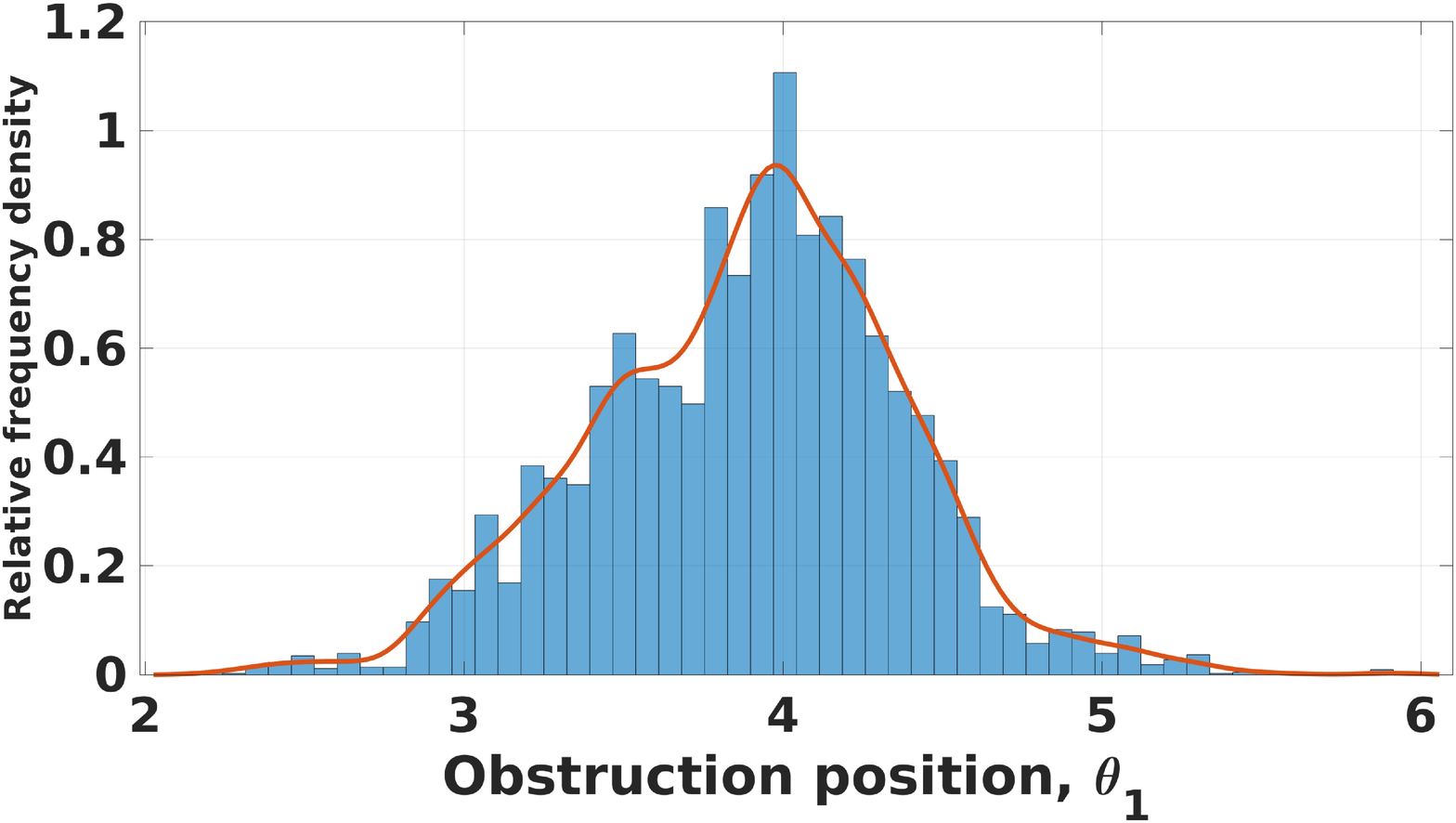}}
\subfloat[]{\centering{}\includegraphics[width=0.5\textwidth]{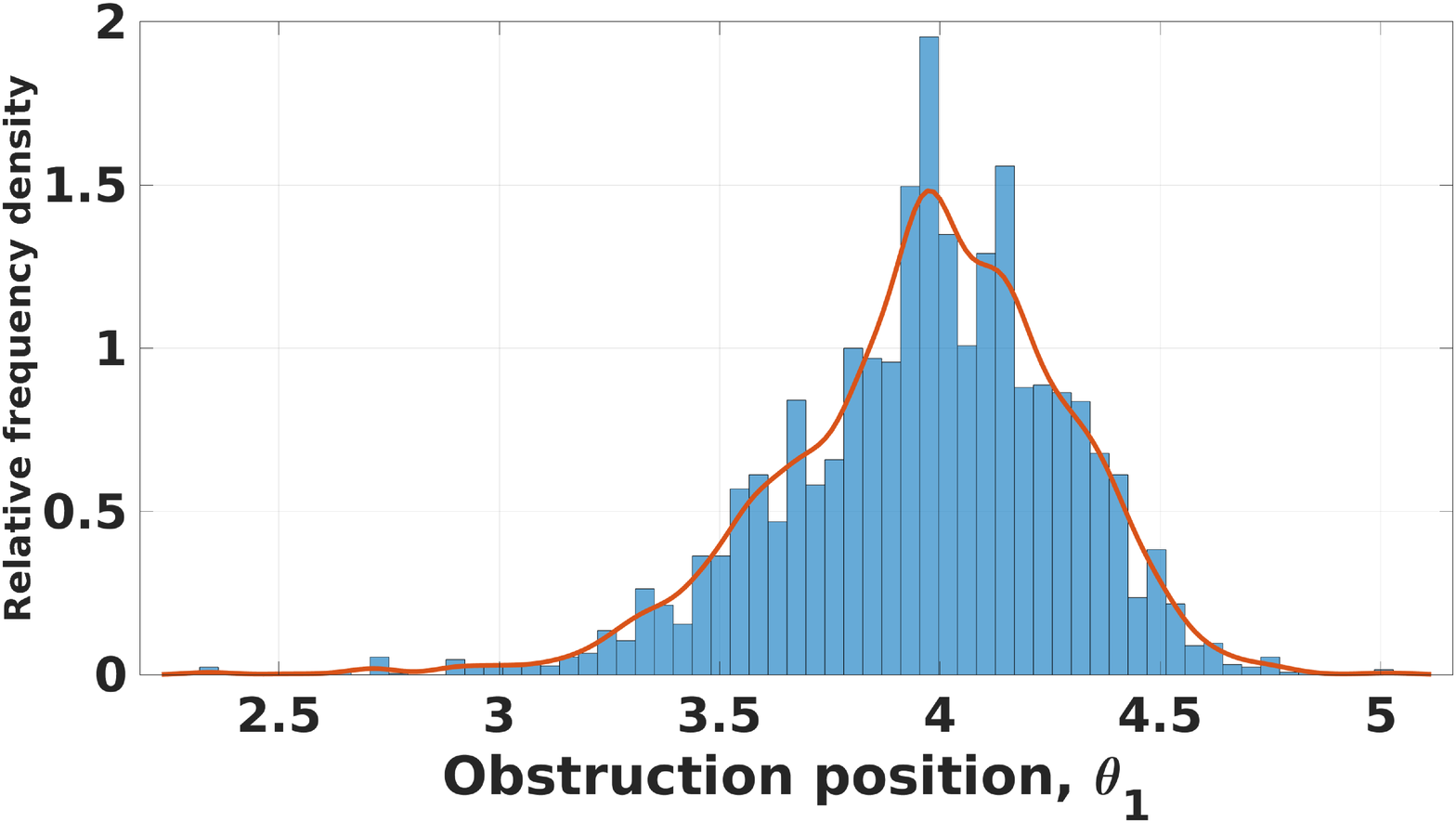}}

\subfloat[]{\centering{}\includegraphics[width=0.5\textwidth]{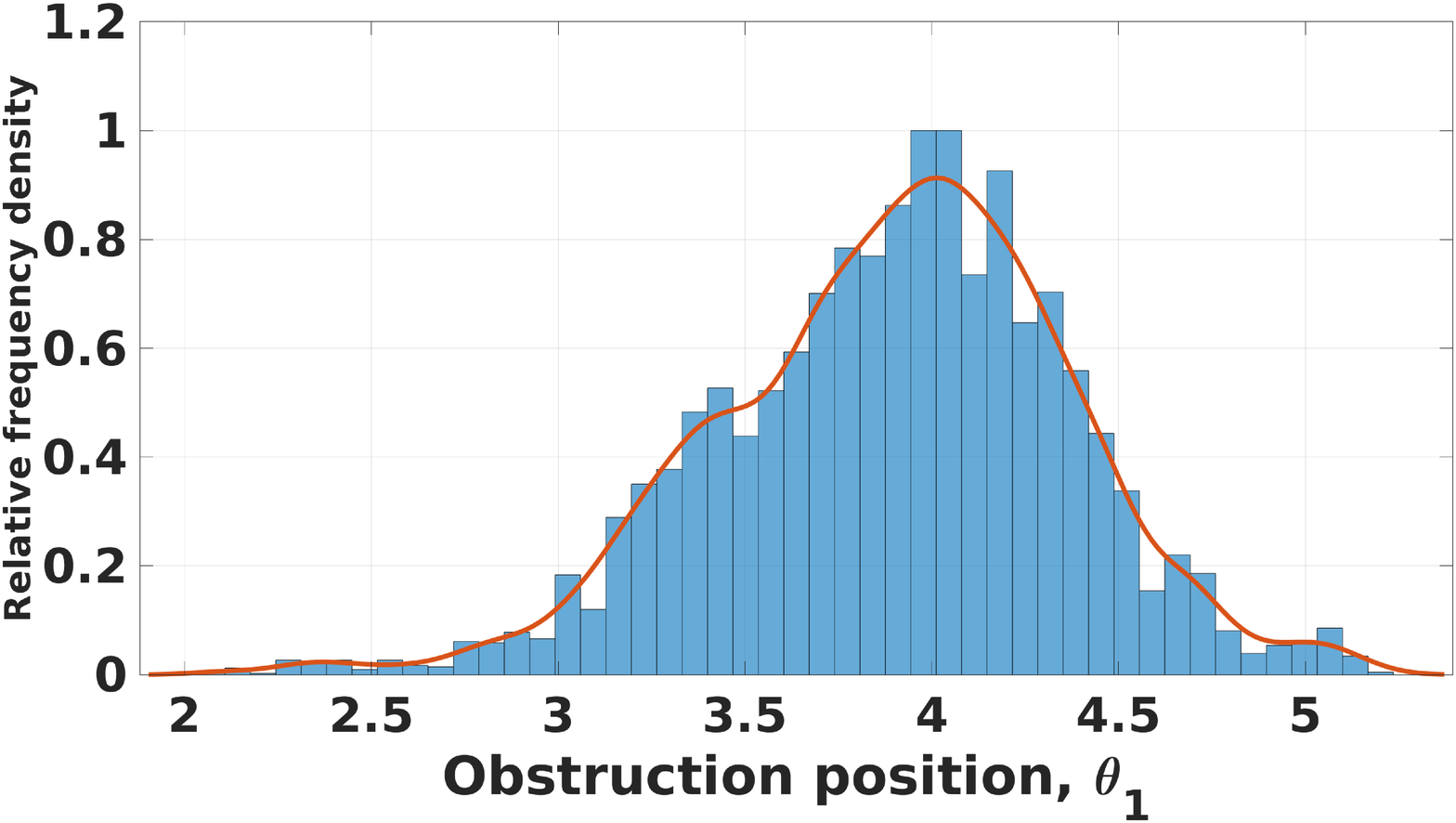}}
\subfloat[]{\centering{}\includegraphics[width=0.5\textwidth]{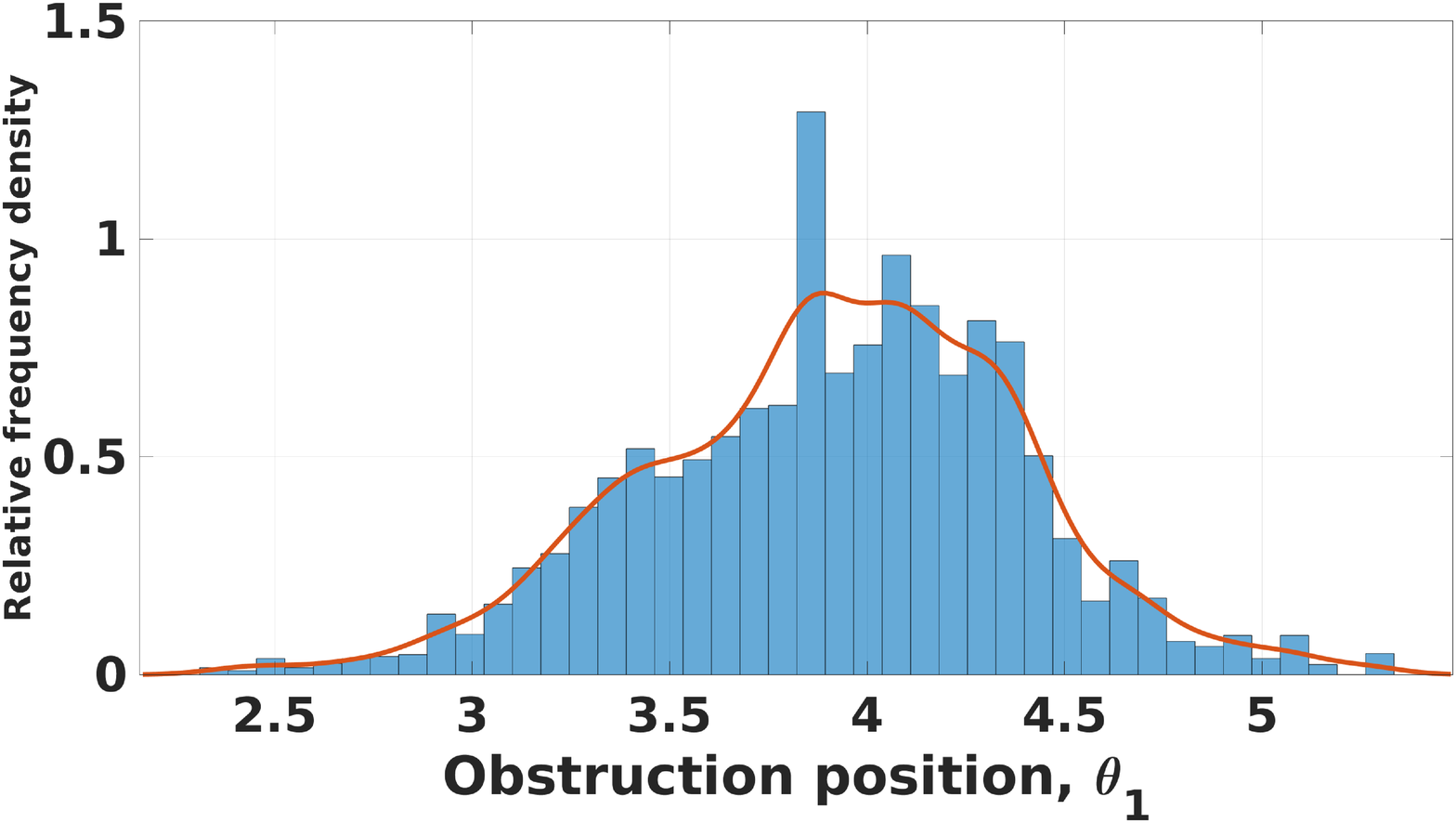}}

\subfloat[]{\centering{}\includegraphics[width=0.5\textwidth]{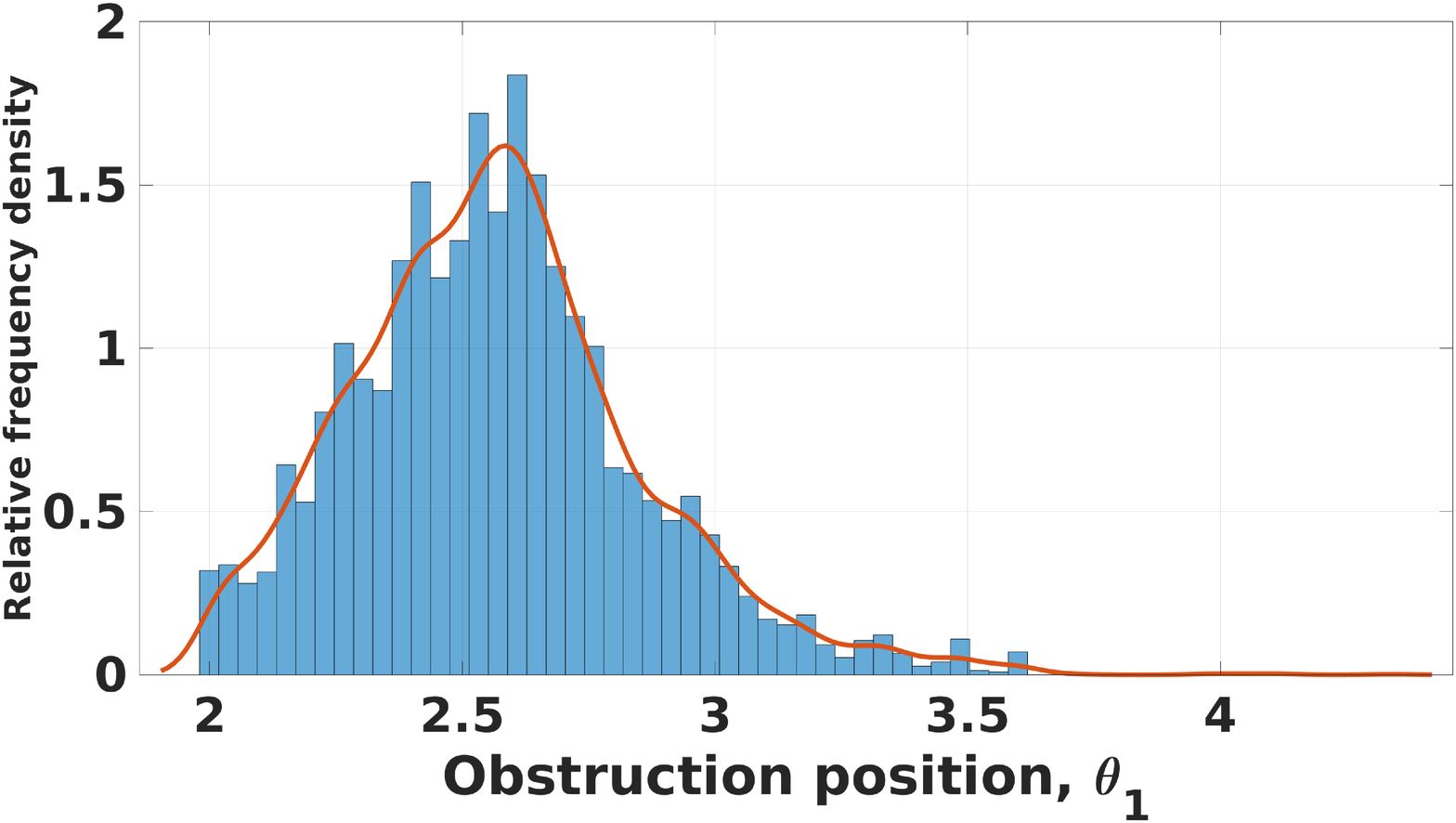}}
\subfloat[]{\centering{}\includegraphics[width=0.5\textwidth]{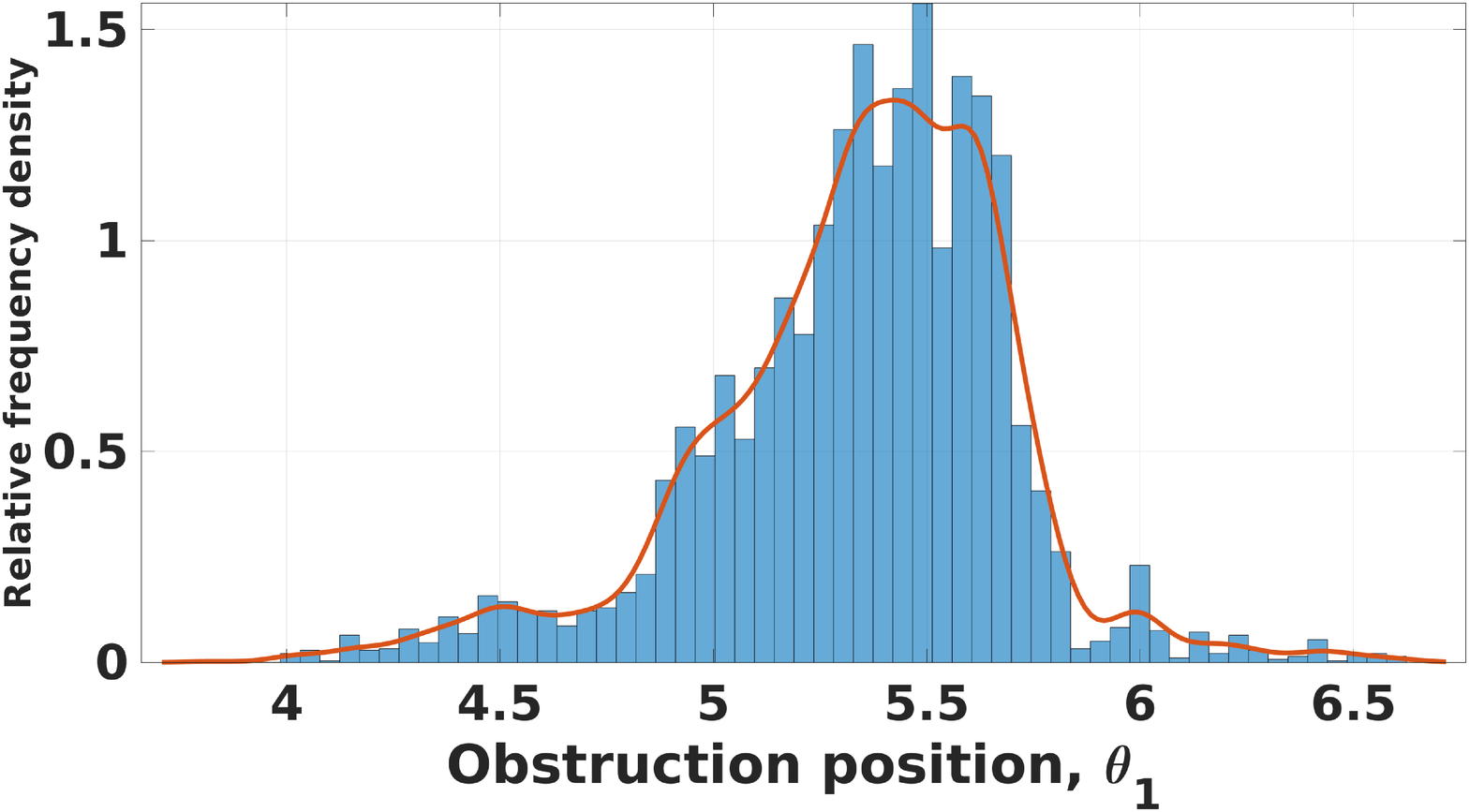}}

\subfloat[]{\centering{}\includegraphics[width=0.5\textwidth]{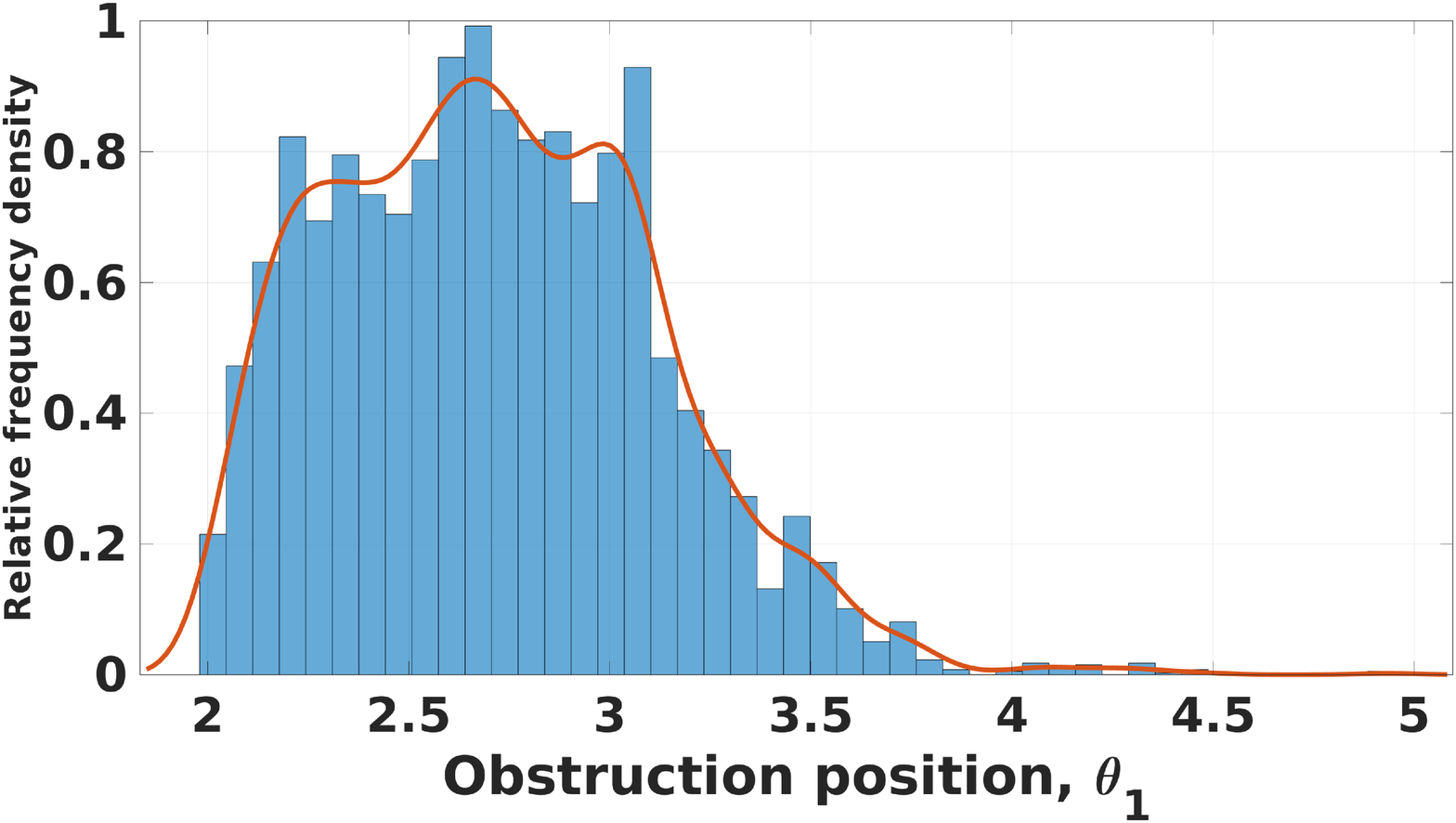}}
\subfloat[]{\centering{}\includegraphics[width=0.5\textwidth]{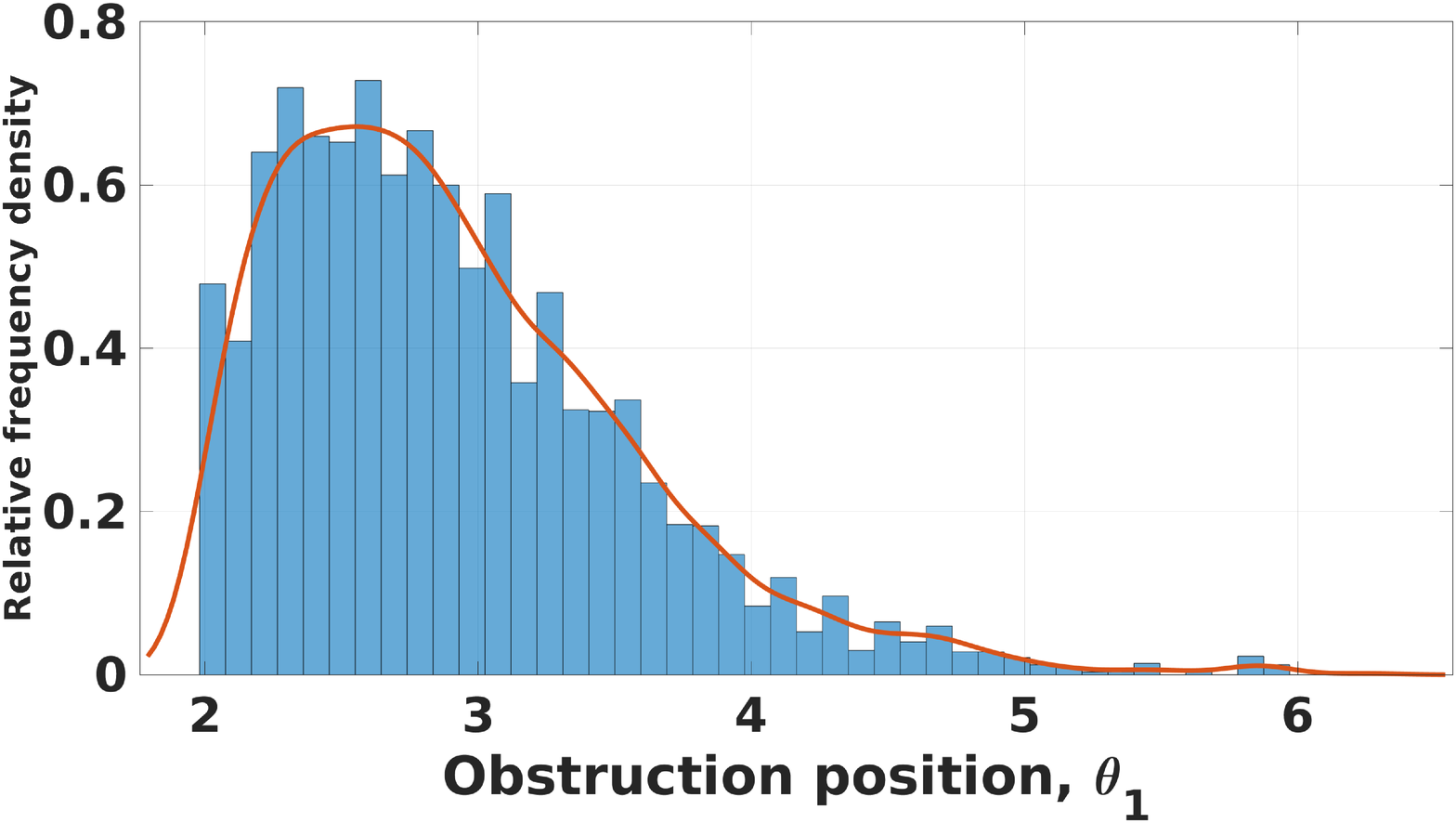}}
\par\end{centering}

\caption{Probability density function estimate of $\theta_1$ (Obstruction position) of experiments 1-8
(a-h)} \label{fig:theta_1_dist}
\end{figure}

\begin{figure}[hbt!]
\begin{centering}
\subfloat[]{\centering{}\includegraphics[width=0.5\textwidth]{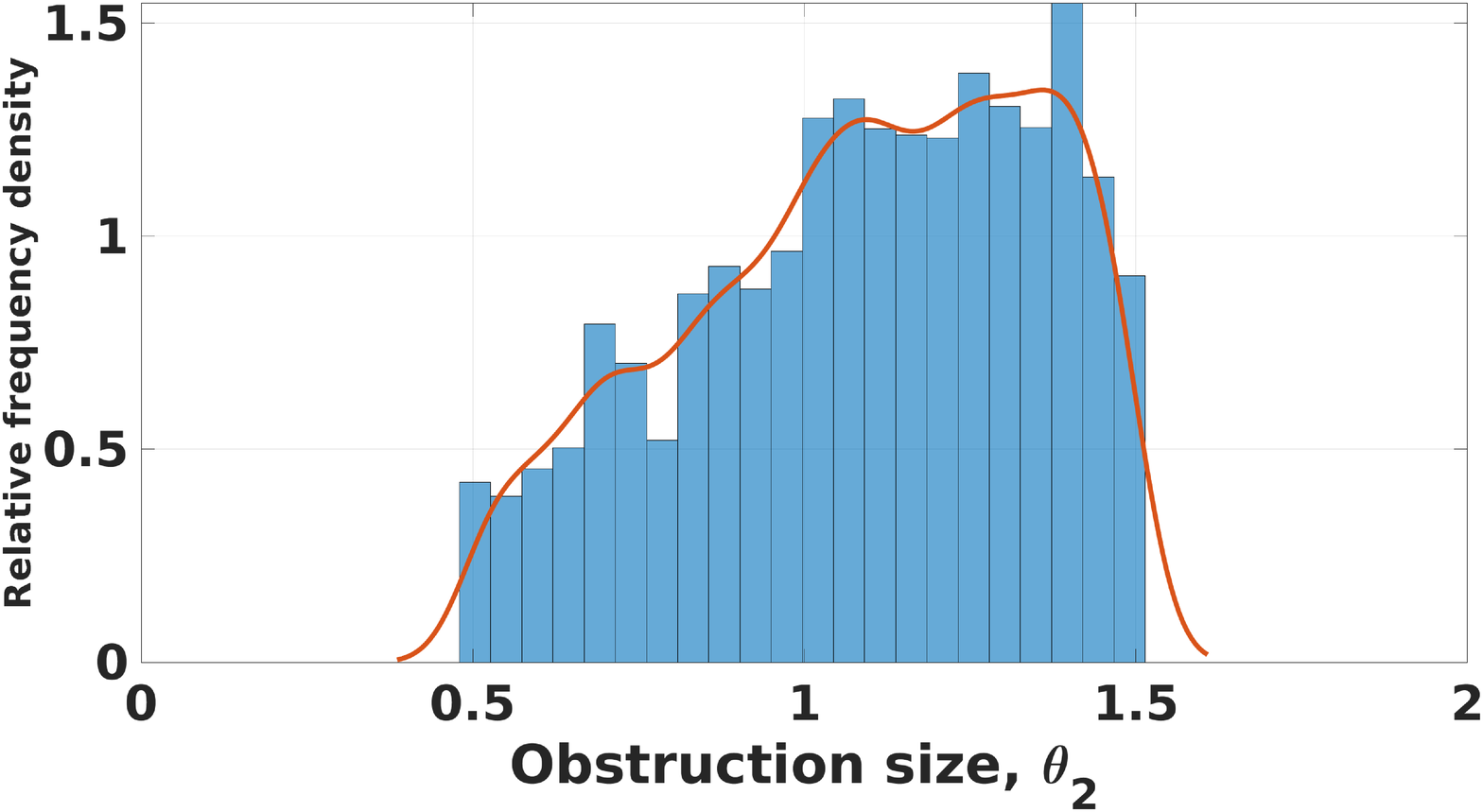}}
\subfloat[]{\centering{}\includegraphics[width=0.5\textwidth]{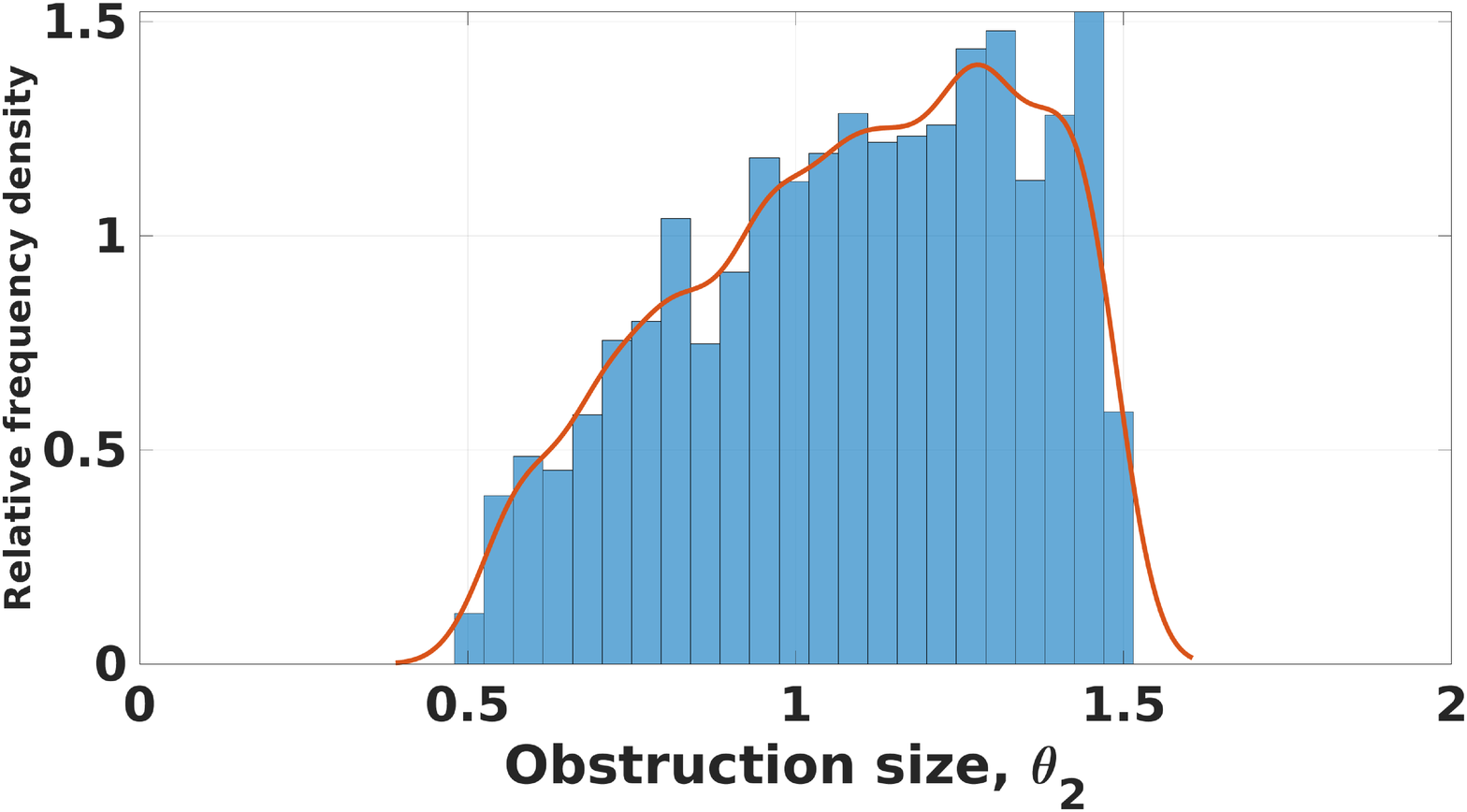}}

\subfloat[]{\centering{}\includegraphics[width=0.5\textwidth]{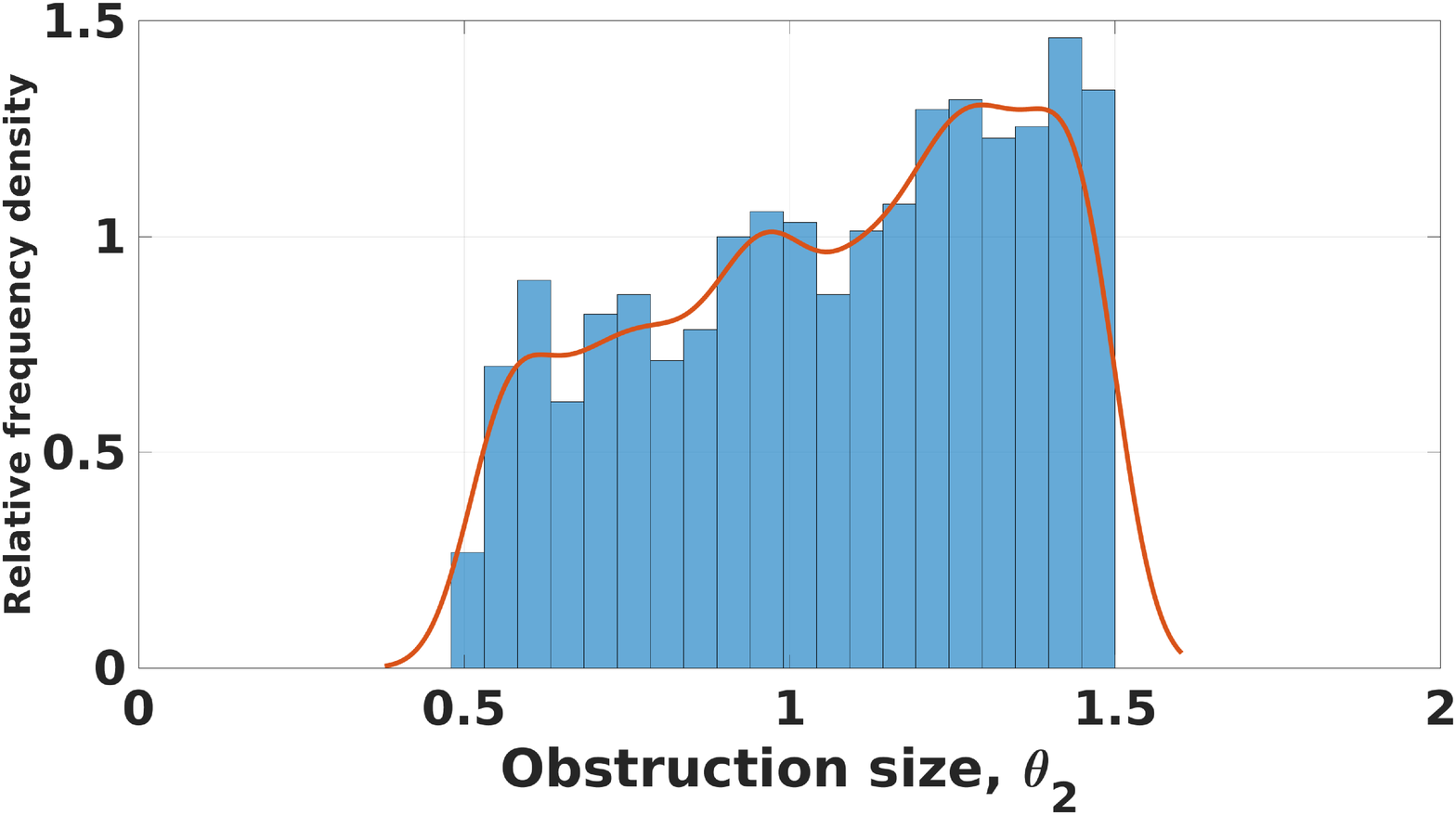}}
\subfloat[]{\centering{}\includegraphics[width=0.5\textwidth]{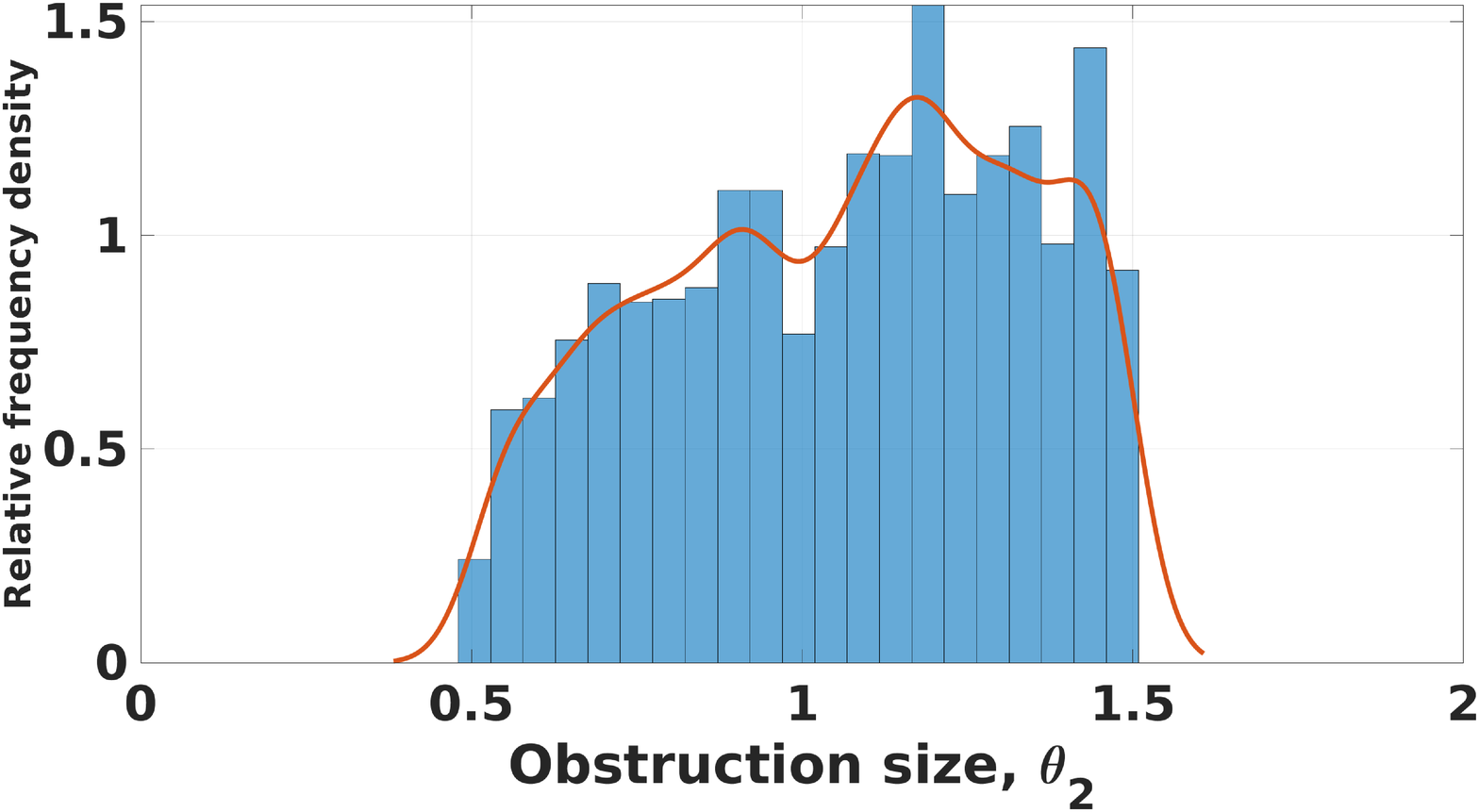}}

\subfloat[]{\centering{}\includegraphics[width=0.5\textwidth]{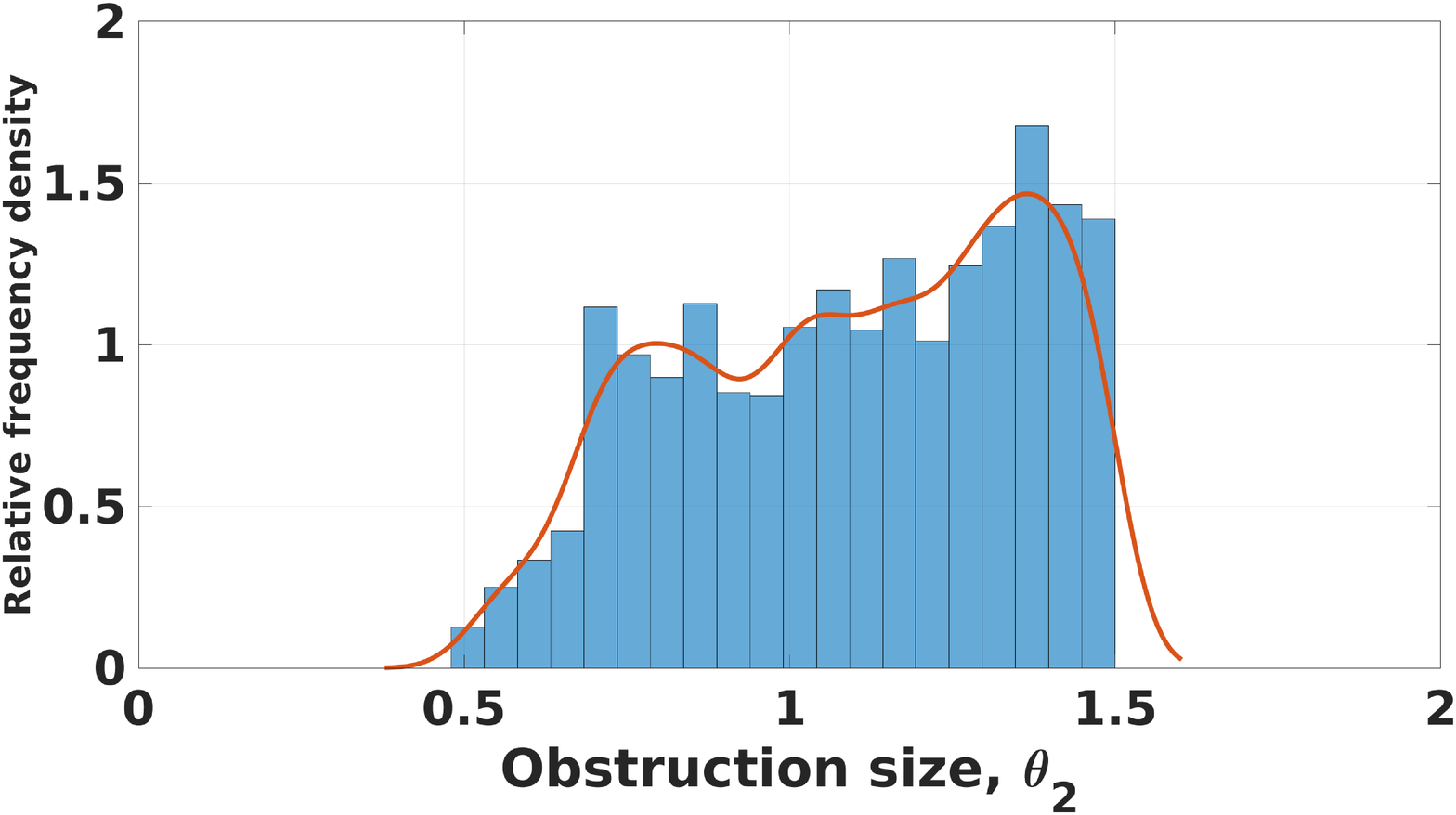}}
\subfloat[]{\centering{}\includegraphics[width=0.5\textwidth]{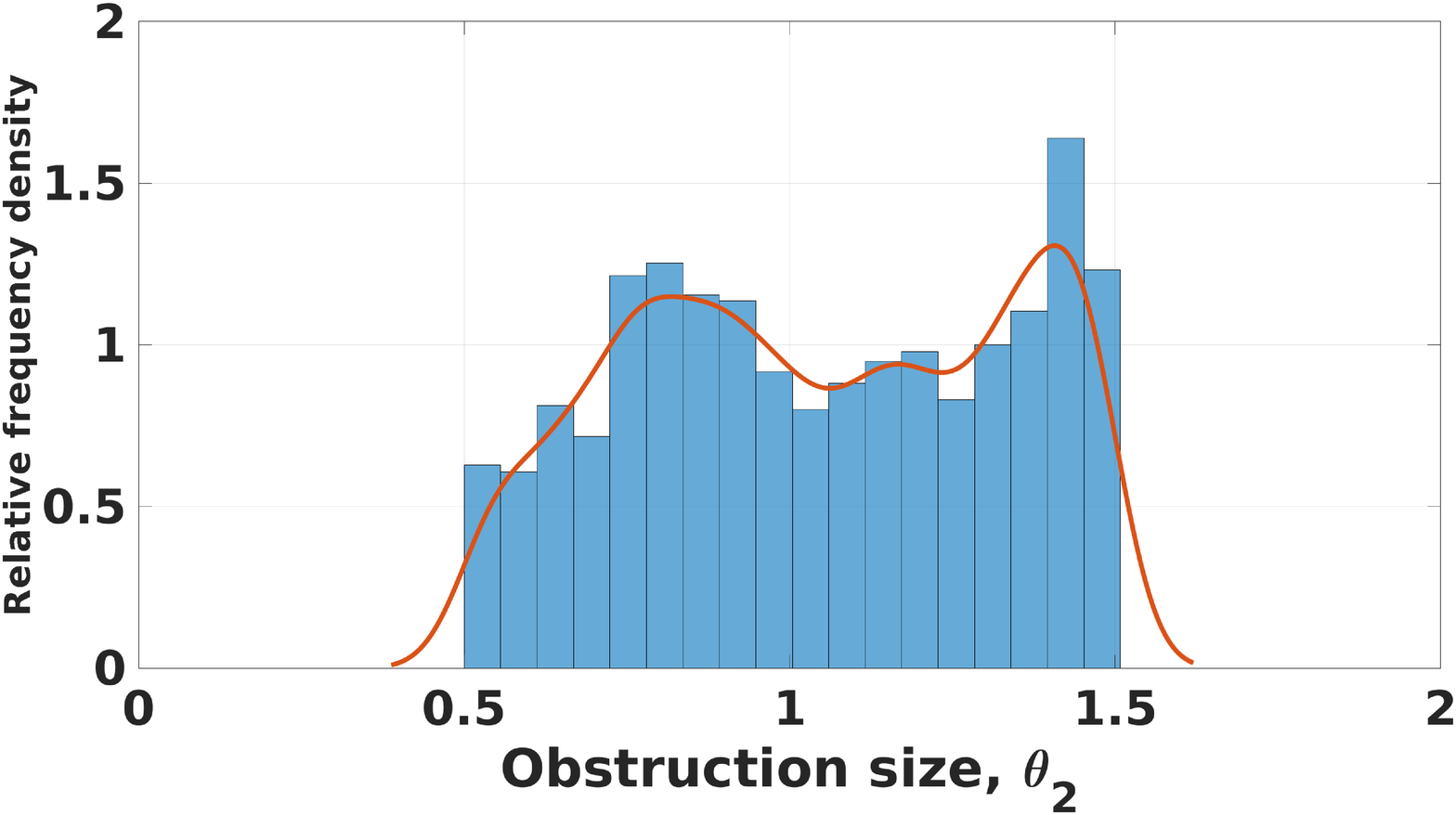}}

\subfloat[]{\centering{}\includegraphics[width=0.5\textwidth]{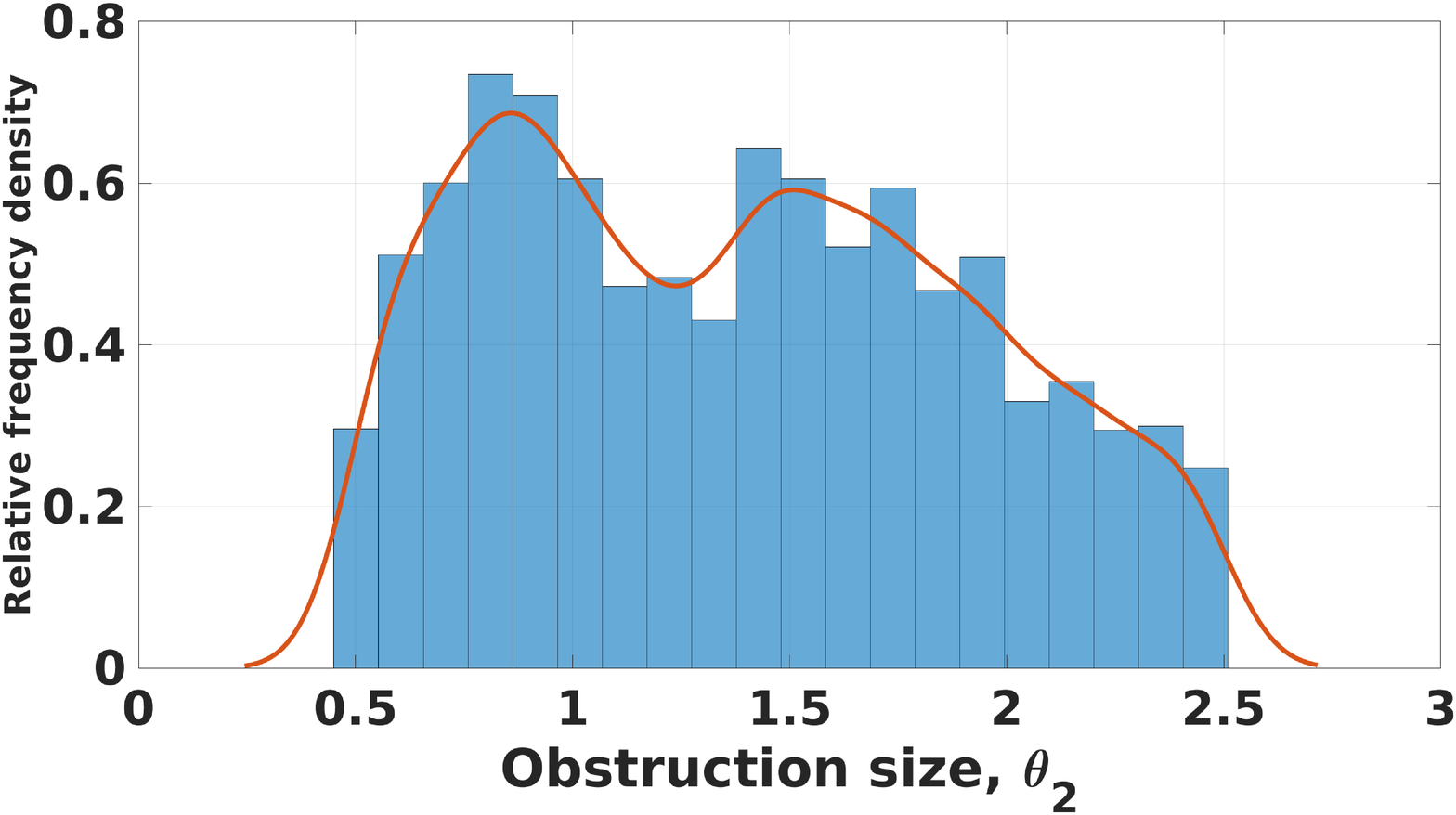}}
\subfloat[]{\centering{}\includegraphics[width=0.5\textwidth]{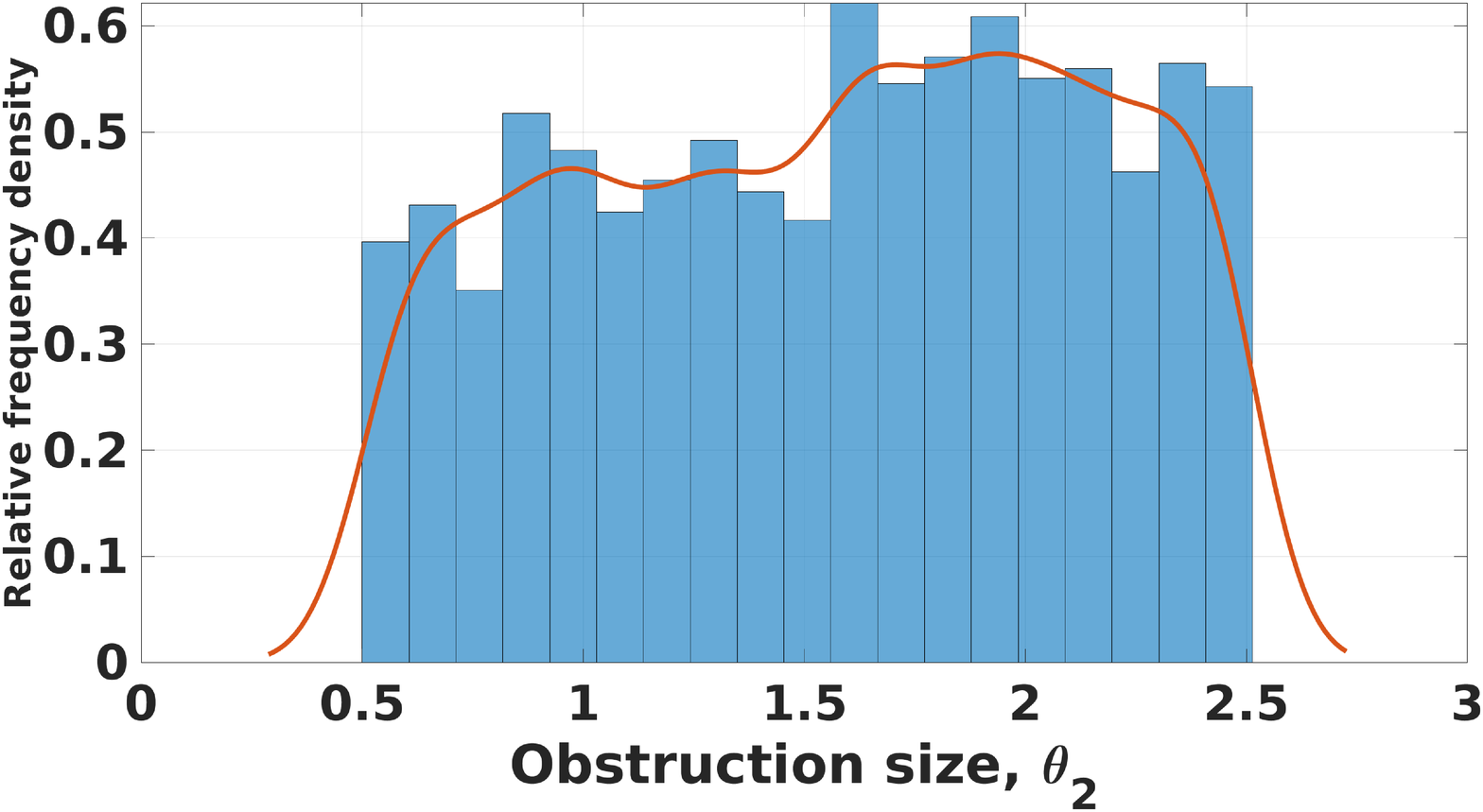}}
\par\end{centering}

\caption{Probability density function estimate of $\theta_2$ (Obstruction size) of experiments 1-8
(a-h)}
\label{fig:theta_2_dist}
\end{figure}

\begin{figure}[hbt!]
\begin{centering}
\subfloat[]{\centering{}\includegraphics[width=0.5\textwidth]{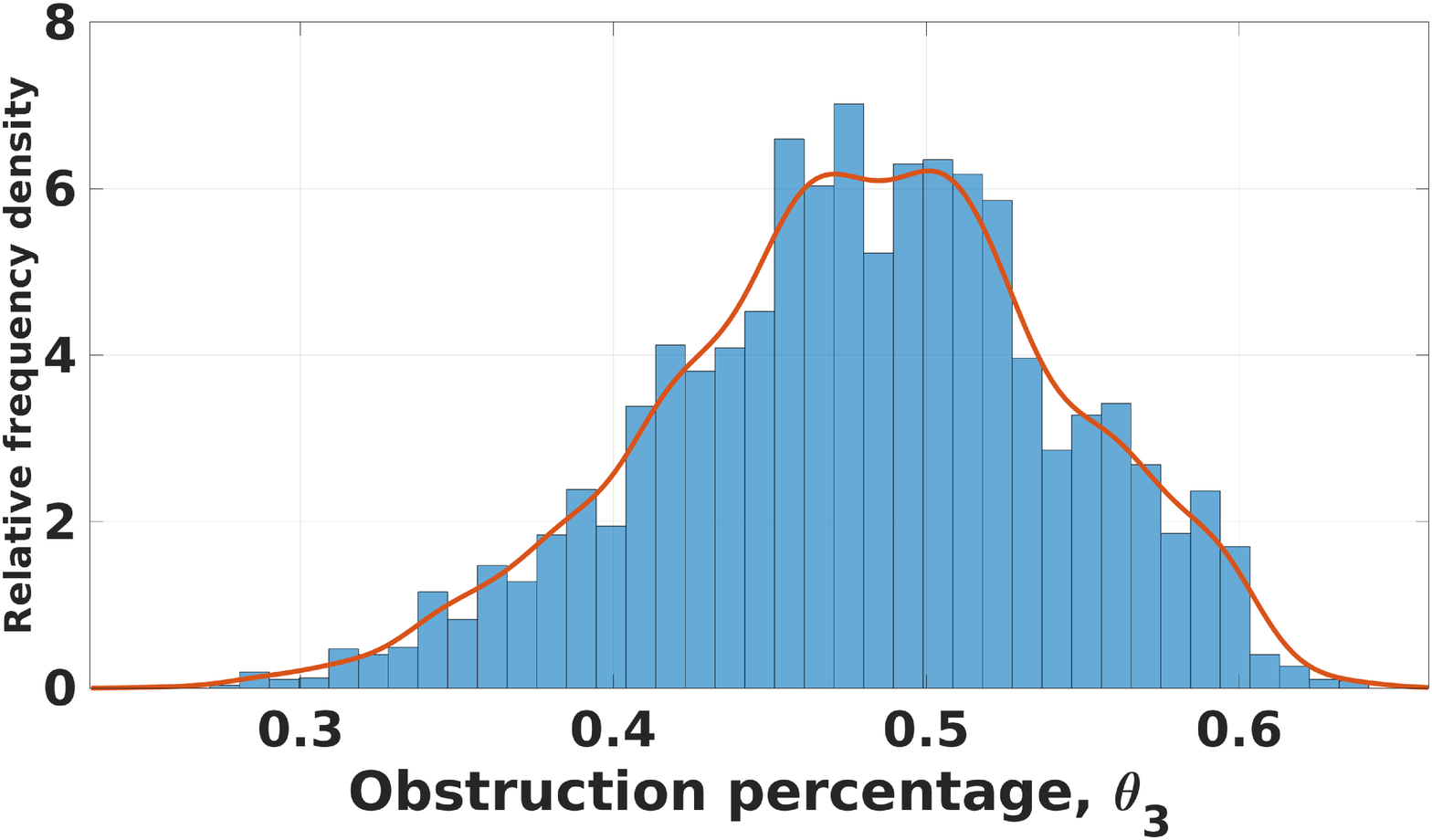}}
\subfloat[]{\centering{}\includegraphics[width=0.5\textwidth]{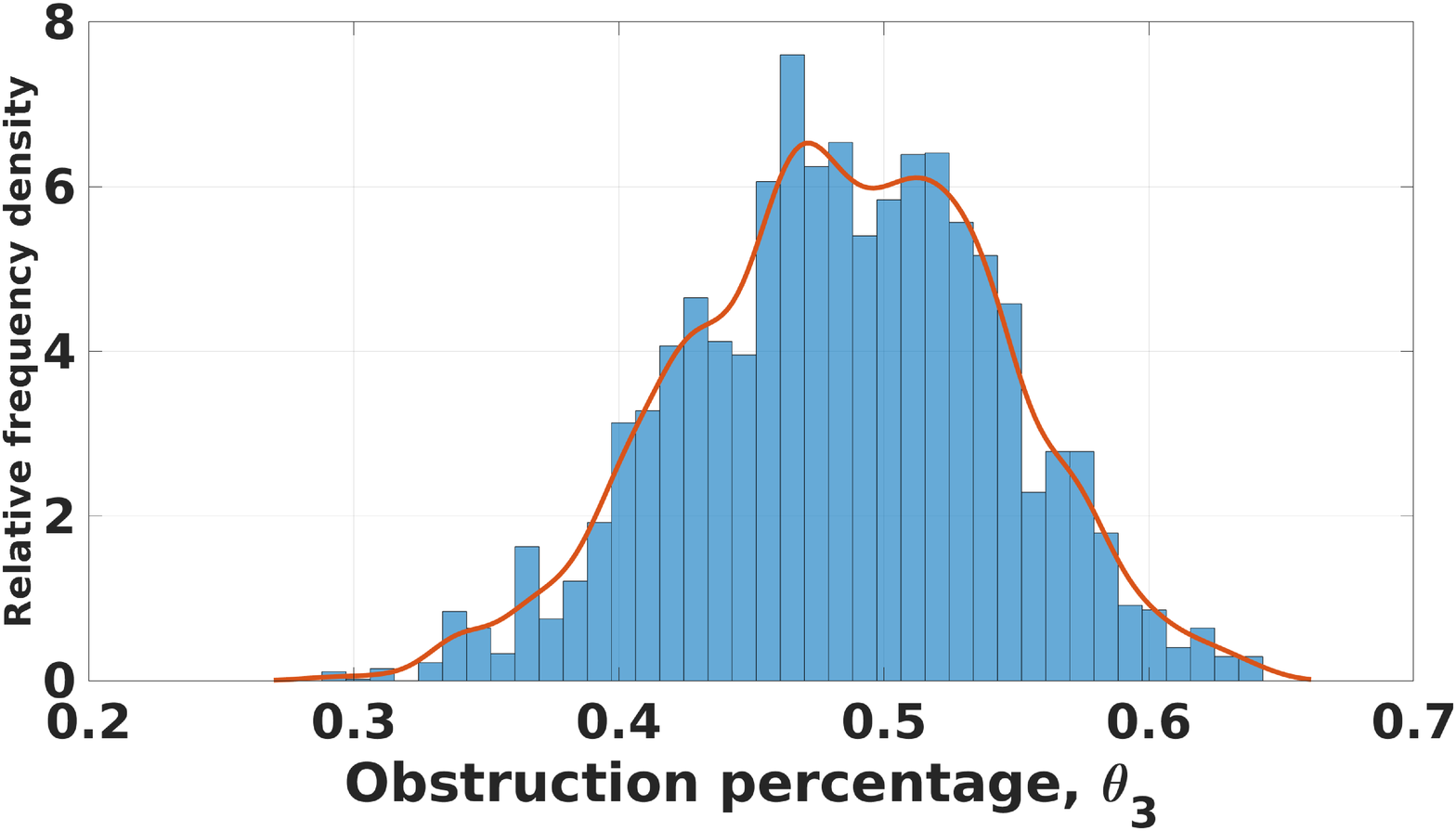}}

\subfloat[]{\centering{}\includegraphics[width=0.5\textwidth]{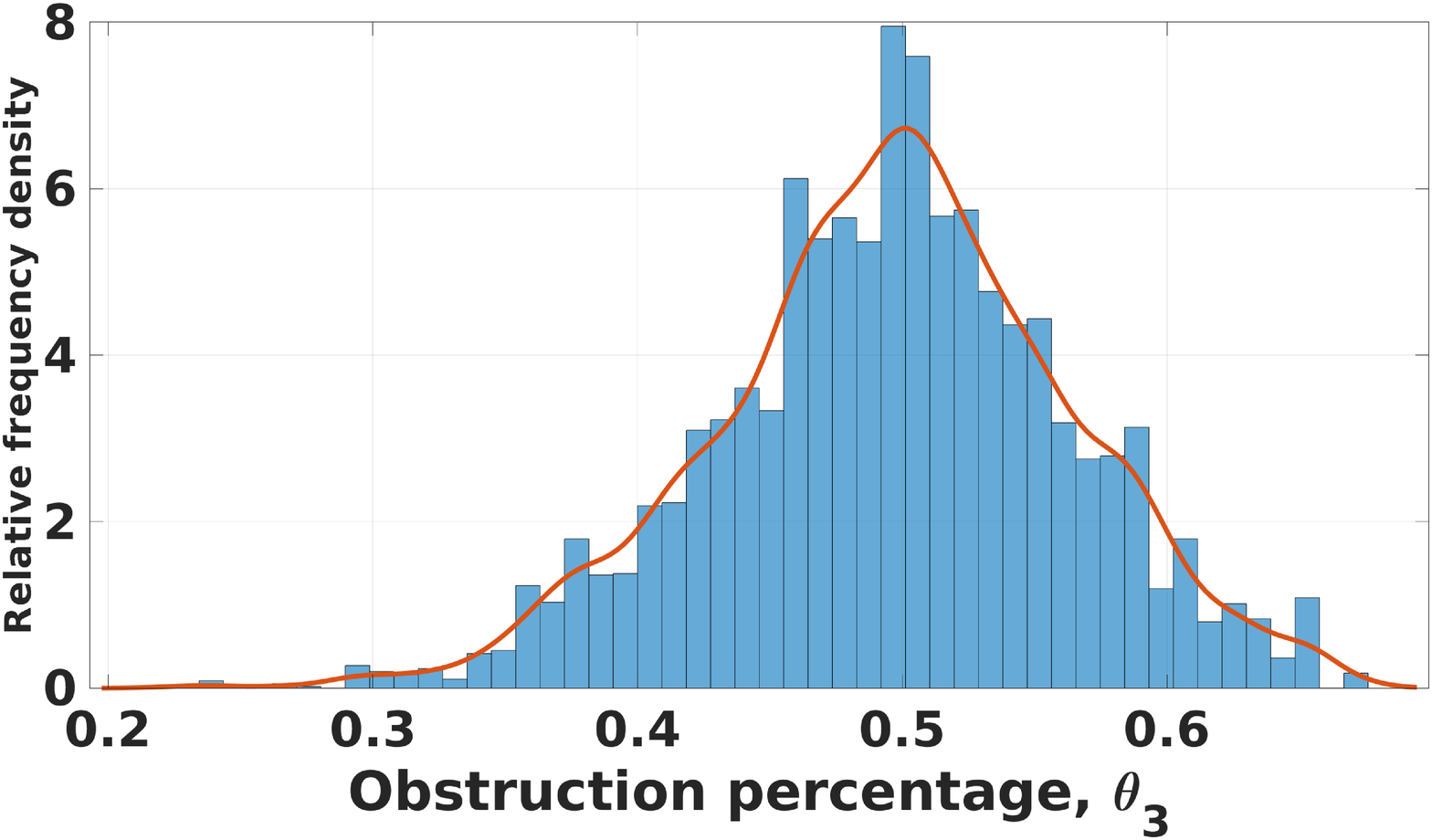}}
\subfloat[]{\centering{}\includegraphics[width=0.5\textwidth]{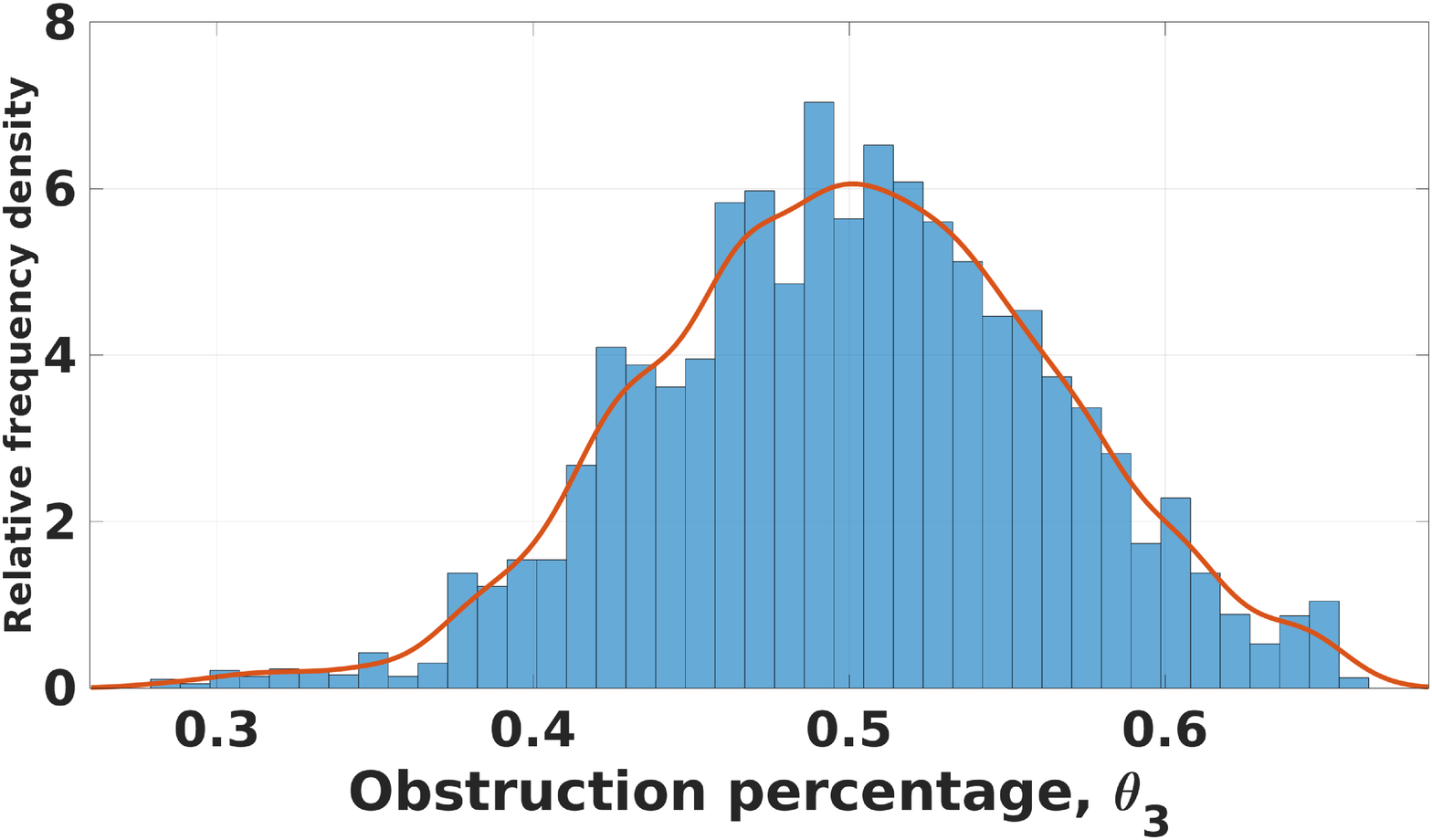}}

\subfloat[]{\centering{}\includegraphics[width=0.5\textwidth]{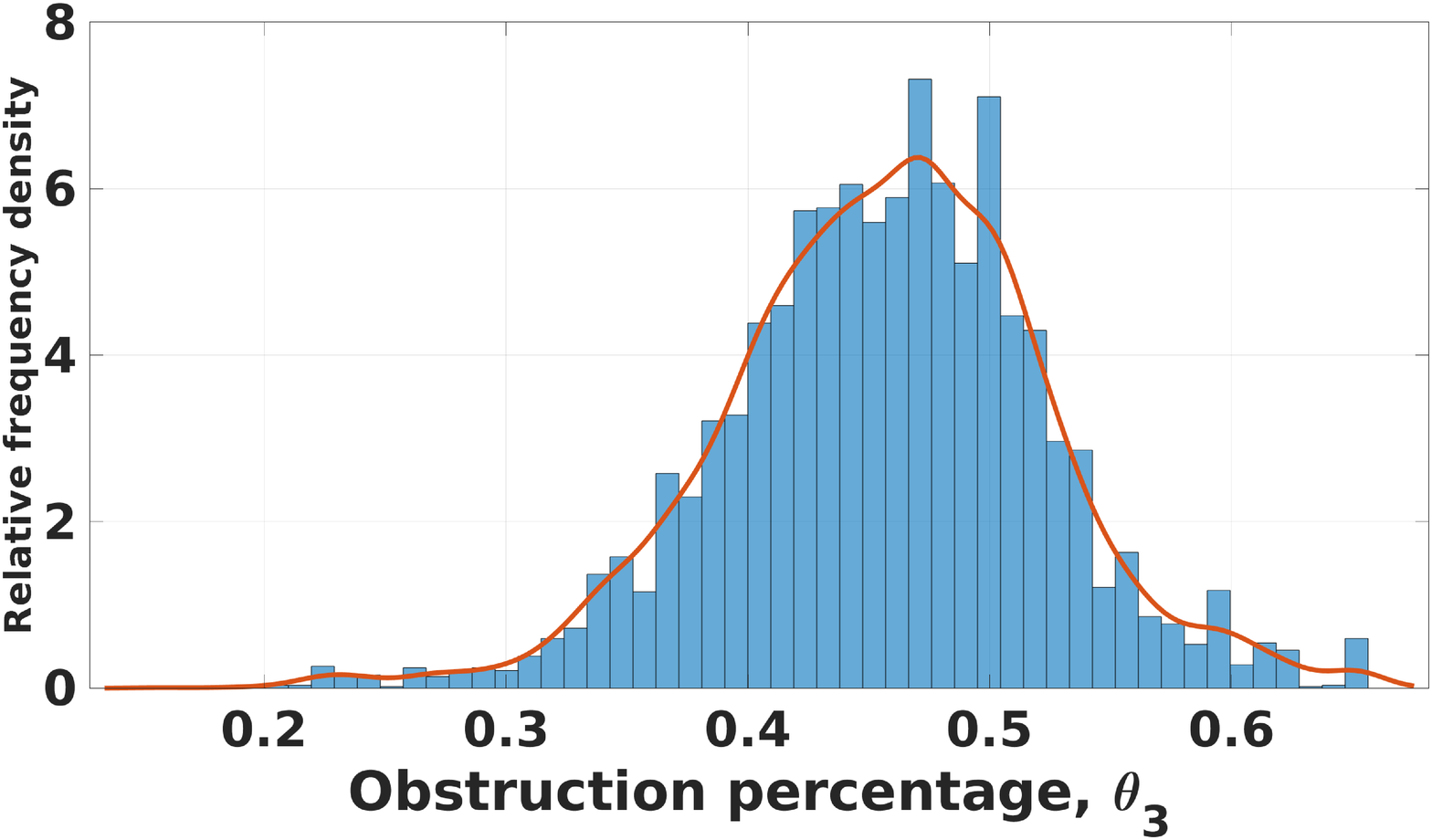}}
\subfloat[]{\centering{}\includegraphics[width=0.5\textwidth]{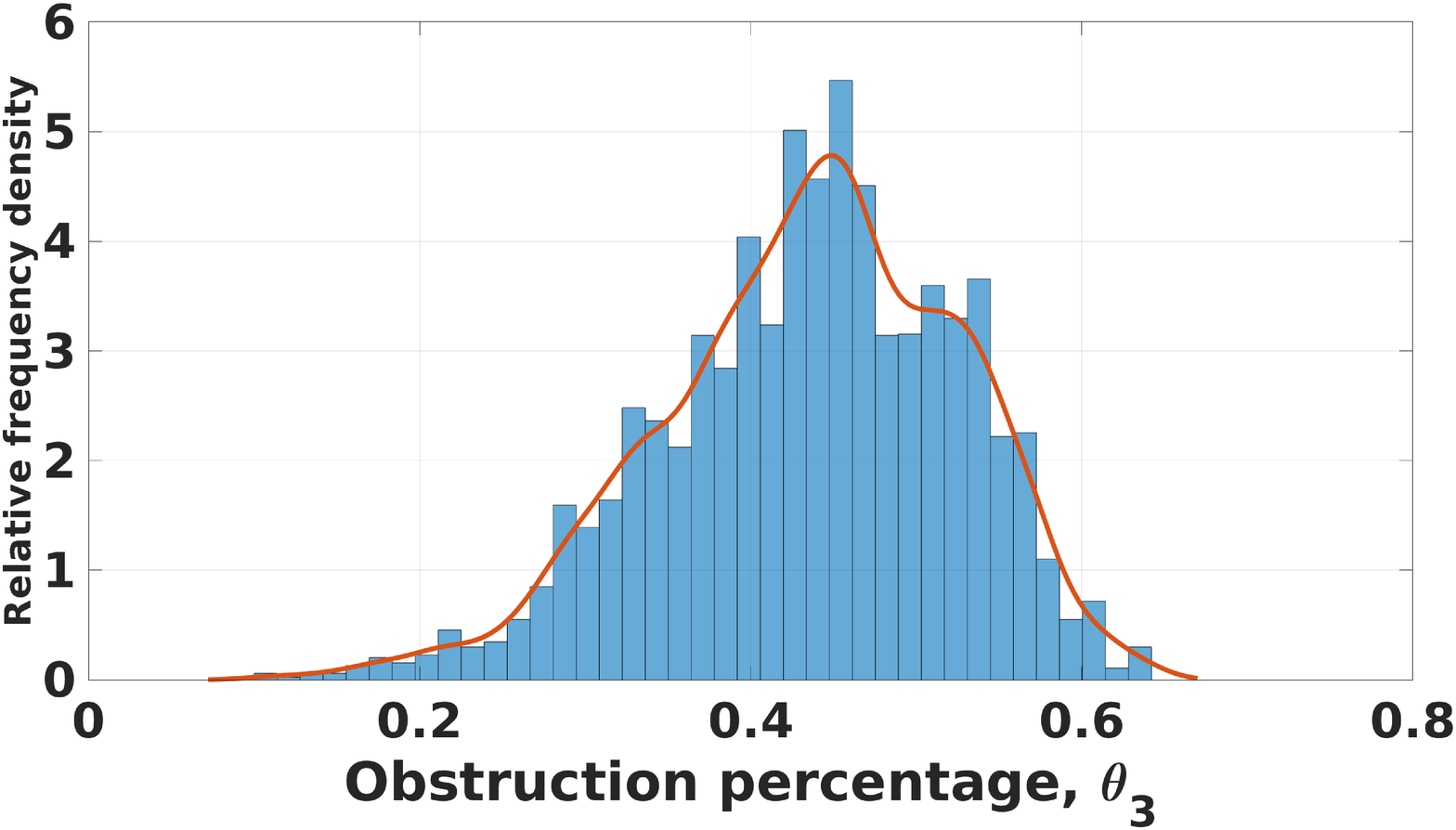}}

\subfloat[]{\centering{}\includegraphics[width=0.5\textwidth]{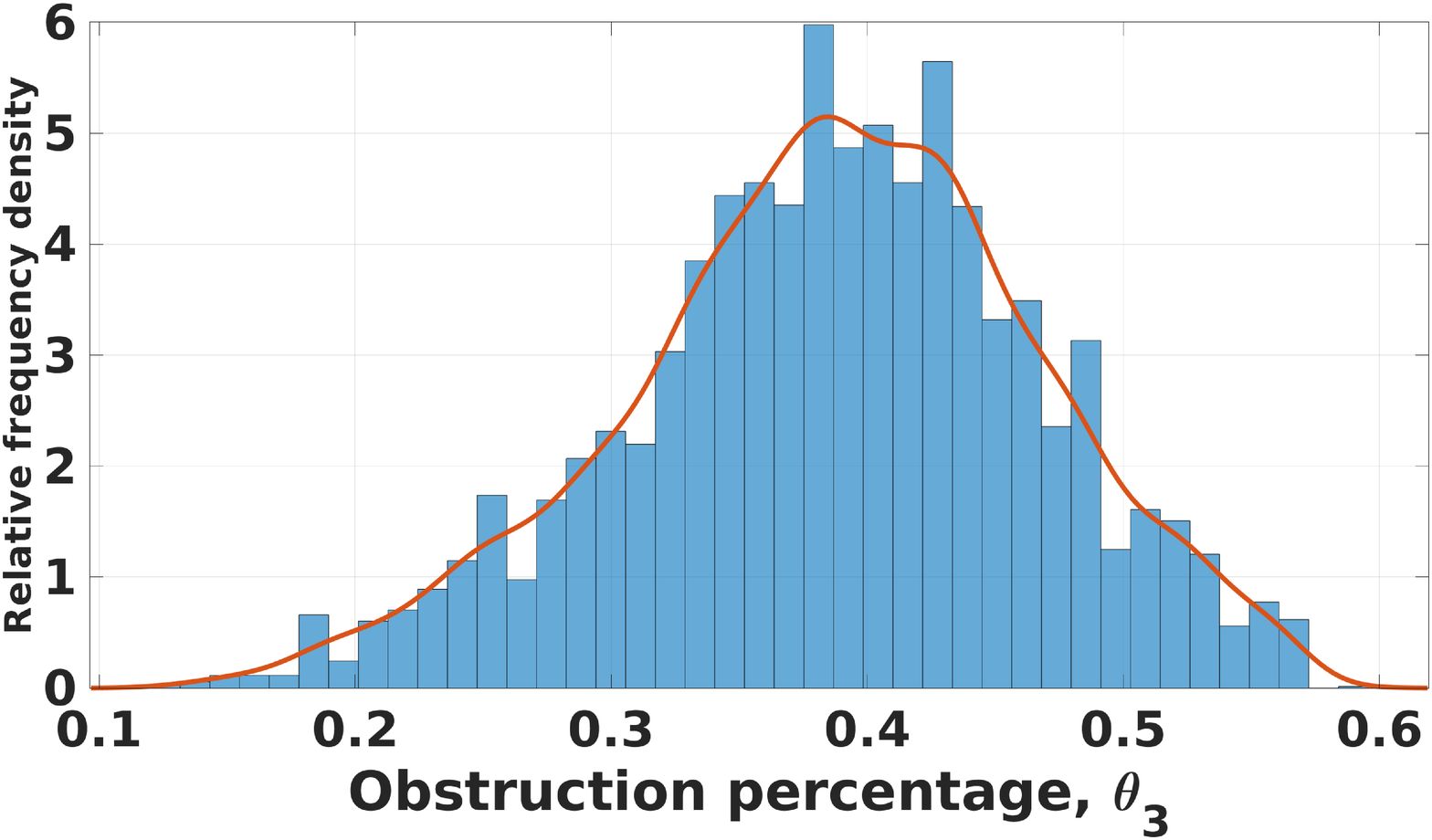}}
\subfloat[]{\centering{}\includegraphics[width=0.5\textwidth]{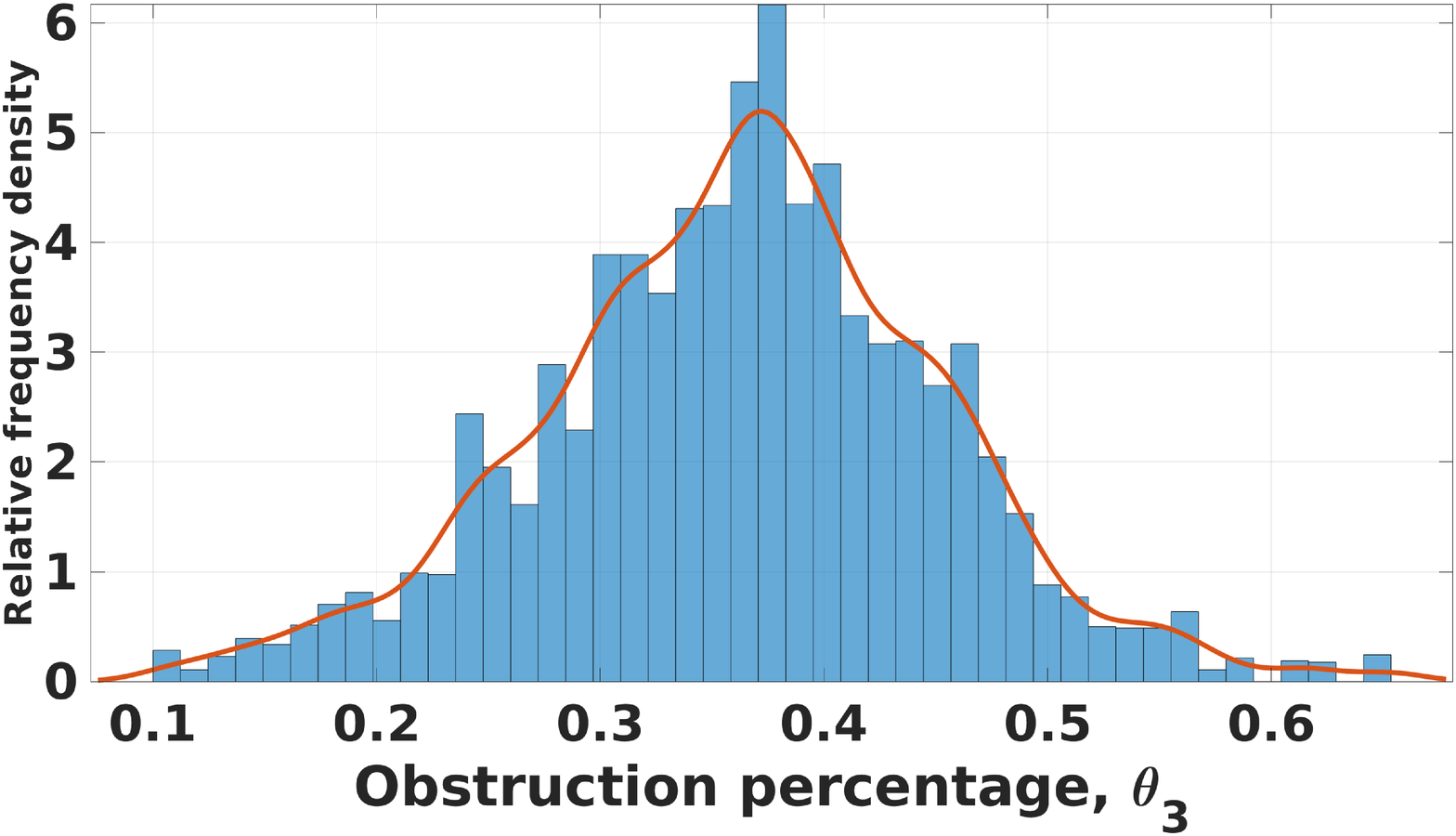}}
\par\end{centering}

\caption{Probability density function estimate of $\theta_3$ (Obstruction percentage) of experiments 1-8
(a-h)}
\label{fig:theta_3_dist}
\end{figure}

Figures \ref{fig:theta_1_dist}-\ref{fig:theta_3_dist} are the probability density functions of the parameters $\theta_{1},\theta_{2},\theta_{3}$ for each experiment, respectively.
As can be seen, most numerical distributions have Gaussian profiles with a mean center almost at the original parameters, with the exception of the distributions corresponding to obstacle size $\theta_{2}$ (Figure \ref{fig:theta_2_dist}).

For the obstacle position $\theta_{1}$ (Figure \ref{fig:theta_1_dist}) and the obstruction blockage percentage $\theta_{3}$ (Figure \ref{fig:theta_3_dist}), we notice that if the tangential velocity measurement site $\Gamma_m$ is above or on the right side of the obstacle, the standard deviation decreases, this could be explained by the fact that the obstacle information is transmitted in the direction of flow. However, even without the tangential velocity, we can obtain a nice Gaussian profile (Figures \ref{fig:theta_1_dist}(e), \ref{fig:theta_3_dist}(e)) with the down side of having a larger standard deviation.

For obstacle size $\theta_{2}$ (see Figure \ref{fig:theta_2_dist}), in most cases, the distributions look uniform. But as we can notice, the best results are those from experiments 2,5 and 6, in particular, the ones with the largest observable domain of the tangential velocity. Therefore, the tangential velocity seems to be a more important measurement in determining the correct size of the obstacle.

\begin{figure}[hbt!]
\begin{center}
\subfloat[Experiment 1]{\includegraphics[width=0.7\textwidth]{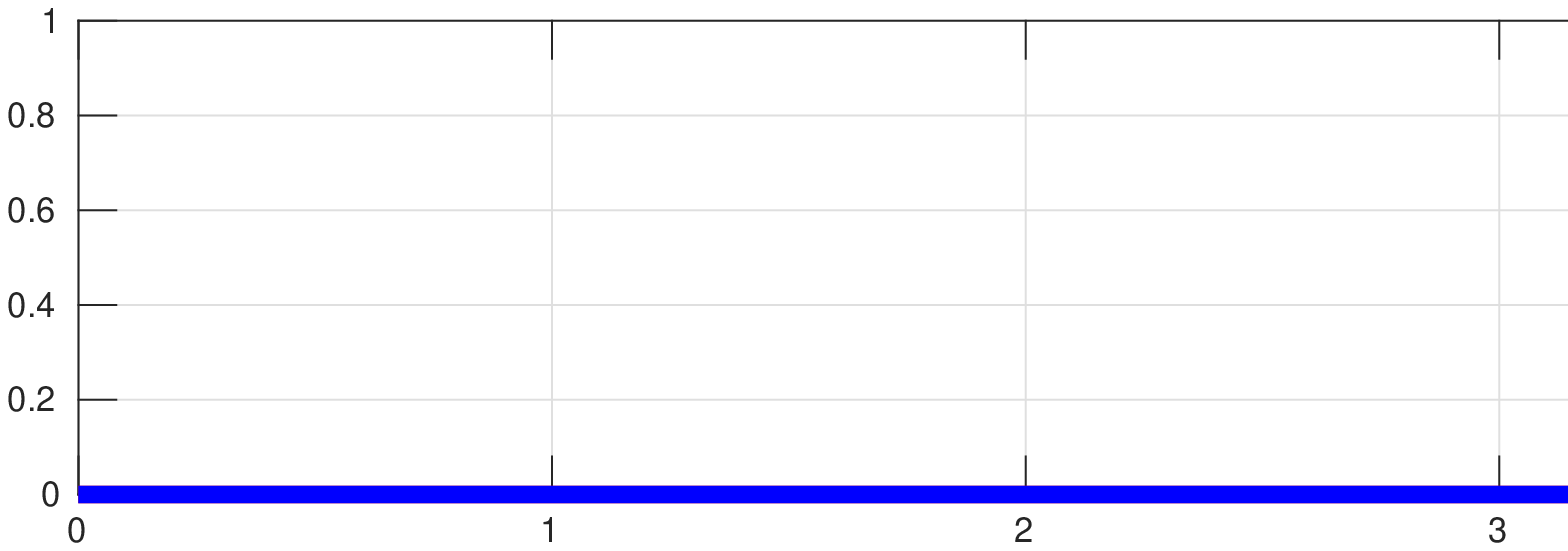}}

\subfloat[Experiment 2]{\includegraphics[width=0.7\textwidth]{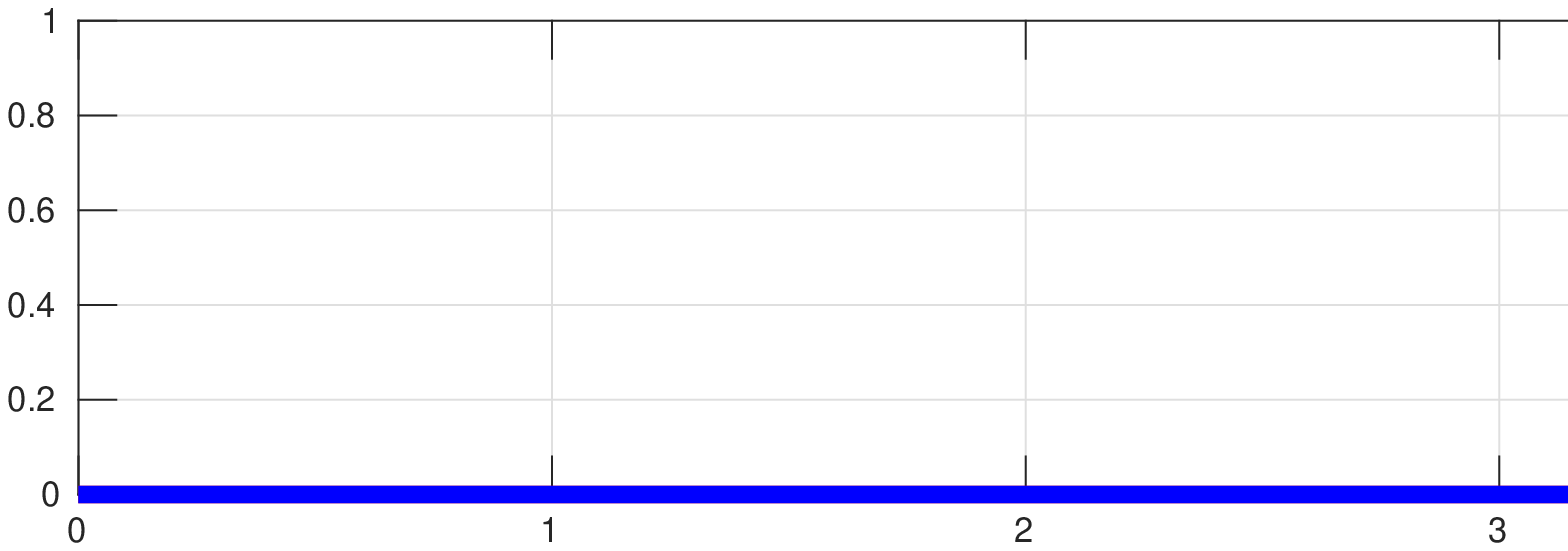}}

\subfloat[Experiment 3]{\includegraphics[width=0.7\textwidth]{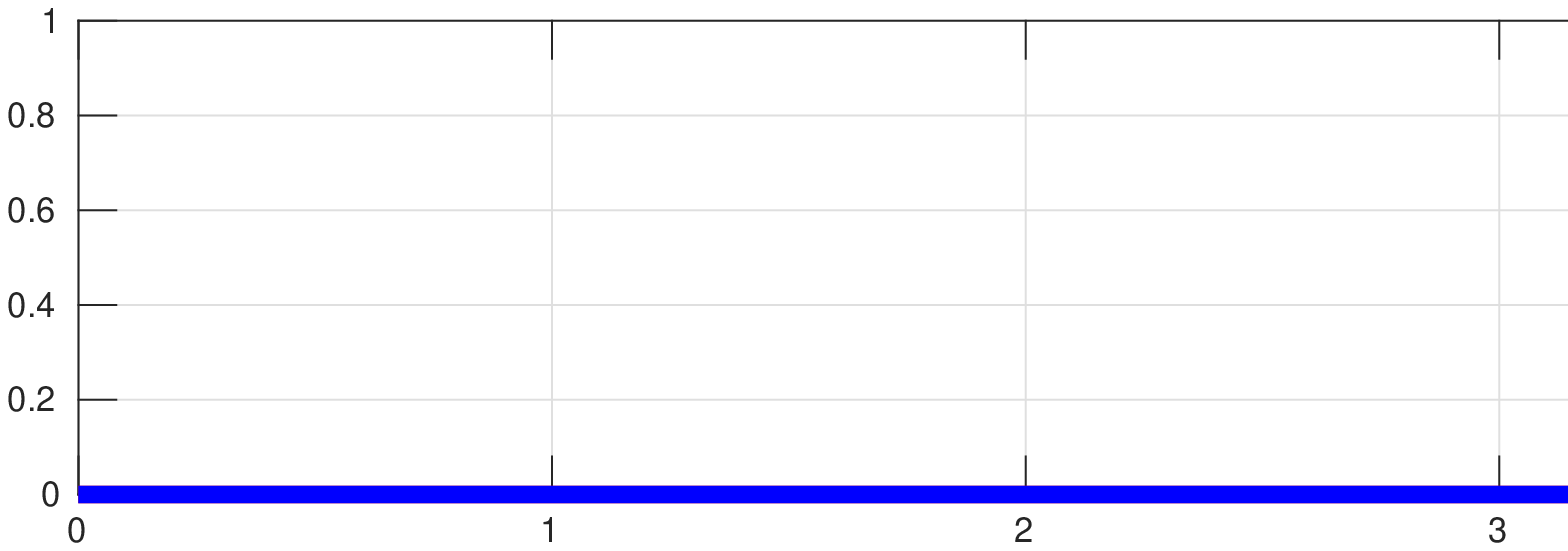}}

\subfloat[Experiment 4]{\includegraphics[width=0.7\textwidth]{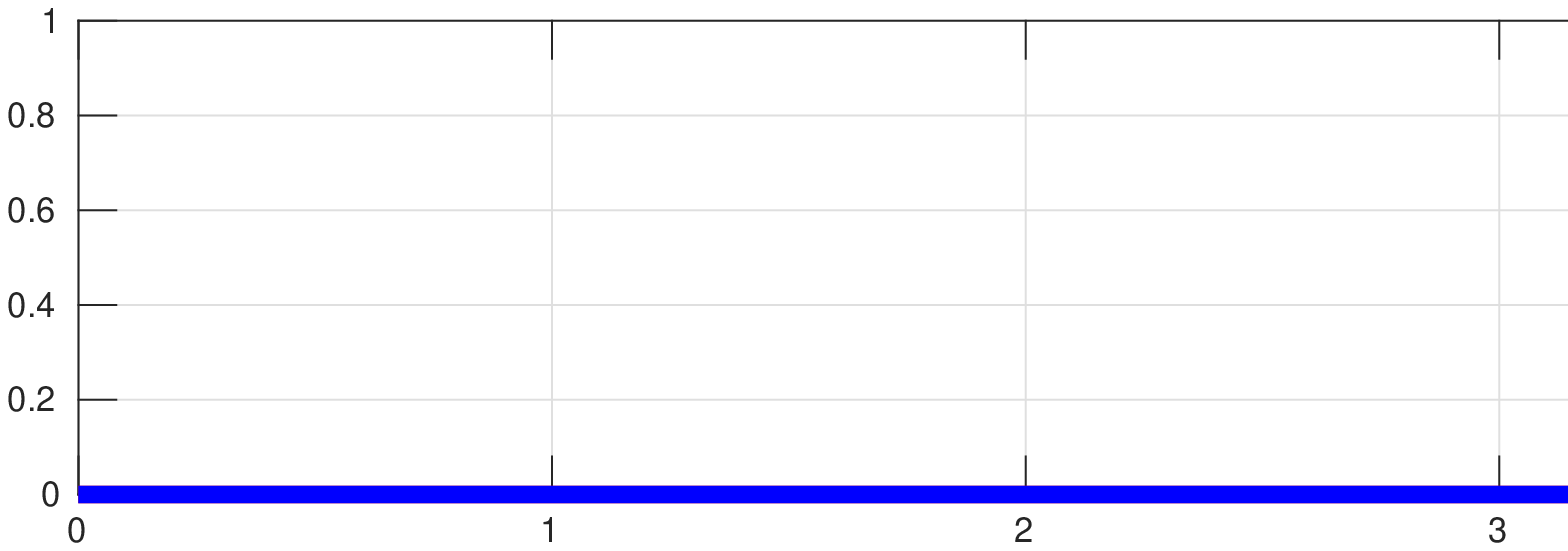}}

\subfloat[Experiment 5]{\includegraphics[width=0.7\textwidth]{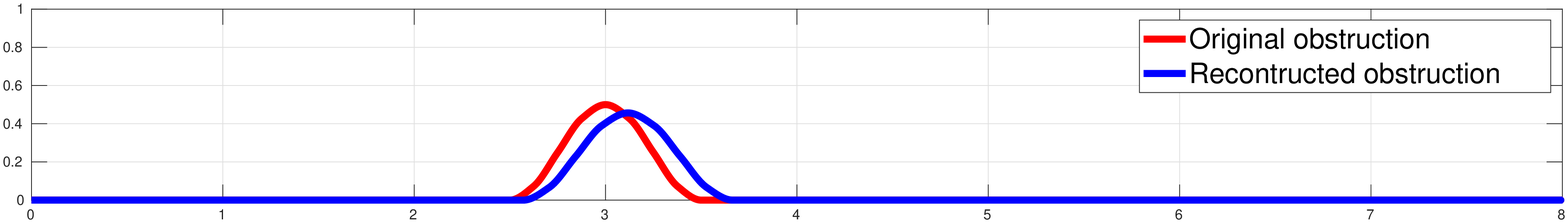}}

\subfloat[Experiment 6]{\includegraphics[width=0.7\textwidth]{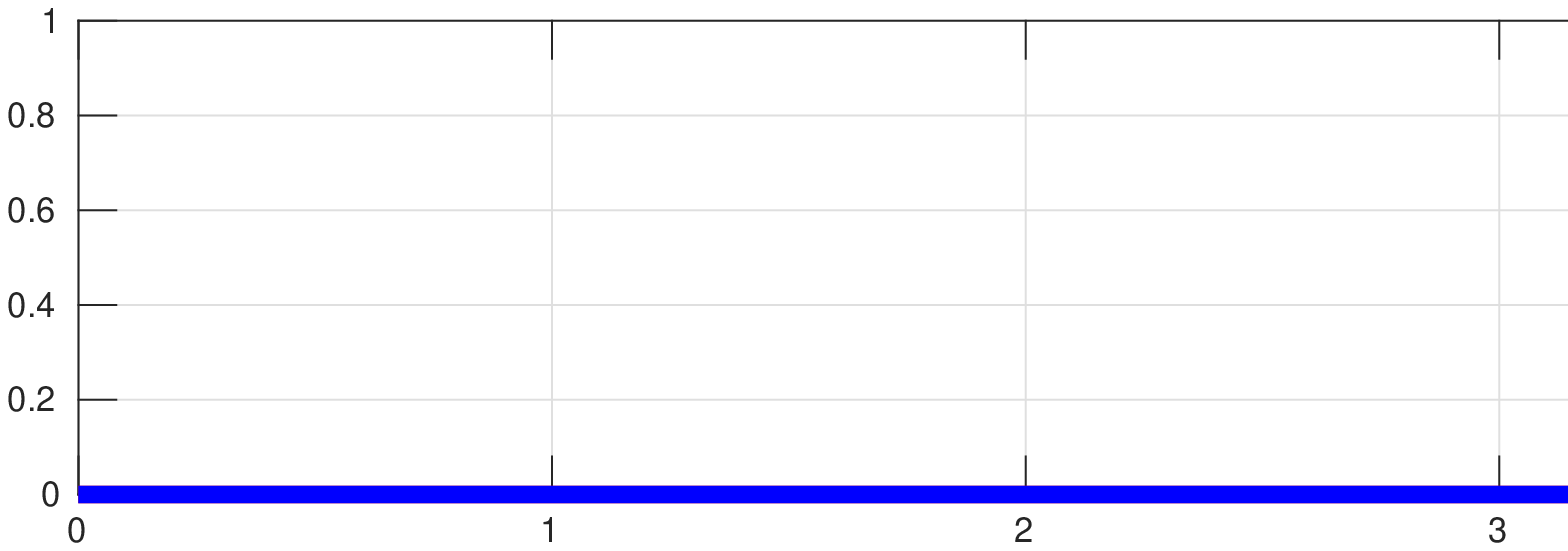}}

\subfloat[Experiment 7]{\includegraphics[width=0.7\textwidth]{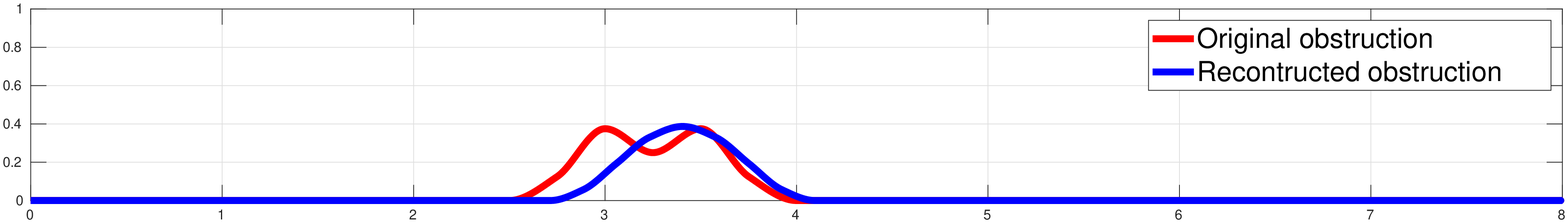}}

\subfloat[Experiment 8]{\includegraphics[width=0.7\textwidth]{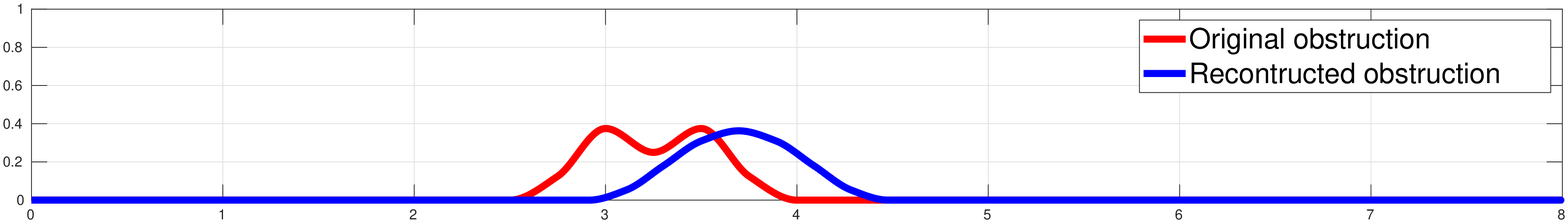}}
\end{center}
\caption{Original domain shape vs. reconstructed shape of experiments 1-8 (a-h)}
\label{fig:Comparacion-Reconstruc-onda-spline}
\end{figure}

In figure \ref{fig:Comparacion-Reconstruc-onda-spline}(a-h) we show a comparison between the original domain and the reconstructed domain for each of the eight numerical experiments. As said previously, the best results are obtained with maximum information about the tangential velocity vector (Figure \ref{fig:Comparacion-Reconstruc-onda-spline}(b)). In all cases, an acceptable depth and size of obstruction is obtained, which are the most relevant parameters to determine whether a stenosis is dangerous or not.

Notice also that figure \ref{fig:Comparacion-Reconstruc-onda-spline}(g-h) shows that even if the obstruction is outside the class of obstruction defined by $(\theta_1,\theta_2,\theta_3)$, we can still obtain some information with our low-parameter model.

It should be noted that the location of the wave measurement domain has no relevant impact on parameter estimation. Numerically, this is because in each case we recover all of the fluid tensor $\sigma(\bb{u},p)n$ on the upper wall of the fluid domain. The only variant is the time interval length of this tensor estimation, as stated in Subsection \ref{wave_subsection}.  


\begin{table}[hbt!]
\caption{Experiment 1-4 parameter summary and numerical result.}\label{Mcmc-table-1} 
{\small
\begin{center}
\begin{tabular}{|c|cc|cc|cc|cc|}
\hhline{~|--|--|--|--|}
\multicolumn{1}{c|}{} & \multicolumn{2}{c|}{Experiment 1} & \multicolumn{2}{c|}{Experiment 2} & \multicolumn{2}{c|}{Experiment 3} & \multicolumn{2}{c|}{Experiment 4}\tabularnewline
\hhline{-|--|--|--|--|}
$S_{m}$  & \multicolumn{2}{c|}{$[0,3]\times\{3\}$ } & \multicolumn{2}{c|}{$[0,8]\times\{3\}$ } & \multicolumn{2}{c|}{$[5,8]\times\{3\}$ } & \multicolumn{2}{c|}{$[0,8]\times\{3\}$}\tabularnewline
$\Gamma_{m}$  & \multicolumn{2}{c|}{$[0,2]\times\{1\}$ } & \multicolumn{2}{c|}{$[0,8]\times\{1\}$ } & \multicolumn{2}{c|}{$[6,8]\times\{1\}$ } & \multicolumn{2}{c|}{$\emptyset$}\tabularnewline
\hline 
Original  & Mean  & S.D  & Mean  & S.D  & Mean  & S.D  & Mean  & S.D\tabularnewline
\hline 
$\theta_{1}=4.0$ & $3.8929$ & $0.4871$ & $3.9652$ & $0.3207$ & $3.8983$ & $0.4698$ & $3.9191$ & $0.4722$ \tabularnewline
$\theta_{2}=1.0$ & $1.0854$ & $0.2653$ & $1.0879$ & $0.2567$ & $1.0667$ & $0.2834$ & $1.0599$ & $0.2731$\tabularnewline
$\theta_{3}=0.5$ & $0.4788$ & $0.0636$ & $0.4841$ & $0.0593$ & $0.4961$ & $0.0669$ & $0.5009$ & $0.0639$\tabularnewline
\hline 
\multicolumn{1}{|c|}{$\left\Vert \theta-\theta^{*}\right\Vert _{2}$ } & \multicolumn{2}{c|}{0.1386 } & \multicolumn{2}{c|}{0.0959 } & \multicolumn{2}{c|}{0.1217 } & \multicolumn{2}{c|}{0.1007 }\tabularnewline
\hline 
\end{tabular}
\end{center}
}
\end{table}

\begin{table}[hbt!]
\caption{Experiment 5-6 parameter summary and numerical result.}\label{Mcmc-table-2}
{\small
\begin{center}
\begin{tabular}{|ccc|ccc|}
\hline 
\multicolumn{6}{|c|}{$S_{m}=[0,8]\times\{3\}$ }\tabularnewline
\hline 
\multicolumn{6}{|c|}{$\Gamma_{m}=[0,8]\times\{1\}$ }\tabularnewline
\hline 
\multicolumn{3}{|c|}{Experiment 5} & \multicolumn{3}{c|}{Experiment 6}\tabularnewline
\hline 
Original  & Mean  & S.D  & Original  & Mean  & S.D \tabularnewline
\hline 
$\theta_{1}=2.5$ & $2.5692$ & $0.2872$ & $\theta_{1}=5.5$ & $5.3433$ & $0.3564$ \tabularnewline
$\theta_{2}=1.0$ & $1.0955$ & $0.2623$ & $\theta_{2}=1.0$ & $1.0433$ & $0.2871$ \tabularnewline
$\theta_{3}=0.5$ & $0.4566$ & $0.0667$ & $\theta_{3}=0.5$ & $0.4353$ & $0.0894$ \tabularnewline
\hline 
\multicolumn{3}{|c|}{$\left\Vert \theta-\theta^{*}\right\Vert _{2}=$ 0.1257} & \multicolumn{3}{c|}{$\left\Vert \theta-\theta^{*}\right\Vert _{2}=$0.1749}\tabularnewline
\hline 
\end{tabular}
\end{center}
}
\end{table}

\begin{table}[hbt!]
\caption{Experiment 7-8 parameter summary and numerical result.}\label{Mcmc-table-3}
{\small
\begin{center}
\begin{tabular}{|ccc|ccc|}
\hhline{---|---}
\multicolumn{3}{|c|}{Experiment 7} & \multicolumn{3}{c|}{Experiment 8}\tabularnewline 
\hhline{---|---}
\multicolumn{3}{|c|}{$S_{m}=[0,8]\times\{3\}$} & \multicolumn{3}{c|}{$S_{m}=[0,8]\times\{3\}$}\tabularnewline
\multicolumn{3}{|c|}{$\Gamma_{m}=[0,8]\times\{1\}$} & \multicolumn{3}{c|}{$\Gamma_{m}=\emptyset$}\tabularnewline
\hline
Parameter & Mean & S.D & Parameter & Mean & S.D\tabularnewline
\hline 
$\theta_{1}$ & $2.7139$ & $0.4047$ & $\theta_{1}$ & $2.9135$ & $0.6575$\tabularnewline 
$\theta_{2}$ & $1.3822$ & $0.5430$ & $\theta_{2}$ & $1.5646$ & $0.5698$\tabularnewline
$\theta_{3}$ & $0.3868$ & $0.0799$ & $\theta_{3}$ & $0.3629$ & $0.0886$\tabularnewline
\hline 
\multicolumn{3}{|c|}{$\left\Vert f-f(\theta^{*})\right\Vert _{2}=0.1251$} & \multicolumn{3}{c|}{$\left\Vert f-f(\theta^{*})\right\Vert _{L^{2}}=0.2958$}\tabularnewline
\hline 
\end{tabular}
\end{center}
}
\end{table}

\FloatBarrier

\section{Conclusions}
\noindent
In this article, we formulate and study a geometric inverse problem for the identification of a boundary obstruction immersed in a Stokes fluid contained in a 2D rectangular elastic duct. A novelties of this work are the mixed boundary condition and that we introduced an exterior approach, which is defined in the following sense. The Stokes flow, contained in an interior domain, becomes "turbulent" after hitting the boundary obstruction, producing a movement of the elastic boundary of the duct, this movement in turn generates acoustic waves.  The inverse problem, that is, the identification of the location, extension and height of the obstruction (i.e. lumen reduction), is then formulated and solved by using, far enough, external wave measurements relative to the duct, that is, the measurements are placed in what can be regarded as the exterior domain. It is worth pointing out that, to the best of our knowledge, all the former related works reported in the literature assume that the obstruction is contained inside the domain and that the boundary measurements are not external, that is, they are located at the unique boundary of the problem.

\vskip 0,3cm
\noindent
The problem presented in this work was studied and analyzed both from a continuous and numerical point of view. From the theoretical perspective, we proved the 
well--posedness of the Stokes flow with mixed boundary conditions, to our knowledge, this problem has only been treated in the literature with Dirichlet boundary conditions. In addition, we prove the uniqueness of the obstruction when the data are given by the normal component of the Cauchy stress tensor and the tangential velocity in a subset of the boundary. 
\vskip 0,3cm
\noindent
From the numerical point of view, we solved, with a fairly reasonable degree of accuracy, the identification of a boundary obstruction immersed in a Stokes flow by measuring the acoustic wave at the external boundary of the problem. The fluid obstruction inverse problem, which was solved using a Monte Carlo Markov chain (MCMC) method, involves the solution of two direct problems. A direct Stokes problem, solved using mesh-free hybrid radial basis function kernels, and a direct wave system solved by using the finite element method. 
\vskip 0,3cm
\noindent
Several related open problems, such as the regularity of the Stokes flow with these particular mixed boundary conditions or the well-posedness of this flow-acoustic inverse problem, remain to be proved.

Further work is in progress on these topics.

\section*{Acknowledgement}
The work of C.Montoya was supported by the Croatian Science Foundation under the project ConDyS IP-2016-06-2468 Dynamic systems management. The work of L.Breton was supported by the Grant 624497 from the National Council of Science and Technology of Mexico. This work was also supported by the National Autonomous University of Mexico [grant: PAPIIT, IN102116].

\begin{appendices}	
\section{}\label{Appendix.Green_identity}

In this appendix, we recall some useful result about the $H_{00}^{1/2}$ space and also we prove Theorem \ref{Th:Green_identity}. Until further notice, we consider a bounded open $\Omega\subset\mathbb{R}^{2}$ such that $\partial\Omega$ is a curvilinear polygone of class $C^{k,1}$
in the sense of \cite[Definition 1.4.5.1]{grisvard2011elliptic}.
Furthermore, we assume that there is $M\in\mathbb{N}$ such that:
\[
\partial\Omega=\bigcup_{j=1}^{2M}\overline{\Gamma_{j}}\quad,
\]
were each $\Gamma_{j}$ is a curve of class $C^{k,1}$, $\Gamma_{j+1}$
follows $\Gamma_{j}$ according to the positive orientation and we
have that $\Gamma_{j}\subset\partial\Omega$ , $\Gamma_{j}\bigcap\Gamma_{j+1}=\emptyset$
. We also define the following Lion-Maganese space:
\[
H_{00}^{1/2}(\Gamma_{i})=\left\{ u\in H^{1/2}(\Gamma_{i}):\ d(x,\partial\Gamma_{i})^{-1/2}u\in L^{2}(\Gamma_{i})\right\} \ ,
\]
equipped with the norm: $\normsq{u}{H_{00}^{1/2}(\Gamma_{i})}=\normsq{u}{{H^{1/2}}(\Gamma_{i})}+\normsq{d^{-1/2}u}{L^{2}(\Gamma_{i})}$. For a smooth function $u\in{\cal D}(\overline{\Omega})$ we denote
$\gamma_{j}(u),$the corresponding trace mappings (i.e):
\begin{equation}
\gamma_{j}(u):u\rightarrow u\lvert_{\Gamma_{j}}\quad.\label{eq:Trace restriction}
\end{equation}
\noindent
An important question is when does the trace of a function in $u\in H^{1/2}(\Omega)$ belong to $H_{00}^{1/2}(\Gamma_{j})$. In fact, from \cite[Theorem 1.5.2.3]{grisvard2011elliptic} we have the following Theorem.
\begin{thm}\label{th.trace_h_00} 
Let $\Omega$ be a bounded open subset of $\mathbb{R}^{2}$ whose
boundary is a curvilinear polygon of class at least $C^{1}$. Let
$u\in H^{1}(\Omega)$ such that:
\[
\gamma_{j-1}(u)=\gamma_{j+1}(u)=0\quad,
\]
 then the mapping define by :
\[
u\to\gamma_{j}(u)\quad,
\]
is a linear continuous mapping from $H^{1}(\Omega)\rightarrow H_{00}^{1/2}(\Gamma_{j})$. 
\end{thm}

Also from within the prove of \cite[Theorem 1.5.2.3]{grisvard2011elliptic} we have the following result.
\begin{thm}\label{th.extension}
Let $\Omega$ be a bounded open subset of $\mathbb{R}^{2}$ whose
boundary is a curvilinear polygon of class at least $C^{1}$. Let
$f_{j}\in H_{00}^{1/2}(\Gamma_{j})$, and define $\widetilde{f}$
as the extension by zero en $\Gamma$ (i.e):
\[
\widetilde{f}(x)=\begin{cases}
f_{j}(x) & x\in\Gamma_{j}\\
0 & x\notin\Gamma_{j}
\end{cases}\quad,
\]
then $\widetilde{f}\in H^{1/2}(\Gamma)$ and there exist a constant
$C>0$ such that:
\[
\norm{\widetilde{f}}{H^{1/2}(\Gamma)}\leq C\norm{f_{j}}{H_{00}^{1/2}(\Gamma_{j})}\quad.
\]
\end{thm}
Using  the same argument as in \cite[Theorem 1.5.3.10]{grisvard2011elliptic} and an analogous prove of \cite[Lemma 2.4]{amrouche2014lp} we have the following Theorem:
\begin{thm}\label{th.dual_tensor}
Let $\Omega$ be a bounded open subset of $\mathbb{R}^{2}$ whose
boundary is a curvilinear polygon of class at least $C^{1}$. Let
$\bb E(\Omega)=\{\bb u\in\bb H^{1}(\Omega):\ \Delta\bb u\in\bb H_{0}(div,\Omega)'\}$,
then the mapping:
\[
\bb u\rightarrow2[\bb D(\bb u)\bb n]_{tg}\lvert_{\Gamma_{j}}\quad,
\]
which is defined on ${\cal D}(\overline{\Omega})$ has a unique continuous
extension as an operator from:
\[
\bb E(\Omega)\rightarrow\left(\bb H_{00}^{1/2}(\Gamma_{j})\right)'
\]
\end{thm}
\noindent
We are now ready to prove our green identity 
\begin{proof}[Proof of theorem \ref{Th:Green_identity}]
Notice that for every $\bb u\in\bb{{\cal \bb D}}(\overline{\Omega})$
and $\varphi\in\{\bb{u\in}\bb H^{1}(\Omega):div(\bb u)=0\}$ we have
the following Green's indentity:
\begin{equation}\label{eq.green_idendityu_prove}
\left\langle \Delta\bb u,\bb{\varphi}\right\rangle _{\Omega}=\int_{\Omega}2\bb D(\bb u):\bb D(\bb{\varphi})\,dx-\sum_{j=1}^{2M}{\left\langle 2[\bb D(\bb u)n]_{tg},\gamma_{j}(\bb u)\right\rangle _{\Gamma_{j}}}
\end{equation}
Thus, assuming that $\gamma_{j}(\bb u)\in\prod_{j=1}^{2M}H_{00}^{1/2}(\Gamma_{j})$
and using the Korn inequality, we can see that all terms of \eqref{eq.green_idendityu_prove}
are continuous in $\bb u$ for the norm of $\bb E(\Omega)$. Therefore,
the result follows from the density of $\bb{{\cal \bb D}}(\overline{\Omega})$
in $\bb E(\Omega)$.
\end{proof}

\end{appendices}

\bibliography{Bibliografia}

\end{document}